\documentclass[stsy,nonblindrev]{informsnosub}

\OneAndAHalfSpacedXI
\usepackage{longtable}
\usepackage{endnotes}
\usepackage{color,soul}
\usepackage{arydshln}
\let\footnote=\endnote

\usepackage{natbib}
\usepackage{amsmath,,amssymb}
\usepackage{scalerel}
\usepackage{xcolor}
\usepackage{booktabs}
\usepackage{enumitem}
\usepackage{tcolorbox}
\usepackage{mathtools}

 \NatBibNumeric
 \def\bibfont{\small}%
 \bibpunct[, ]{[}{]}{,}{n}{}{,}%

\let\footnote=\endnote

\newtheorem{theorem}{Theorem}

\newtheorem{corollary}{Corollary}

\newtheorem{lemma}{Lemma}

\newcommand{\FCBA}{\mathbb A}
\newcommand{\FCBG}{\mathbb G}

\newcommand{\FCBE}{\mathbb E}
\newcommand{\FCBH}{\mathbb H}

\newcommand{\hD}{\hat{D}^{(L)}}

\newcommand{\irange}{1 \leq i \leq m}
\newcommand{\jrange}{1 \leq j \leq n}

\newcommand{\ipp}{\mathbf{IP}^{k-}}
\newcommand{\ip}{\mathbf{IP}}

\newcommand{\bal}{\boldsymbol \alpha}
\newcommand{\bnu}{\boldsymbol \nu}
\newcommand{\bk}{\boldsymbol \kappa}
\newcommand{\ipos}{IP}

\newcommand{\E}{\mathbf E}
\newcommand{\jD}{\mathbf d}
\newcommand{\BZ}{\mathbf {Z}}
\newcommand{\B}{\mathbf {B}}

\newcommand{\y}{\mathbf {y}}

\newcommand{\BD}{\mathbf {D}}
\newcommand{\z}{\mathbf {z}}

\newcommand{\Q}{\mathbf {Q}}
\newcommand{\x}{\mathbf {x}}
\newcommand{\aj}{\mathbf A_j}

\newcommand{\mA}{\mathcal A}
\newcommand{\mD}{\mathcal{D}}

\newcommand{\mG}{\mathcal G}

\newcommand{\mK}{\mathcal K}

\newcommand{\mR}{\mathcal{R}}
\newcommand{\mT}{\mathcal T}

\newcommand{\mZ}{\mathcal{Z}}

\newcommand{\sL}{\hat{s}^{(L)}}

\newcommand{\tL}{t^{(L)}}

\newcommand{\tIP}{{\mathbb{IP}}}

\newcommand{\sD}{\mathcal D^{(L)}}
\newcommand{\tsD}{\tilde{\mathcal D}^{(L)}}
\newcommand{\sBD}{\mathbf D^{(L)}}
\newcommand{\hBD}{\hat{\mathbf D}^{(L)}}
\newcommand{\sjD}{\mathbf d^{(L)}}

\newcommand{\hjD}{\hat{\mathbf d}^{(L)}}
\newcommand{\hejD}{\hat{d}^{(L)}}
\newcommand{\sT}{\mathcal T^{(L)}}
\newcommand{\sE}{\mathcal E^{(L)}}
\newcommand{\sWp}{\mathcal W^{(L)-}}
\newcommand{\sW}{\mathcal W^{(L)}}
\newcommand{\hT}{\hat{\mathcal T}^{(L)}}
\newcommand{\hW}{\hat{\mathcal W}^{(L)}}
\newcommand{\hE}{\hat{\mathcal E}^{(L)}}
\newcommand{\stL}{\tilde{t}^{(L)}}
\newcommand{\stLp}{\tilde{t}^{(L)-}}
\newcommand{\hG}{\hat{\mathcal G}^{(L)}}
\newcommand{\sIP}{IP^{(L)}}                    
\newcommand{\hIP}{\widehat{IP}^{(L)}}

\newcommand{\thIP}{\widehat{{\mathbb{IP}}}^{(L)}}
\newcommand{\sbY}{\mathbb Y^{(L)}}
\newcommand{\hbY}{\hat{\mathbb Y}^{(L)}}

\newcommand{\bB}{\mathbb B}
\newcommand{\sB}{B^{(L)}}
\newcommand{\hB}{\hat{B}^{(L)}}
\newcommand{\sbB}{\mathbb B^{(L)}}
\newcommand{\hbB}{\hat{\mathbb B}^{(L)}}
\newcommand{\sbbB}{\mathbb B^{*(L)}}
\newcommand{\hbsB}{\hat{\mathbb B}^{*(L)}}
\newcommand{\nsB}{\underline{B}^{(L)}}
\newcommand{\hnsB}{\underline{\hat{B}}^{(L)}}
\newcommand{\mincost}{\mathcal C_{\min}^{(L)}}

\newcommand{\sttL}{\hat{t}^{(L)}_0}
\newcommand{\LL}{L^{(L)}}
\newcommand{\hLL}{\hat{L}^{(L)}}

\begin{document}
 \RUNAUTHOR{Reiman, Wan, and Wang}
 \RUNTITLE{Independent Base-Stock Policies}

 \TITLE{Asymptotically Optimal Inventory Control for Assemble-to-Order Systems}

 \ARTICLEAUTHORS{

 \AUTHOR{Martin I. Reiman}
 \AFF{Industrial Engineering and Operations Research, Columbia University, New York, NY 10027, \EMAIL{mr3598@columbia.edu}} %, \URL{}}
 \AUTHOR{Haohua Wan}
 \AFF{Industrial and Enterprise Systems Engineering, University of Illinois at Urbana-Champaign, Urbana, IL 61801, \EMAIL{hwan3@illinois.edu}}
 \AUTHOR{Qiong Wang}
 \AFF{Industrial and Enterprise Systems Engineering, University of Illinois at Urbana-Champaign, Urbana, IL 61801, \EMAIL{qwang04@illinois.edu}}
 }
 
\ABSTRACT{We consider Assemble-to-Order (ATO) inventory systems with a general Bill of Materials and general deterministic lead times. Unsatisfied demands are always backlogged. We apply a four-step asymptotic framework to develop inventory policies for minimizing the long-run average expected total inventory cost.  Our approach features a multi-stage Stochastic Program (SP) to establish a lower bound on the inventory cost and determine parameter values for inventory control.  Our replenishment policy deviates from the conventional constant base stock policies to accommodate non-identical lead times. Our component allocation policy differentiates demands based on backlog costs, Bill of Materials, and component availabilities. We prove that our policy is asymptotically optimal on the diffusion scale, that is, as the longest lead time grows, the percentage difference between the average cost under our policy and its lower bound converges to zero. In developing these results, we formulate a broad Stochastic Tracking Model and prove general convergence results from which the asymptotic optimality of our policy follows as specialized corollaries.  
}

\maketitle

\section{Introduction}
\label{sec:intro}

Optimal control of Assemble-to-Order (ATO) inventory systems is a canonical problem in inventory theory.
In an ATO system, product assembly is assumed to take a negligible amount of time, so there is no need to keep inventories of final products. Component supplies are not capacity-constrained, but it is necessary to hold inventories for them to accommodate replenishment lead times (i.e., delays between ordering and receiving components). 
The goal of ATO inventory management is to minimize the total inventory cost,
which consists of both holding and backlog costs,
by controlling the timing and quantities of component orders and allocation of available components to different product demands.  
 
While optimal control of single-product ATO systems has long been settled (see \cite{Karlin1958} and \cite{Rosling1989} for systems with single and multiple components respectively), optimizing multi-product ATO systems is an immensely more difficult problem that remains unsolved.  The complexity arises from the need to allocate components that are  used by multiple products. Optimal allocation depends on component availabilities, which are outcomes of replenishment decisions. Optimal replenishment needs to take into account how components will be allocated. A joint optimal solution of both decisions is contingent upon not only the current inventory and backlog levels, but also the times of all outstanding replenishment orders, giving rise to an enormous state space.   The problem becomes even more complicated in systems 
where components do not have
identical lead times, because the ordering of  components often needs to be coordinated with the availabilities of other components with longer lead times. 

There has been a large body of studies on managing multi-product ATO systems. One may find a thorough review of the related literature in \cite{Song2003}, and a more recent one in \cite{Atan2017}.  Many previous studies focus on particular types of policies, such as base stock replenishment policies, FIFO or No-Hold-Back allocation policies (see \cite{Lu2005}, \cite{Lu2010}, \cite{Huang2015} for some samples), and for periodic-review systems, allocation policies that always satisfy demands from previous periods first (e.g., \cite{Agrawal2001}, \cite{Akccay2004}, \cite{Hausman1998}, \cite{Zhang1997}).  Restricting the consideration to these types of policies makes the problem more tractable, but also sacrifices optimality, since both analytical and numerical studies (e.g., \cite{Dogru2010},\cite{Dogru2017}) have shown that there are often better alternatives in other types of policies. 

Asymptotic analysis has recently emerged as a powerful tool for analyzing inventory systems. Though a weaker standard than being exactly optimal, asymptotic optimality provides vital guidance to formulating novel inventory policies and evaluating their performance when meeting the latter criterion is analytically intractable. Asymptotic optimality has been used to justify the use of base stock policies in lost sale inventory systems (in the regime of high penalty costs) \cite{Huh2009}, ATO inventory systems with identical lead times  \cite{Reiman2015} or when differences in lead times  are small relative to the lead times themselves 
\cite{Reiman2016}, or systems with non-stationary demands and probabilistic service level constraints \cite{Wei2018}. It also supports the use of No-Hold-Back allocation policies in assemble-to-order N and W- systems \cite{Lu2015}. Asymptotic analysis has  led to several surprising discoveries. For instance, simple constant ordering policies can be highly effective in managing lost sale inventory systems with long lead times (\cite{Reiman2004}, \cite{Goldberg2016}, \cite{Xin2016}, \cite{Bu2017}); and randomness of lead times is a useful feature that can be exploited to reduce the inventory cost in backlog systems when orders can cross in time \cite{Stolyar2018}. 

In this paper, we study ATO inventory systems with \emph{non-identical, deterministic lead times and a general Bill of Materials (BOM)}. We develop an inventory policy that includes both replenishment and allocation decisions. We prove that our policy is asymptotically optimal, i.e., as the longest lead time increases, the percentage difference between the long-run average inventory cost   under our policy and the optimal policy converges to zero. 

%\textcolor{red}{
Our analysis follows a general four-step asymptotic framework that has produced many ground-breaking results in the study of stochastic processing networks (see, e.g., \cite{Harrison1988}, \cite{Harrison1990}, \cite{Harrison1996}, \cite{Harrison1999}). It has also been adopted for optimizing pricing, capacity, and allocation decisions in  ATO production-inventory systems \cite{Plambeck2006}. In \cite{Reiman2015}, the  framework is formally summarized where each step is explicitly described and specialized to the study of ATO inventory systems with identical lead times. As a continuation of the latter work, here we apply the framework to address ATO inventory systems with general deterministic lead times. Below is a step-by-step discussion of previous results and our new contributions.

 \begin{itemize}
 \item[]
 \emph{Step 1: relax some feasibility constraints of the original system to formulate a proxy model that is easier to solve and provides a lower bound on the cost achievable under any feasible policy}.

 In \cite{Reiman2012}, a multi-stage stochastic program (SP) is formulated as a static analog to dynamic ATO systems. The optimal objective value of the SP is proven to be a lower bound on the average inventory cost of the latter systems under any feasible policy. 
Unfortunately this lower bound takes the form of an infimum that is often not attained at finite values of the decision variables.

As a contribution of this paper, we transform this SP into an equivalent form where the optimal objective value is the same,
while the optimal decision variables are finite. This paves the way for using these optimal decision variables as parameters for ATO inventory control policies. 

\item[]
\emph{Step 2: solve the proxy model}. 

There have been many studies on the model structure and computational efficiency of the aforementioned SP, mainly on special cases that have two stages and correspond to systems with particular BOMs (\cite{DeValve2018},\cite{Dogru2017},\cite{Nadar2014},\cite{Jaarsveld2015},\cite{Zipkin2016}).
Developing efficient algorithms for multi-stage SPs is an active area of research.
Although we solve the SP for some simple examples in Section \ref{sec:simulation},
solving more general cases is well beyond the scope of this paper, and we focus
on a different aspect of the problem that is critical for proving asymptotic optimality of our policy. 
To use optimal solutions of the SP as parameters of inventory control policies, we need to show that these values are stable, i.e., they do not change drastically in response to small fluctuations of inputs. The analysis should be applicable to the SP with a general number of stages and  general BOMs. 

To this end, we  show that with a negligible loss of the accuracy, the SP can be approximated by a finite-dimension Linear Program (LP), which has a unique optimal solution that is Lipschitz continuous in input parameters.  By characterizing and making use of the structure of constraints of the approximating LPs, we prove that their optimal solutions are Lipschitz continuous in inputs that change with the state of the ATO system, and importantly, the Lipschitz constant is independent of the problem size. 

\item[]
 \emph{Step 3: use the optimal solution of the proxy model to formulate inventory control policies for the original systems}. 

In \cite{Reiman2015},  a two-stage SP  is used to set asymptotically optimal inventory policies for systems with identical lead times.  Replenishment decisions follow a base stock policy with the base stock levels specified by the first stage SP solution. Allocation decisions are controlled by an Allocation Principle, which uses the second stage SP solution to set backlog targets to dynamically choose the amount of each demand to serve. However, it is well known that base stock policies are inefficient for general ATO systems considered in this paper, which allow significantly different lead times \cite{Zipkin2000}. 

To the best of our knowledge, except for single-product ATO systems \cite{Rosling1989}, an asymptotically optimal replenishment policy remains to be developed for systems with general deterministic lead times and general BOMs.  We fill the gap in this paper by formulating such a policy, which uses the optimal solutions of the aforementioned multi-stage SP to set dynamic targets for inventory positions that determine order quantities. The policy generalizes previous work: it specializes to the policy in \cite{Rosling1989} in systems with a single product and the base stock policy in \cite{Reiman2015} in systems with identical lead times. We adopt the same Allocation Principle in \cite{Reiman2015} with minor twists to fit with the new replenishment policy.

\item[]
\emph{Step 4: show the performance benefit of the policy by proving that it is asymptotically optimal}.

For the special case of identical lead times, an SP-based policy has been shown to be asymptotically optimal \cite{Reiman2015}. However, the proof relies on the fact that the constant inventory positions prescribed by the two-stage SP can be exactly followed by a base stock policy. Also, the state of the system under this policy is completely determined by the history within the previous lead time. In contrast, for systems with non-identical lead times, the ideal inventory positions prescribed by the multi-stage SP may not always be attainable. The state of system can be affected by the history in the unbounded past. Therefore, new techniques are needed to prove asymptotic optimality of our policy for  general systems. 

We introduce a broadly-defined \emph{Stochastic Tracking Model}, which features a general target process and a general tracking process.   We prove that the expected difference between these two processes converges to zero. We apply this model to compare the inventory position and backlog targets prescribed by the SP for reaching the cost lower bound with their actual levels under our policy. We show that these comparisons are  special cases of the convergence results of this tracking model, and use this to establish asymptotic optimality of our policy. We also perform simulations on several special ATO systems to illustrate that our policy performs well, even under ``non-asymptotic'' conditions. 

\end{itemize}

The rest of the paper is organized as follows. We define the problem in Section \ref{sec:formulation}.
In Section \ref{sec:SPdisc}, we formulate the SP, develop bounds on its optimal solutions (Theorem \ref{thm:solbound}), and show that its optimal objective value is the same as that of the SP in \cite{Reiman2012}, and thus is a lower bound on the average inventory cost of the ATO system (Theorem \ref{thm:sameobj}).
In Section \ref{sec:policy}, we present our inventory policy. In Section \ref{sec:asymopt}, we show that the policy is asymptotically optimal if the resulting inventory positions and backlog levels converge to their respective SP-based targets (Theorem \ref{thm:final}). 
The latter convergence results are proven  in Section \ref{sec:sufficond} with the formulation and analysis of a general Stochastic Tracking Model. There,  Theorem \ref{thm:master} proves convergence for the general Stochastic Tracking model. Its corollaries show that inventory positions and backlog levels  converge to their targets if the latter satisfy stability conditions that require them to be asymptotically Lipschitz continuous in changes of demands in some previous periods. In Section \ref{sec:SPsolution}, we prove that the target stability conditions are indeed satisfied with the aforementioned development of finite-dimensional LP approximations.
%(Theorem \ref{thm:SPconti})
We support our analysis with numerical studies in Section \ref{sec:simulation} and conclude the paper in Section \ref{sec:conclusion}.

As for notation,
$\mathbb R^l$ and $\mathbb R_+^l$ are respectively sets of $l$-dimensional real vectors and non-negative real vectors ($l \ge 1$).
Their superscripts are omitted when $l=1$. 
We define $\mathbf 1 \{\}$ to be an indicator function,
which equals 1 if the statement inside the bracket is true and 0 otherwise. 
The maximum and minimum of $x_1$ and $x_2$ are denoted by $x_1 \vee x_2$ and $x_1 \wedge x_2$ respectively, and $\max(x,0)$ is denoted by $x^+$.
Vector symbols are always in bold,
and as two special vectors,
$\mathbf e_j$ is the unit vector with the $j^{th}$ element taking the value of unity,
and $\vec{\mathbf 1}$ is the vector of all $1$s
(dimensions of both vectors depend on the context). 
%The $L$ norm of a vector $\x$ is denoted by $||\x||_L$ ($L=1,\cdots,\infty)$.
The norm $||\x||_{\beta}$ on $\mathbb R^l$ is defined by 
$||\x||_{\beta}= (\sum_{i=1}^l |x_i|^{\beta} )^{1/\beta}$ for $l \ge 1$ and $\beta \ge 1$.
For each pair of vectors $\x_1$ and $\x_2$, 
the maximum and minimum,
$\x_1 \vee \x_2$ and $\x_1 \wedge \x_2$,
are taken componentwise.
Between any pair of vectors,
$\x_1 \ge (\le) \x_2$ if every component of $\x_1$ is greater (less) than or equal to its corresponding component in $\x_2$,
and $\x_1 \ne \x_2$ unless each component in $\x_1$ equals the corresponding component in $\x_2$.

To improve the flow and avoid distracting from our main results, we leave proofs of all lemmas and some theorems in Appendix A. To highlight the major ideas of our work, we  include proofs of key results, Theorems \ref{thm:final} and \ref{thm:master}, as well as corollaries \ref{cor:invtra} and \ref{cor:backlog} of the latter theorem, in the main body of the paper.

 \section{Problem formulation}
\label{sec:formulation}

We consider continuous-review ATO inventory systems.  Inventories are non-perishable and unserved demands are always backlogged.  Component lead times are deterministic but may differ from each other. Our policy and analysis apply to periodic-review systems,  but do not extend to cases with perishable inventories, lost sales, or stochastic lead times. 

\subsection{The System}
\label{subsec:atodesc}

There are $m$ products and $n$ components. The Bill of Materials, given as an $n \times m$ non-negative integer matrix, $A$, specifies the usage of components by different products.  Elements of $A$, $a_{ji}$, represents the amount of component $j$ needed to assemble product $i$ ($\irange$). Thus the $j^{\mbox{th}}$ row of $A$, $\mathbf A_j$, specifies the amounts of component $j$ needed by all products ($\jrange$).  Without  loss of generality, we assume that there is at least one non-zero entry in every row and every column of $A$. 

\begin{figure}
\begin{center}
\includegraphics[scale=0.75]{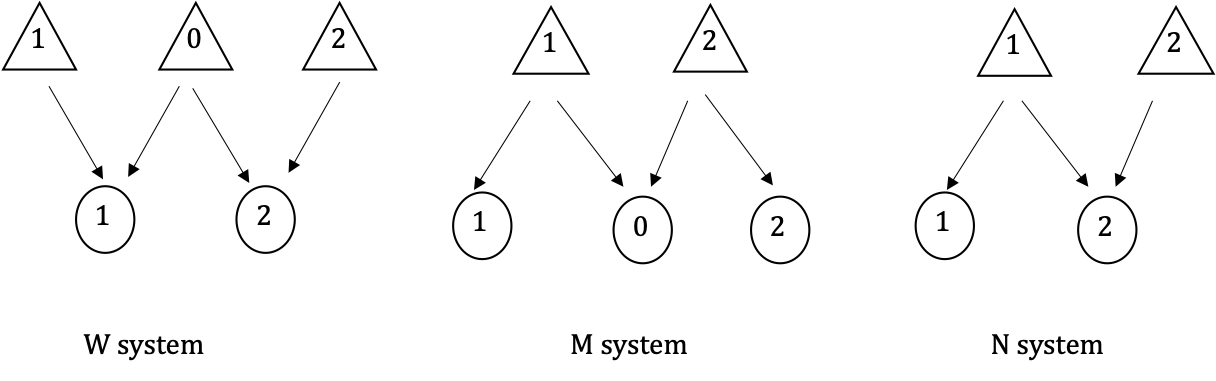}
\caption{Examples: the $W$, $M$, and $N$ systems}
\label{fig:spesys}
\end{center}
\end{figure}

Figure \ref{fig:spesys} shows three  ATO systems that are common in the literature.  The $W$ system has two products, and three components, one of which is used by both products.  Each product is assembled from one unit of the common component and one unit of a product-specific component.   With a slight deviation from the aforementioned notation, the common component is referred to as component $0$.  The $M$ system uses two components to build three products. Products 1 and 2 use one unit of components $1$ and $2$ respectively.   The third product uses one unit of both components and, also with a slight deviation from the aforementioned notation, is referred to as product 0.  The $N$ system is a special case of both the W and M systems. It has two products and two components. Product 0 uses one unit of both components $0$ and $1$ and Product 1 uses only one unit of component $0$. 

There are $K$ distinct component lead times, $L_1<\cdots<L_K$.  Define $L_0=0$ for notational convenience.  Let $n_k$ be the number of components with lead time $L_k$ ($1 \le k \le K$). Components are indexed according to an ascending order of their lead times. Let $\bar{n}_0=0$ and $\bar{n}_k=\sum_{k' \le k} n_{k'}$ ($1 \le k \le K$). Observe that $\bar{n}_K=n$.  Thus $\{\bar{n}_{k-1}+1,\cdots,\bar{n}_k\}$ are the indexes of components with lead times $L_k$ ($1 \le k \le K$).  We associate each component $j$ with an index $k_j$ ($1 \le k_j \le K$) such that $L_{k_j}$ is the lead time of component $j$ ($\jrange$).  Without loss of generality, we arrange the rows of $A$ in an order such that the submatrix $A^k$,  composed of rows $\bar{n}_{k-1}+1,..,\mbox{and}~\bar{n}_k$ of $A$, specifies the usage of components with lead time $L_k$ ($1 \le k \le K$).

Without loss of generality and for simplicity,  we let the system start at a time when there is no inventory on-hand, no order in transit, and no backlog, and define that time to be $t=-L_K$.
Demand arrives according to an integer vector valued compound Poisson process 
\[
\mathcal D(t) = (\mathcal D_1(t), \cdots, \mathcal D_m(t)),~t \ge -L_K,
\]
where $\mD_i(t)$ is the amount of demand for product $i$ ($\irange$) that arrives during $[-L_K,t]$. 
Note that the process $\{ \mathcal D(t), t \ge 0 \}$ is right continuous.
The number of demand orders arriving  during  $[ -L_K,t]$ ($t \ge -L_K$) is a Poisson process  $\Lambda=\{ \Lambda(t), t \ge -L_K \}$, and there is an associated i.i.d sequence of random vectors that give order sizes. A generic element of this sequence is denoted by  $\mathbf S=(S_1,S_2,...,S_m)$, where $S_i$ is the order size for product $i$. Although the elements of the sequence (the order size vectors)  are independent, the components $(S_1,S_2,...,S_m)$ within each vector can be dependent. Let  
\[
\underline{\lambda}=\E[\Lambda(1)-\Lambda(0)]
\]
denote the order arrival rate. Mean demand arriving within a unit of time is
\[
\boldsymbol \mu =(\mu_1,\cdots,\mu_m):=\E[\mD(1)-\mD(0)]= \underline{\lambda} \E[\mathbf S].
\]
The covariance matrix of $(\mD(1)-\mD(0))$ is denoted by $\boldsymbol \Sigma$, of which the diagonal elements, $\sigma_{ii}$, are variances of demand $i$ ($\irange$) over a unit of time.  Since the demand process is stationary, $\boldsymbol \mu$ and $\boldsymbol \Sigma$ are also respectively the means and the covariance matrix of demands over $[t,t+1]$ ($t \ge -L_K$).

We assume that $S_i$ has a finite moment of order $6$, i.e.,
\begin{equation*}
\eta_i :=\E[S_i^6]<\infty, \quad \irange. 
\end{equation*}
As will be evident from the proof of Theorem \ref{thm:master}, we assume a finite moment of order $6$ for a similar reason as that in \cite{Huh2009b}, which is to use the moment to bound the difference between two stochastic processes. (Beyond that, both the time horizon and stochastic processes involved are completely different between the two studies.) 

\subsection{The Inventory Control Problem}
\label{subsec:invprodef}

Inventory control includes both replenishment and allocation decisions. For components with lead time $L_k$, the replenishment starts from time $-L_k$ and is defined by $\mathcal R^k(t)$ $(t \ge -L_k)$, an integer-valued vector process with $n_k$ components. Each component of $\mathcal R^k(t)$, $\mathcal R_j(t)$, represents the amount of component $j$ ordered over the period $[-L_k,t]$ $(\bar{n}_{k-1}+1 \le j \le \bar{n}_k$). In this context, $\mathcal R^k(t)$ $(t \ge -L_k)$ should be non-negative and non-decreasing, which we assume in the paper. Orders are placed at distinct points of time. Hence we assume the process is right-continuous. 

The allocation decision is specified via an $m$-dimensional, integer-valued process $\mathcal Z(t)$.  Each component of $\mathcal Z(t)$, $\mathcal Z_i(t)$ $(t \ge -L_K)$, represents the amount of demand for product $i$ $(\irange)$ served over the period $[-L_k,t]$.  We assume that $\mathcal Z(t)$ $(t \ge -L_K)$ is also non-negative, non-decreasing, and right- continuous. 

Both $\mathcal R^k(t)$ ($t \ge -L_k$, $1 \le k \le K$) and $\mathcal Z(t)$ ($t \ge 0$) must be adapted to the filtration generated by the initial states of the system, as well as $\mathcal D(s)$ ($-L_K \le s \le t$), $\mathcal R^k(s)$ ($-L_k \le s <t$, $1 \le k \le K$), and $\mathcal Z(s)$ ($-L_K \le s <t$), which means that policies under our consideration are non-anticipating. 

For the inventory control to be feasible, at any point of time, the total amounts of demand served cannot exceed the amounts that have arrived. The amounts of components used cannot exceed the amounts that have been received. Since the replenishment of components with the lead time $L_k$ starts at time $-L_k$ $(1 \le k \le K)$,  no unit is received and thus no demand can be served before time $0$. In our model, these restrictions are formalized by the following assumptions:
\begin{align}
\label{eq:feasibility}
\begin{split}
& \mathcal Z(t)=0,~~t<0, \quad \mathcal Z(t) \le \mathcal D(t),~~ t \ge 0,
\\
\mbox{and} \quad & A^k \mathcal Z(t) \le \mathcal R^k(t-L_k), ~1\le k \le K,~t \ge 0.
\end{split}
\end{align}
Observe that by definition, $\mR^k(t-L_k)$ $(1 \le k \le K)$ are the amounts of components ordered at least one lead time before time $t$.  

For $t \ge 0$, the inventory level of components with the lead time $L_k$  is
\begin{equation}
\label{eq:invlvl}
\mathbf I^k(t):=\mathcal R^k(t-L_k) - A^k \mathcal Z(t),\quad k=1,\cdots, K,
\end{equation}
and the backlog level is 
\begin{equation}
\label{eq:blklvl}
\mathbf B(t) :=\mathcal D(t) -\mathcal Z(t). 
\end{equation}
By (\ref{eq:feasibility}), both $\mathbf I^k(t)$ $(1 \le k \le K)$ and $\mathbf B(t)$ $(t \ge 0)$ are nonnegative. Also by (\ref{eq:feasibility}) and our assumption about the replenishment starting times, 
\[
\mathbf I^k(t)=0,~-L_k \le t <0, \quad \mbox{and} \quad \mathbf B(t)=\mathcal D(t),~~t<0.
\]

Let $b_i$ be the cost of backlogging one unit of demand of product $i$ ($\irange$) per unit of time. Let $h_j$ be the cost of holding one unit of inventory of component $j$ ($\jrange$) per unit of time. We assume that $h_j$ ($\jrange$) and $b_i$ ($\irange$) are strictly positive.  Let
\[
\mathbf b=(b_1,\cdots,b_m)~~~\mbox{and}~~~\mathbf h^k=(h_{\bar{n}_{k-1}+1},\cdots,h_{\bar{n}_k}),~1 \le k \le K. 
\]
Then at each time $t$ ($t \ge 0$), the system incurs the total inventory cost at the expected rate 
\begin{equation}
\label{eq:timecost}
\mathcal C(t) = \sum_{k=1}^K \mathbf h^k \cdot \E[\mathbf I^k(t)]+\mathbf{b}\cdot \E[\mathbf{B}(t)].
\end{equation}
The problem of inventory control is to develop integer-valued, non-negative, non-decreasing, right-continuous, and non-anticipating vector processes $\mathcal R^k(t)$ ($t \ge -L_k$, $1 \le k \le K$) and $\mathcal Z(t)$ ($t \ge -L_K$), subject to (\ref{eq:feasibility})-(\ref{eq:blklvl}), to minimize the following long-run average expected total inventory cost:
\begin{equation}
\label{eq:LAC}
\mathcal C =\limsup_{T \to \infty} {1 \over T} \int_0^{T} \mathcal C(t) dt.
\end{equation}

\subsection{Additional Variables, Processes, and Relationships}

For the needs of our analysis,  we introduce some additional variables and processes, and specify their relationships. 

\subsubsection{Additional Variables and Processes}
The additional variables and processes, defined in terms of the more 'primitive' processes introduced in Sections \ref{subsec:atodesc} and \ref{subsec:invprodef} are listed in Table \ref{tbl:notations} along with their definitions. For the convenience of the reader, Table \ref{tbl:notations} also includes the definitions of the inventory and backlog processes defined in Section 2.2.
% lists key variables and their values. 
Further discussions of these quantities are as follows:

The first processes we introduce here denote increments of the processes $\mathcal D(t)$ $(t \ge -L_K)$, $\mR^k(t)$ $(t \ge -L_k, 1 \le k \le K)$, and $\mathcal Z(t)$ $(t \ge -L_K)$ over particular time intervals. Here $\mathbf D(t_1,t_2)$ denotes demand that arrives during the time interval $(t_1,t_2]$ ($-L_K \le t_1 <t_2)$, $\mathbf D^k(t)$ are special cases of $\BD(t_1,t_2)$ with $t_1=t-L_k$ and $t_2=t-L_{k-1}$ $(t \ge -L_k+L_K$, $1 \le k \le K$), and 
\[
 \bar{\BD}^k(t)=\sum_{k'=k}^K \BD^{k'}(t),\quad t \ge 0.
 \]
We also use $\bar{\BD}(t)$ as a shorthand of $\bar{\BD}^1(t)$ $(t \ge 0)$.   (Recall that $L_0=0$.) 
We also introduce some random vectors. Let $\BD^k$, $\bar{\mathbf D}^k$, and $\underline{\BD}^k$ be random vectors that follow the same distributions of  $\BD^k(t)$, $\bar{\BD}^k(t)$, and $\underline{\BD}^k(t)$, $1 \le k \le K$, $t \ge 0$, respectively. Let $\bar{\BD}$ be a random vector that follows the same distribution of $\bar{\BD}(t)$.   By stationarity of $\mD(t)$, these distributions do not depend on $t$. 
Recall that we assumed that the demand process is right continuous. We let $\mathbf d(t)$ denote the demand that arrives at time $t, t \ge -L_K$ (if any).

Similarly, $\mathbf R^k(t)$ $(1 \le k \le K)$ denotes the amounts of components ordered over the previous lead time, or for $t<0$, since the beginning of the replenishment process. Each component of $\mathbf R^k(t)$, $R_j(t)$ $(t \ge -L_k, \bar{n}_{k-1}+1 \le j \le \bar{n}_k)$, represents the amount of component $j$ that is in-transit: ordered from the supplier but not yet arrived. The amounts of demand served over the interval $(t_1,t_2]$ is denoted by $\mathbf Z(t_1,t_2)$ $(-L_K \le t_1 <t_2)$. 

\begin{table}[ht]
\begin{center}
\begin{tabular}{|l|l|l|}
\hline & & \\ [-0.3cm]
&~~~~~~~~~variables & ~~~~~~~~~~~definitions \\ [0.3cm] \hline & & \\ [-0.3cm]
demand 
 & $\mathbf D(t_1,t_2)$, $-L_K \le t_1 < t_2$ & $:=\mD(t_2)- \mD(t_1)$ \\ [0.3cm]
 %\cline{2-3}
& $\mathbf D^k(t)$, $t \ge -L_K+L_k$, $1 \le k \le K$ & $:=\mD(t-L_{k-1}) -\mD(t-L_k)$  \\ [0.3cm]
%\cline{2-3}
& $\bar{\mathbf D}^k(t)$, $t \ge 0$, $1 \le k \le K$ & $:=\mD(t-L_{k-1})-\mD(t-L_K)$ \\ [0.3cm]
%\cline{2-3}
& $\bar{\BD}(t)$, $t \ge 0$ & $:=\mD(t)-\mD(t-L_K)$ \\ [0.3cm]
%\cline{2-3}
& $\mathbf d(t)$, $t \ge -L_K$ &  $:=\left \{\begin{array}{ll} \mD(-L_K), & t=-L_K \\ \mD(t)-\mD(t^-), & t >-L_K \end{array} \right .$ \\ [0.5cm]
%\cline{2-3}
%\cline{2-3}
& $\mathbf D^k$, $1 \le k \le K$ &  $\stackrel{d}{=} \mathbf D^k(t)$ \\[0.3cm]
%\cline{2-3}
& $\bar{\BD}^k$, $1 \le k \le K$ & $\stackrel{d}{=} \bar{\mathbf D}^k(t)$ \\ [0.3cm] 
& $\bar{\BD}$ & $\stackrel{d}{=} \bar{\mathbf D}(t)$ \\[0.3cm] 
& $\underline{\BD}^k$, $1 \le k \le K$ & $:=\BD_k+\cdots \BD_1 \stackrel{d}{=} \BD(t-L_k,t)$ \\ [0.3cm] \hline & & \\[-0.3cm]
replenishment 
&  $\mathbf R^k(t)$, $t \ge -L_k,~1 \le k \le K$ & $:= \left \{\begin{array}{ll} \mR^k(t)-\mR^k(t-L_k), &t \ge 0 \\
 \mR^k(t), &t <0 \end{array} \right.$  \\ [0.3cm] \hline  & &  \\[-0.3cm]
 allocation 
 & $\mathbf Z(t_1,t_2)$, $-L_K \le t_1 <t_2$ & $:=\mathcal Z(t_2)- \mathcal Z(t_1)$ \\ [0.3cm] \hline  & & \\  [-0.3cm]
 inventory & $\mathbf I^k(t)$, $t \ge 0$, $1 \le k \le K$ & $:=\mR^k(t-L_k) - A^k \mZ(t)$ 
 %$\left \{\begin{array}{ll} 0 \\ \mR(t-L_k)-A^k \mZ(t), & t \ge 0 \end{array} \right .$ 
 \\ [0.3cm] \hline  & & \\ [-0.3cm]
 backlog & $\mathbf B(t)$, $t \ge 0$ & $:=\mD(t)-\mZ(t)$\\  [0.3cm] 
 & $\mathbf B^-(t)$, $t \ge 0 $ & $:=\left \{\begin{array}{ll} \mD(t), & t=0 \\  \mathbf B(t^-)+ \mathbf d(t), & t >0 \end{array} \right .$ \\ [0.3cm] \hline  & & \\ [-0.3cm]
  inventory position & $\ip^k(t)$, $t \ge -L_k$, $1 \le k \le K$ & $:=\mR^k(t)-A^k \mD(t)$  \\ [0.3cm]
& $\ipp(t)$, $t \ge -L_k$, $1 \le k \le K$ & $:=\left \{\begin{array}{ll}  -A^k \mD(-L_k), & t=-L_k \\  \ip^k(t^-)-A^k \mathbf d(t), & t >-L_k \end{array} \right .$ \\ [0.3cm] \hline & & \\ [-0.2cm]
inventory balance & $\mathbf Q^k(t)$, $t \ge 0$, $1 \le k \le K$, & $:=A^k \mD(t)-\mR^k(t-L_k)$ \\ [0.3cm] \hline
\end{tabular}
\end{center}

\smallskip
\caption{Key Variables: Notation and Value}
\label{tbl:notations} 
\end{table}

Next we introduce two additional processes: inventory position $\ip^k(t)$ and component balance $\mathbf Q^k(t)$ $(t \ge -L_k$, $1 \le k \le K)$. For $k=1,\cdots, K$, entries of $\ip^k(t)$, $IP_j(t)$ ($\bar{n}_{k-1} +1 \le j \le \bar{n}_k$), are the differences between the total amounts of components that have been ordered up to time $t$ $(t \ge -L_k)$ and the total amounts needed to serve all demands that have arrived by that time. Entries of $\mathbf Q^k(t)$, $Q_j(t)$ ($\bar{n}_{k-1}+1 \le j \le \bar{n}_k)$, are the differences between the amounts needed to serve demands and the total amounts that have been ordered and received.  

%It also helps to distinguish values of a same variable at any particular time before and after relevant actions have been taken.
Recall that the replenishment and allocation processes (which are the control processes) are assumed to be integer valued and right continuous. Thus the controls are exercised at discrete time points. It is helpful to distinguish the values of certain processes immediately before these actions are taken.
To this end, we let $\ipp(t)$ $(1 \le k \le K$) denote the inventory positions at time $t$ $(t \ge -L_k)$ before placing any order at that time, and let $\mathbf B^-(t)$ $(t \ge -L_K)$ denote the backlog levels at time $t$ $(t \ge -L_K$) before serving any demand at that time. 
In Table \ref{tbl:notations}, $\mathbf d(t)$, $\ipp(t)$, and $\mathbf B^-(t)$ are defined using the left limits $\mD(t^-)$ $(t > -L_K)$, $\ip^k(t^-)$ $(t > -L_k$, $1 \le k \le K$), and $\mathbf B(t^-)$ $(t >-L_K)$.  These left limits exist because in (\ref{eq:invlvl}) and (\ref{eq:blklvl}),  $\mD(t)$ $(t \ge -L_K)$, $\mR^k(t)$ $(t \ge -L_k$, $1 \le k \le K$), and $\mZ(t)$ $(t \ge -L_K)$ are non-decreasing and right-continuous. One may also observe that by definition, 
\begin{align*}
\begin{split}
& \ip^k(t)  =\ipp(t)+\mR^k(t)-\mR^k(t^-), \quad t \ge -L_k,~1 \le k \le K \\
\mbox{and} \quad & \mathbf B(t)=\mathbf B^-(t)- (\mZ(t)-\mZ(t^-)), \quad t \ge -L_K. 
\end{split}
\end{align*}

\subsubsection{Relationships} 
The key definitions (\ref{eq:invlvl}) and (\ref{eq:blklvl}), along with the definitions introduced in Section 2.3.1 allow us to obtain relationships that are useful in our analysis.
Using (\ref{eq:invlvl}) and (\ref{eq:blklvl}), and the definitions of $\mathbf R^k(t)$ $(t \ge -L_k$, $1 \le k \le K$) and $\mathbf Z(t)$ $(t \ge -L_K)$ in Table \ref{tbl:notations},  changes of inventory and backlog levels over a lead time $L_k$ can be determined by 
%Over a component's lead time, the change of its on-hand inventory is determined by its starting inventory level (both on-hand and in-transit) and its usage during that period. Changes of backlog levels over any period are determined by the differences between demands arrived and served during that period. Therefore at each time $t \ge 0$ and for all $k=1,\cdots, K$, 
\begin{align}
\label{eq:flow}
\begin{split}
\mathbf I^k(t)&=\mathbf I^k(t-L_k)+\mathbf R^k(t-L_k)-A^k \BZ(t-L_k,t), \quad t \ge 0, ~1 \le k \le K,
\\
\mbox{and} \quad
\mathbf B(t)&=\mathbf B(t-L_k)+\mathbf D(t-L_k,t)-\BZ(t-L_k,t), \quad t \ge -L_K+L_k. 
\end{split}
\end{align}

Again using (\ref{eq:invlvl}) and (\ref{eq:blklvl}), along with  definitions of $\mathbf R^k(t)$ and $\ip^k(t)$ $(t \ge -L_k$, $1 \le k \le K)$ in the table, 
\begin{align}
\begin{split}
\label{eq:IPdef}
\ip^k(t) &= \mathbf I^k(t) +(\mR^k(t)-\mR^k(t-L_k)) -A^k \mathbf B(t) \\
& = \mathbf I^k(t)+\mathbf R^k(t)-A^k \mathbf B(t),~~t\ge -L_k. 
\end{split}
\end{align}
The second line of the above corresponds to another definition of the inventory position in the literature: the total amount of the component in the system, including both on-hand inventory and orders in-transit, minus the amount needed to clear all existing backlog at a given time.  

Yet again using (\ref{eq:invlvl}), (\ref{eq:blklvl}), along with the definition of $\mathbf Q^k(t)$ in the table, and then applying (\ref{eq:flow}) and (\ref{eq:IPdef}) to replace $\mathbf B(t)$ and $\mathbf I^k(t)$, 
\begin{align}
\label{eq:balancedef}
\begin{split}
\mathbf Q^k(t)&=A^k \mathbf B(t)-\mathbf I^k(t)\\
&=A^k \BD(t-L_k,t)-\ip^k(t-L_k), \quad t \ge 0,~1 \le k \le K.
\end{split}
\end{align}
Observe that the negative  of an entry of $\mathbf Q^k(t)$, 
\begin{equation*}
-Q_j(t)=I_j(t)-\mathbf A_j \cdot \mathbf B(t),\quad \bar{n}_{k-1}+1 \le j \le \bar{n}_k,
\end{equation*}
is commonly referred to as the net inventory of component $j$ at time $t$ $(t \ge 0)$:
there is a shortage (surplus) of component $j$ to clear existing backlogs at time $t$ if $Q_j(t)>(<)0$. 

\subsubsection{Other parameters}
Define
\begin{equation}
\label{eq:defC}
\mathbf c=(c_1,\cdots,c_m) \equiv \mathbf b+ \sum_{k=1}^K (A^k)' \mathbf h^k,
\end{equation}
where element $c_i$ of $\mathbf c$ represents the amount of inventory cost that can be removed from the system by serving one unit of demand $i$ ($\irange$).

For convenience, we define additional variables to denote the smallest (non-zero) and largest components of $A$, $\mathbf b$, $\mathbf h$, and $\mathbf c$ in the following table:

\begin{table}[ht]
\begin{center}
\begin{tabular}{|c|c|c|c|c|c|c|c|}
\hline 
& & & & & & & \\[-0.3cm]
symbol & $\underline{a}$ & $\bar{a}$ & $\bar{h}$ & $\underline{h}$ & $\bar{b}$ & $\underline{b}$ & $\bar{c}$ \\[0.1cm] \hline
& & & & & & & \\[-0.3cm]
definition & $\min_{i,j:a_{ij}>0} \{a_{ij}\}$ & $\max_{i,j} \{a_{ij}\}$ & $\max_{j} \{h_j\}$ & $\min_{j} \{h_j\}$ & $\max_{i} \{b_i\}$ & $\min_{i} \{b_i\}$ & $\max_{i} \{c_i\}$ \\[0.1cm] \hline
\end{tabular}

\medskip
\caption{Smallest and Largest Components of Relevant Vectors and Matrix}
\label{tbl:minmaxpara}
\end{center}
\end{table}

\section{Stochastic Program}
\label{sec:SPdisc}

As we outlined in Section \ref{sec:intro}, the first step of our analysis is to develop a multi-stage SP that provides a lower bound on the average inventory cost and sets the stage for developing inventory control policies to drive the cost towards that bound.  To this end, we first present a relevant SP from the literature and discuss the intuition that underlies its formulation in  Section \ref{subsec:preSP}. In Section \ref{subsec:newdev} we transform this SP into an alternative SP that we use for policy development in Section \ref{sec:policy}. 

\subsection{Previous Result}
\label{subsec:preSP}

For ATO systems formulated in the last section, 
Theorem 1 in \cite{Reiman2012} proves that 
\begin{equation}
 \label{eq:lbound}
\underline{C}= \inf_{\bal \ge 0} \{\mathbf b \cdot \bal+\Phi^K(\bal)\}+\mathbf b \cdot \E[\bar{\BD}],
\end{equation}
is a lower bound on $\mathcal C$, the long-run average cost defined in (\ref{eq:LAC}). 
Here $\bar{\BD}$ is defined in Table \ref{tbl:notations} and $ \Phi^K(\bal)$ is the minimum objective value of the $K+1$ stage stochastic program (SP):
 \begin{align}
 \label{eq:MSP1}
 \begin{split}
 \Phi^K(\bal) &= \inf_{\mathbf{y}^K \ge 0} \{\mathbf{h}^K \cdot \mathbf{y}^K + \E [\Phi^{K-1}( \y^K,\bal + \BD^K)  ]\},
 \\
 \Phi^k(\y^{k+1},\cdots, \y^K,\x) &= \inf_{\mathbf{y}^k \ge 0} \{\mathbf{h}^k \cdot \mathbf{y}^k + \E [ \Phi^{k-1}(\y^k,\cdots, \y^K,\x + \BD^k)]\},~~k=K-1,\cdots,1,
 \\
 \Phi^0(\y^1,\cdots, \y^K,\mathbf {x})& = - \max_{\mathbf{z} \ge 0} \{ \mathbf{c}\cdot \mathbf{z}  |  \mathbf{z} \le \mathbf{x}, A^k \mathbf{z} \le \y^k, 1 \le k \le K \},
 \end{split}
 \end{align}
where $\BD^k$ ($1 \le k \le K$) are also defined in Table \ref{tbl:notations}. 

We now  take the rest of this subsection to provide some intuition as to why the SP (\ref{eq:lbound})-(\ref{eq:MSP1}) yields a lower bound on the cost in the ATO inventory control problem. The basic idea is that the SP (\ref{eq:lbound})-(\ref{eq:MSP1}) is a relaxation of the actual ATO inventory control problem. The SP corresponds to myopically focusing on a single point of time, with no concern for the effect that any replenishment or allocation decision might have on costs at other points of time. Thus, in the presence of random demand, all decisions are put off until the last possible moment, allowing as much information about actual demand as possible to be used for the decisions.
%(The discussion here is not meant to be rigorous. For a rigorous proof, see \cite{Reiman2012}.)

To explain the SP in more detail, using (\ref{eq:flow}) and (\ref{eq:defC}), with $k=K$ when substituting for $\mathbf B(t)$ and 
\[
\mathbf Z(t-L_k,t)=\mathbf Z(t-L_K,t)-\mathbf Z(t-L_K,t-L_k)
\]
when substituting for $\mathbf I^k(t)$ $(1 \le k \le K)$, we can write the inventory cost at any time $t$ $(t \ge 0)$ as
\begin{align*}
\begin{split}
\sum_{k=1}^K \mathbf h^k \cdot \mathbf I^k(t)+\mathbf b \cdot \mathbf B(t) &=\mathbf b \cdot \mathbf B(t-L_K)- \mathbf c \cdot \mathbf Z(t-L_K,t)
\\
+ &\sum_{k=1}^K \mathbf h^k \cdot \left [\mathbf I^k(t-L_k)+ \mathbf R^k(t-L_k)+A^k \mathbf Z(t-L_K,t-L_k)\right] +\mathbf b \cdot \BD(t-L_K,t).
 \end{split}
 \end{align*}
Notice that the cost at $t$ depends only on states and processes over the period $[t-L_K,t]$. Comparing (\ref{eq:lbound})-(\ref{eq:MSP1}) with this expression for the cost, $\bal$ corresponds to the initial backlog levels $\mathbf B(t-L_K)$ and $\z$ corresponds to $\mathbf Z(t-L_K,t)$, the amounts of demand served in the period $[t-L_K,t]$. For components with lead time $L_k$, $\mathbf I^k(t-L_k)+ \mathbf R^k(t-L_k)+A^k \mathbf Z(t-L_K,t-L_k)$ includes the amounts in inventory at time $t-L_k$, the amounts that will arrive between $t-L_k$ and $t$ from the pipeline, and the amounts used in the period $[t-L_K,t-L_k]$ $(1 \le k \le K)$. Corresponding to the sum of these three quantities, $\y^k$ represents the total amounts of these components that can be used to serve demands in period $[t-L_K,t]$. 

In (\ref{eq:lbound})-(\ref{eq:MSP1}), $\bal$, $\z$, and $\y^k$ $(1 \le k \le K)$ are chosen to minimize the inventory cost with no constraint other than the ones that must be satisfied by the aforementioned counterparts of these variables in the ATO system. By (\ref{eq:feasibility})-(\ref{eq:blklvl}), $\mathbf B(t)$ and $\mathbf I^k(t)$ $(1 \le k \le K$, $t \ge 0$) are nonnegative, so using (\ref{eq:flow}), 
\[
 \mathbf Z(t-L_K,t)  \le \mathbf B(t-L_K) +\BD(t-L_K,t),
\]
and
\begin{align*}
\begin{split}
A^k \mathbf Z(t-L_k,t)  & \le \mathbf I^k(t-L_k) + \mathbf R^k(t-L_k) \\ 
\mbox{i.e.}, \quad A^k \mathbf Z(t-L_K,t)  &\le \mathbf I^k(t-L_k)+\mathbf R^k(t-L_k)+A^k \mathbf Z(t-L_K,t-L_k),\quad 1 \le k \le K.
\end{split}
\end{align*}
Correspondingly,  in (\ref{eq:MSP1}),
\[
\mathbf z \le \bal+\bar{\BD}= \bal +\BD^K+\cdots +\BD^1 \quad \mbox{and} \quad A^k \mathbf z \le \y^k~ (1 \le k \le K).
\]

The information constraint in the SP is less obvious, as the sequential/recursive nature of its formulation both encodes and potentially obscures the information available when certain decisions are made. Although it is not needed in the definition of the SP, the information available when decisions are made in the SP can be described via a discrete filtration
\begin{equation}
\label{eq:filtration}
\mathcal F_K \subseteq \mathcal F_{K-1} \subseteq \cdots \subseteq \mathcal F_0,
\end{equation}
where $\mathcal F_K=\{\emptyset,\Omega\}$ and $\mathcal F_k$ is the $\sigma$-field generated by $\{\BD^{K},\cdots \BD^{k+1},\y^K,\cdots \y^{k+1}\}$ ($0\le k <K$). The set of optimal values of $\y^k$ is $\mathcal F_k$-measurable ($1\le k \le K$) and that of $\z$ is $\mathcal F_0$-measurable.
 
The filtration $\mathcal F_k$ $(0 \le k \le K$) imitates the information available in the associated ATO system, and is generated by both $\mathbf D^{k'}$ and $\y^{k'}$, corresponding respectively to histories of demand arrivals in periods $[t-L_{k'},t-L_{k'-1}]$ and decision making at times $t-L_k'$ $(k<k' \le K)$. In ATO systems, $t-L_k$ is the last moment of decision-making that can affect the value of $\mathbf I^k(t-L_k)+ \mathbf R^k(t-L_k)+A^k \mathbf Z(t-L_K,t-L_k)$, so letting $\y^k$ be adapted to $\mathcal F_k$ allows its value to be chosen with the maximum amount of information.  Note that this adaptedness is not explicitly enforced in the SP, but is implicit in the recursive structure of the SP. When $ \Phi^k$ (which corresponds to stage $K+1-k$) is solved to obtain $\y^k$, the prior decisions $\y^K, \ldots,\y^{k+1}$ are all known, and $\x=\BD^{K}+\cdots +\BD^{k+1}$ is known as well. Similarly, $t$ is the last moment of decision-making that can affect $\mathbf Z(t)$, so the choice of $\z$ is adapted to $\mathcal F_0$. It is not surprising that, with $\bal$, $\z$, and $\y^k$ $(1 \le k \le K)$ chosen under the minimum constraints and maximum information, (\ref{eq:lbound}) is a lower bound on the expected cost of the ATO system at any given time $t$, and hence a lower bound on its average cost $\mathcal C$. 

\subsection{New Development}
\label{subsec:newdev}

The SP in (\ref{eq:MSP1}) is not directly applicable to our analysis. As one may observe from (\ref{eq:lbound})-(\ref{eq:MSP1}), and as is shown in the proof of Theorem \ref{thm:sameobj},  $\mathbf b \cdot \bal+\Phi^K(\bal)$ decreases in $\bal$, and strictly so in many cases.  Therefore, the lower bound $\underline{C}$ is not directly computable by solving (\ref{eq:MSP1}) with any fixed $\boldsymbol \alpha$. To reach the infimum in (\ref{eq:lbound}), $\bal$, and thus $\y^k$ ($1 \le k \le K$) and $\z$ often have to approach infinity, so one cannot make use of the optimal values of these decision variables for policy development. 

To address this issue, we transform (\ref{eq:MSP1}) into an alternative SP by replacing $\y^k$ in (\ref{eq:MSP1}) with $\y^k -A^k \bal$ $(1 \le k \le K)$ and $\z$ with $\z-\bal$ and, following (\ref{eq:lbound}), letting $\bal$ approach infinity. The transformation removes $\bal$ and the non-negativity constraints in (\ref{eq:MSP1}) and leads to:
%This gives rise to the following SP, which differs from (\ref{eq:MSP1}) by the absence of $\bal$ and non-negativity constraints:
\begin{align}
\label{eq:MSP2}
\begin{split}
\varphi^K &= \inf_{\y^K \in \mathbb R^{n_{\scaleto{K}{2.5pt}}}} \{\mathbf{h}^K \cdot \mathbf{y}^K + \E [\varphi^{K-1}( \y^K, \BD^K)  ]\},
\\
\varphi^k(\y^{k+1},\cdots, \y^K,\x) &= \inf_{\y^k \in \mathbb R^{n_{\scaleto{k}{2.5pt}}}} \{\mathbf{h}^k \cdot \y^k + \E [ \varphi^{k-1}(\y^k,\cdots, \y^K,\x+\mathbf{D}^k)]\},~1 \le k <K,
\\
\varphi^0(\y^1,\cdots, \y^K,\x)& = - \max_{\z \in \mathbb R^m} \{ \mathbf{c}\cdot \mathbf{z}  |  \mathbf{z} \le \mathbf{x}, A^k \mathbf{z} \le \y^k,~1 \le k \le K \}.
\end{split}
\end{align}
%Similar to (\ref{eq:filtration}),  we define a filtration that starts with $\mathcal F_K=\{\emptyset,\Omega\}$ and generates $\sigma-$field $\mathcal F_k$ by $\{\BD^{K},\cdots \BD^{k+1},\y^K,\cdots \y^{k+1}\}$ ($0\le k <K$). The set of optimal values of $\y^k$ is $\mathcal F_k$-measurable ($1\le k \le K$) and that of $\z$ is $\mathcal F_0$-measurable.

For any feasible solution $\bal$, $\y^k$ $(1 \le k \le K)$, and $\z$ in (\ref{eq:MSP1}), $\y^k-A^k \bal$ $(1 \le k \le K)$ and $\z-\bal$ is a feasible solution to (\ref{eq:MSP2}), so $\varphi^K \le \Phi^K(\bal)$ for all $\bal$. Therefore (\ref{eq:MSP2}) can also be used to set a lower bound on the inventory cost of an ATO system, a result that we will formally present in Theorem \ref{thm:sameobj}. Below we first prove that the optimal solution of (\ref{eq:MSP2}) is bounded. 

%Let\begin{equation}\label{eq:ubd}\underline{\BD}^k=\BD^k+\cdots+\BD^1.\end{equation}%where $\BD^k,\cdots,\BD^1$ are demand inputs to the subsequent stages of the SP. 
\begin{theorem}{(see Appendix A for the Proof)}
\label{thm:solbound}
There exist finite constants $\underline{\beta}>0$ and $\bar{\beta}>0$, depending only on $A$, $\mathbf h$, and $\mathbf b$, such that for any given values of $\y^{k+1},\cdots, \y^K$, and $\x$,  if $\y^{k*}=(y_{\bar{n}_{k-1}+1}^*,\cdots,y_{\bar{n}_k}^*)$ is an optimal solution to the SP in (\ref{eq:MSP2}) at stage $k$, i.e., 
 \[
 \mathbf{h}^k \cdot \y^{k*} + \E [ \varphi^{k-1}(\y^{k*}, \y^{k+1}\cdots, \y^K,\x+\mathbf{D}^k)]
  \]
 attains the value of $\varphi^k(\y^{k+1},\cdots,\y^{K}, \x)$ $(k=K,\cdots,1)$, then 
\begin{equation}
\label{eq:solbound}
-\underline{\beta}(||\x||_1+\E\left[||\underline{\BD}^k||_1\right]+\sum_{l=k+1}^K ||\y^{l}||_1) 
\le y_j^* \le \bar{\beta}(||\x||_1+\E\left[||\underline{\BD}^k||_1\right]+1),~~\bar{n}_{k-1} <j \le \bar{n}_k,
\end{equation}
where $\underline{\BD}^k$ $(1 \le k \le K$) are random vectors defined in Table \ref{tbl:notations}. 
\end{theorem}
%\begin{remark}The proof of the theorem only requires that $\BD^k$ $(1 \le k \le K)$ have finite expected values. This allows us to apply (\ref{eq:solbound}) to the proof of \ref{thm:} where $\BD^k$ $(1 \le k \le K)$ are not  compound Poisson variables,  and to the asymptotic analysis where $\BD^k$ $(1 \le k \le K)$ are centered and scaled, so their values and values of $\x$ can be negative.\end{remark}

The formulation in (\ref{eq:MSP2}) helps us to avoid the aforementioned issue with (\ref{eq:lbound})-(\ref{eq:MSP1}). By (\ref{eq:solbound}), the optimal solutions to the new SP are finite, so they can and will be used as parameters of our inventory policy. Moreover instead of taking the infimum in (\ref{eq:lbound}), we can use the optimal objective value of (\ref{eq:MSP2}), which is directly computable, to set the same lower bound on the cost objective. 
\begin{theorem}{(see Appendix A for the Proof)}
\label{thm:sameobj}
Let $\underline{C}$ be the lower bound defined in (\ref{eq:lbound}) and $\varphi^K$ be the optimal objective value of the SP defined in (\ref{eq:MSP2}). Then
\begin{equation}
\label{eq:sameobj}
\underline{C}=\varphi^K+\mathbf b \cdot \E[\bar{\BD}].
\end{equation}
\end{theorem}

For future analysis, we replace $\z$ with $\mathbf B=\x-\z$ in (\ref{eq:MSP2}) to transform the last stage LP into
\begin{align}
\begin{split}
\label{eq:phi0pi}
\varphi^0(\y^1,\cdots, \y^K,\x)&=\tilde{\varphi}^0(\y^1,\cdots, \y^K,\x)-\mathbf c \cdot \x ,
\\
\mbox{where}~\tilde{\varphi}^0(\y^1,\cdots, \y^K,\x)
&=\min_{\mathbf B} \{\mathbf c \cdot \mathbf B |\mathbf B \ge 0,~ A^k \mathbf B \ge A^k \x-\y^k,~1 \le k \le K \}.
\end{split}
\end{align}
In $\tilde{\varphi}^0(\y^1,\cdots, \y^K,\x)$, $\mathbf B$ represents the amounts of unserved demands, and thus is analogous to $\mathbf B(t)$ ($t \ge -L_K$), backlog levels in the ATO system. Let $\y^{K*}$ be an optimal solution to $\varphi^K$, $\y^{k*}$ be an optimal solution to $\varphi^k(\y^{k+1 *},\cdots,\y^{K*},\bar{\BD}^{k+1})$ ($1 \le k < K$),
$\z^*$ be an optimal solution  to $\varphi^0(\y^{1 *},\cdots,\y^{K*},\bar{\BD})$, 
and $\mathbf B^*=\bar{\BD}-\z^*$. Then by Theorem \ref{thm:sameobj},
\begin{align}
\begin{split}
\underline C &=\sum_{k=1}^K \mathbf h^k \cdot \E[\y^{k*}]-\mathbf c \cdot \E[\mathbf z^*]+\mathbf b \cdot \E[\bar{\BD}]
\\
&=\sum_{k=1}^K \mathbf h^k \cdot \E[\y^{k*}]+\mathbf c \cdot \E[\mathbf B^*]-\sum_{k=1}^K [(A^k)' \mathbf h^k] \cdot \mathbf \E[\bar{\BD}].
\label{eq:lbound2}
\end{split}
\end{align}

Since both components and products are measured in discrete units, (\ref{eq:MSP1}) and (\ref{eq:MSP2}) should be a discrete optimization problems that can be difficult to solve exactly.  For the purpose of setting a lower bound, if suffices to solve (\ref{eq:MSP2}) without the integrality constraint, which is a continuous relaxation of the discrete problem. However, as will be shown in the next section, we will also use the SP to formulate inventory control policies, and for that purpose, the solutions need to be integer-valued, and also satisfy certain uniqueness and continuity conditions.  These issues will be addressed in Section \ref{sec:SPsolution}.

\section{Inventory Policy}
\label{sec:policy}
 
Below we will first describe the general idea that motivates our approach in Section \ref{subsec:GIdea}, followed by developments of our replenishment and allocation policies in Sections \ref{subsec:r-policy} and \ref{subsec:a-policy} respectively.  We will then make a few important observations in Section \ref{subsec:comments}.
% observations and highlight a related issue that will be addressed in Section \ref{sec:SPsolution}.
A simple example that illustrates the application of our inventory policy, involving the N system, is given in Appendix B.

\subsection{General Idea}
\label{subsec:GIdea}

To drive the average cost of the ATO system towards its SP-based lower bound in (\ref{eq:MSP2}), we formulate inventory policies that mimic the optimal solution of this SP.
As is commonly known and easily verifiable from (\ref{eq:flow})-(\ref{eq:IPdef}), in an ATO system, the net inventory levels satisfy
\begin{align*}
\begin{split}
\mathbf I^k(t) - A^k \mathbf B(t)&= \ip^k(t-L_k)- A^k \BD(t-L_k,t)
\\
&= \ip^k(t-L_k)- A^k [\BD(t-L_K,t)-\BD(t-L_K,t-L_k)],\quad 1 \le k \le K, \quad t \ge 0.
\end{split}
\end{align*}
The second equality allows us to rewrite the expected cost rate as the following function of inventory positions and backlog levels: 
\begin{align}
\begin{split}
\mathcal C(t) &=\sum_{k=1}^K \mathbf h^k \cdot \E[\mathbf I^k(t)]+\mathbf b \cdot \E[\mathbf B(t)]
\\
&=\sum_{k=1}^K \mathbf h^k \cdot \left(\E[\ip^k(t-L_k)]+A^k \E[\BD(t-L_K,t-L_k)]\right)+\mathbf c \cdot \E[\mathbf B(t)]
\\
&~~~~
-\E[\BD(t-L_K,t)]\sum_{k=1}^K (A^k)'\mathbf h^k , \quad t \ge 0,
\label{eq:costdcom}
\end{split}
\end{align}
where $\mathbf c$ is defined in (\ref{eq:defC}). 
Since $\bar{\BD} \stackrel{d}{=}\BD(t-L_K,t)$, $\mathbf E[\bar{\BD}]=\E[\BD(t-L_K,t)]$. Thus by a direct comparison of the cost rate in (\ref{eq:costdcom}) with its lower bound in (\ref{eq:lbound2}), 
at any time $t$ $(t \ge 0)$,
\begin{eqnarray}
\label{eq:pmatch}
\mbox{if}\quad
&~&\E[\ip^k(t-L_k)]+A^k \E[\BD(t-L_K,t-L_k)]=\E[\y^{k*}], \quad 1 \le k \le K,
\\
\label{eq:bmatch}
&~&\mbox{and} \quad\E[\mathbf B(t)]=\E[\mathbf B^*],
\\
\nonumber
\mbox{then} \quad &~&\mathcal C(t) =\underline{C}. 
\end{eqnarray}
Hence the average inventory cost of an ATO system $(\mathcal C)$ reaches its lower bound $(\underline{C})$ if (\ref{eq:pmatch})-(\ref{eq:bmatch}) are both satisfied for all time. Based on this observation, we develop replenishment and allocation policies by using the SP solutions to set targets for inventory positions and backlog levels, keeping actual values ``close'' to these targets under feasibility constraints.

\subsection{Replenishment Policy}
\label{subsec:r-policy}

\subsubsection{Overview:}
The RHS of (\ref{eq:pmatch}) is determined in (\ref{eq:MSP2}). For $k=K$, $\y^{K*}$ is constant, and for $k=K-1,\cdots, 1$, $\y^{k*}$ $(k=K-1,\cdots, 1)$ are determined recursively (going backwards in $k$)  on each sample path of $(\BD^K,\cdots,\BD^{k+1})$ as values of $\y^k$ that minimize
\begin{equation}
\label{eq:ykdef}
\mathbf h^k\y^k + \E\left [\varphi^{k-1}(\y^k,\y^{(k+1)*},\cdots,\y^{K*}, \x+\BD^{k})\right],
\end{equation}
where $\x$ denotes the value of $\BD^K+\cdots +\BD^{k+1}$ on that path.  

\begin{figure}[ht]
\begin{center}
\includegraphics[scale=0.6]{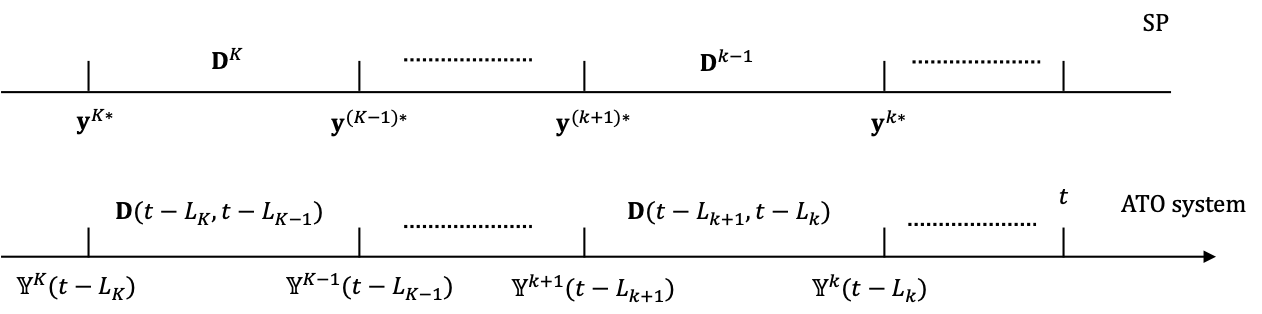}
\caption{Match between processes $\mathbf Y^k(t)$ in the ATO system with SP solutions $\y^{k*}$ $(t \ge -L_k, 1 \le k \le K)$}
\label{fig:match}
\end{center}
\end{figure}

Our replenishment policy is built upon processes $\mathbb Y^k(t)$ $(t \ge -L_k)$, with their distributions at each point of time mimicking those of $\y^{k*}$ $(1 \le k \le K)$. Figure \ref{fig:match} shows a match of these two quantities in any period $[t-L_K,t]$ ($t \ge 0$). Let $\mathbb Y^{K}(t-L_K)=\y^{K*}$, and determine $\mathbb Y^k(t-L_k)$ $(k=K-1,\cdots 1)$ recursively on each sample path of $\BD(t-L_K,t-L_k)$ to minimize (\ref{eq:ykdef}) with $\x=\BD(t-L_K,t-L_k)$. Since $\BD(t-L_k,t-L_{k-1})$ and $\BD^k$ $(1 \le k \le K)$ have the same distribution,
\begin{equation}
\label{eq:pmatch-1}
\E[\mathbb Y^k(t-L_k)]=\E[\y^{k*}],\quad 1 \le k \le K. 
\end{equation} 
Motivated by (\ref{eq:pmatch}), our policy uses $\mathbb Y^k(t-L_k)$ and $A^k \BD(t-L_K,t-L_k)$ to set {\it inventory position targets} $\tIP^k(t)$  $(t \ge -L_k, 1\le k \le K)$ and make ordering decision to move actual inventory positions $\ip^k(t-L_k)$ $(t \ge -L_k, 1 \le k \le K)$ towards these targets.

%For components with the longest lead time $L_K$, targets are constant and our design leads to a base stock policy to keep inventory positions at these constant levels. For components with shorter lead times,
Except for components with the longest lead time $L_K$, inventory position targets typically change over time. When a target rises above the actual inventory position, it is feasible to bring the latter position to its target immediately by ordering an appropriate amount of the component.  When a target falls below the actual position, it is not feasible to close the gap immediately before the target is reduced and/or new demand arrives. Recognizing the latter restriction, our replenishment policy orders a component when and only when its inventory position is below the target, in the exact amount needed to eliminate the deficit.

\subsubsection{Specific Policy Procedure:}
\label{subsubsec:rpolicydef} For components with the longest lead time $L_K$, the inventory position target process starts from time $-L_K$ and is determined by 
\begin{equation}
\label{eq:YK}
\tIP^K(t)= \mathbb Y^K, \quad t \ge -L_K,
\end{equation}
where $\mathbb Y^K$ is the value of $\y^K$ that minimizes 
\begin{equation}
\mathbf h^K \cdot \y^K+\E\left [\varphi^{K-1}(\y^K,\BD^K)\right].
\label{eq:targetK}
\end{equation}
%and we take the ceilings of its values to keep inventory position targets as integers.  
For components with lead time $L_k$ $(k=K-1,\cdots,1)$, the target process starts from $-L_k$ and that its values are given by 
\begin{equation}
\label{eq:target}
\tIP^{k}(t)=\mathbb Y^{k}(t)-A^{k} \BD(t+L_k-L_K,t),\quad t \ge -L_k, \quad 1 \le k \le K-1, 
\end{equation}
where referring to (\ref{eq:ykdef}) and the match shown in Figure \ref{fig:match}, 
$\mathbb Y^k(t)$ is the value of $\y^k$ that minimizes 
\begin{equation}
\label{eq:Yset}
\mathbf h^{k}\cdot\y^{k} + \E\left [\varphi^{k-1}(\y^{k},\mathbb Y^{k+1}(t+L_{k}-L_{k+1}), \cdots, 
\mathbb Y^K,\BD(t+L_{k}-L_K,t)+\BD^{k})\right],\quad t \ge -L_k.
\end{equation}

Starting from time $-L_k$, the actual inventory position of a component with lead time $L_k$ is compared with its target. New replenishment is ordered to keep the actual inventory position at 
\begin{equation}
\label{eq:ipcontrol}
\ip^k(t) = \lceil \tIP^k(t)\rceil \vee \ip^{k-}(t), \quad \quad t \ge -L_k, 1 \le k \le K,
\end{equation}
where recognizing that inventory positions need to be integers, the ceiling function is applied to $\tIP^k(t)$. 
The ``pre-order'' inventory position $\ip^{k-}(t)$ $(t \ge -L_k, 1 \le k \le K)$ is defined in  Table \ref{tbl:notations}.
% the inventory position at time $t$ before placing an order, and its value is determined by the definition in Table \ref{tbl:notations}.

While $\tIP^k(t)$ and $\ip^k(t)$ $(t \ge -L_k, 1 \le k \le K)$ are both continuous-time processes, their values change only at discrete points of time,
corresponding to new demand arriving at $t$, or a demand arrival at one of a particular set of previous times, as follows.
For components with lead time $L_K$, as seen in (\ref{eq:YK}),
the procedure reduces to a constant base stock policy with $\mathbb Y^K$ as base stock levels. For components with lead time $L_{K-1}$, $\tIP^{K-1}(t)$ needs to be updated only when there is a demand arrival at time $t+L_{K-1}-L_K$ or $t$, which changes the value of $\BD(t+L_{K-1}-L_K,t)$, so (\ref{eq:Yset}) needs to be re-solved to update $\mathbb Y^{K-1}(t)$ $(t \ge -L_{K-1}$). 
Similarly, by induction on $k$ in (\ref{eq:Yset}), for $k<K-1$,  $\tIP^k(t)$ can change only when there is a demand arrival at time i) $t$, which changes $\BD(t+L_k-L_K,t)$, or ii) $t+L_k-L_{k'}$ ($k<k'<K$) which might change $\mathbb Y^{k'}(t+L_k-L_{k'})$, or iii) $t+L_k-L_K$, which changes both $\BD(t+L_k-L_K,t)$ and $\mathbb Y^{k'}(t+L_k-L_{k'})$ ($k<k'<K$). Moreover, (\ref{eq:ipcontrol}) implies that only at these times can $\ip^k(t)$ change: because $\tIP^k(t)$ can only change at these times, and as Table \ref{tbl:notations} shows, $\ip^{k-}(t)$ can differ from $\ip^k(t^-)$ only when there is demand arrival at time $t$ $(t \ge -L_k, 1 \le k \le K$). 

Although it may seem that an exact implementation of the above idea would require looking back in time at every moment $t$, to see if there is any demand arrival at times $t+L_k-L_{k'}$ $(k<k'<K)$, to decide whether to update $\tIP^k(t)$ and place orders to reach new $\ip^k(t)$ $(t \ge -L_k, 1 \le k \le K$),
there is an equivalent simpler 
and more intuitive implementation, which allows 
demand arrivals to drive future updates. Specifically, when there is a demand arrival at time $t$, schedule updates of inventory position targets and corresponding ordering decisions at times $t, t+L_{k+1}-L_k,\cdots, t+L_K-L_k$ $(1 \le k \le K-1)$. 

\begin{figure}[ht]
\begin{center}
\includegraphics[scale=0.5]{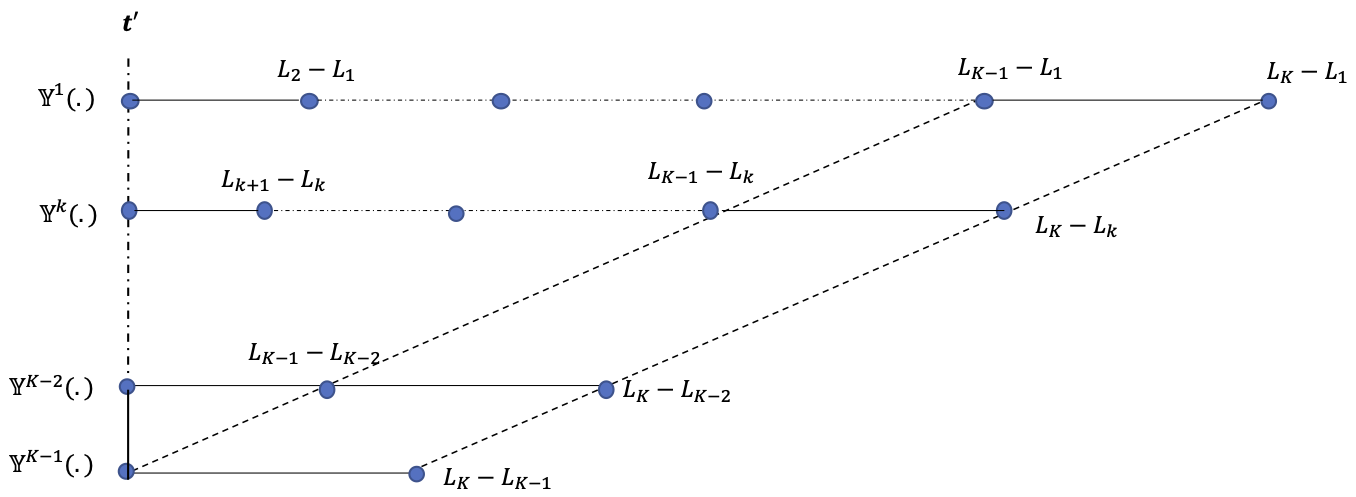}
\caption{Inv. Position Target Updates and Ordering Decisions that Follow Demand Arrival at Time $t$}
\label{fig:updateschdule}
\end{center}
\end{figure}

Figure \ref{fig:updateschdule} illustrates the process. The vertical line represents a particular time $t'$ when there is a demand arrival. Each horizontal line corresponds to a lead time. On the line corresponding to $L_k$, circles mark points of time when (\ref{eq:Yset}) needs to be re-solved to update $\mathbb Y^k(\cdot)$ and $\tIP^k(\cdot)$ $(1 \le k<K$). The expression next to a circle shows the distance in time from time $t'$ to the circle.  Updates of $\mathbb Y^k(\cdot)$ take place at times $t'$ and $t'+L_K-L_k$ as demand arrival at time $t'$ changes $\BD(t+L_k-L_K,t)$ in (\ref{eq:Yset}) when $t=t'$ or $t=t'+L_K-L_k$ $(1 \le k <K)$. Moreover, as shown by two parallel dashed lines, updates of $\mathbb Y^{K-1}(\cdot)$ at times $t'$ and $t'+L_K-L_{K-1}$ change values of $\mathbb Y^{K-1}(t+L_k-L_{K-1})$ in (\ref{eq:Yset}) when $t=t'+L_{K-1}-L_k$ and $t=t'+L_K-L_k$, so $\mathbb Y^k(\cdot)$ needs to be updated at the latter two times $(1 \le k <K-1$). In general updates of $\mathbb Y^k(\cdot)$ at times $t'$ and $t'+L_K-L_k$ trigger updates of $\mathbb Y^l(\cdot)$ at times $t' + L_k-L_l$ and $t'+L_K-L_l$ respectively $(1 \le l <k$).

\begin{table}[ht]
\caption{Replenishment Policy Procedure}
\label{tbl:replenishment}
\begin{tcolorbox}

%\begin{center}{\bf Replenishment Policy Procedure} \end{center}

Initialization: for $k=1,\cdots K$, let
\[
\Gamma^k=\{-L_k\}, \quad \mathcal R^k(-L_k^-)=0, \quad \mbox{and} \quad \ip^{k}(-L_k)=-A^k \mathcal D(-L_k).
\]

For $t \ge -L_K$:

\medskip

\begin{enumerate}
\item
if $\mathbf d(t)>0$, then for $k=1,\cdots,K$ such that $t \ge -L_k$, let

%\smallskip\begin{enumerate}\item

\begin{align*}
\begin{split}
\ip^{k-}(t) & =\ip^k(t^-)-A^k \mathbf d(t).
\\
\mbox{and} \qquad \qquad
\Gamma^k &= \Gamma^k \cup \{t, t+L_{k+1}-L_k,\cdots, t+L_K-L_k\}.
\end{split}
\end{align*}

\item
for $k=1,\cdots, K$, if $t \in \Gamma^k$, then:

\smallskip

\begin{enumerate}
\item
let $\mathbb Y^k(t)$ be the value of $\y^k$ that minimizes (\ref{eq:Yset})  

\smallskip
\item
as in (\ref{eq:target}), let 
\[
\tIP^{k}(t)=\mathbb Y^{k}(t)-A^{k} \BD(t+L_k-L_K,t).
\]
\item
for each component $j$ with lead time $L_k$, let
\[
\mathcal R_j(t)=\mathcal R_j(t^-)+\left(\lceil \tIP_j(t)\rceil-IP_j^-(t)\right)^+,
\]
by ordering $\left(\lceil \tIP_j(t)\rceil-IP_j^-(t)\right)^+$ units, changing its inventory position to
\[
IP_j(t)=IP_j^-(t)+\mathcal R_j(t)-\mathcal R_j(t^-).
\]
\end{enumerate}
\end{enumerate}

%\small{(note that for all $j=1,\cdots n$, given $IP_j^-(-L_{k_j})$ and $\mathcal R_j(t)$ $(t \ge -L_{k_j})$, values of $IP_j(t)$ $(t \ge -L_{k_j})$ and $IP^-_j(t)$ $(t>-L_{k_j})$ can be determined according to their definitions in Table \ref{tbl:notations}).}

\end{tcolorbox}
\end{table}

The complete replenishment procedure is presented in Table \ref{tbl:replenishment}. 
The collection of times at which $\mathbb Y^k(\cdot)$, $\tIP^k(\cdot)$ and $\ip^k(\cdot)$ possibly change is denoted by $\Gamma^k (1 \le k \le K)$.
The ceiling function on $\tIP_j(t)$  in step 2(c) ensures that both $IP_j(t)$ and $\mathcal R_j(t)$ are integer-valued ($t \ge -L_{k_j},1 \le j \le n$).  In Steps 2(a)-(b), inventory position targets $\tIP^k(t)$ $(t \ge -L_k,1 \le k \le K$) are updated repeatedly over time except for $k=K$, in which case $\tIP^K(t)$ are constants that are set once at time $-L_K$. 
%, 2(a) and 2(b) in Table \ref{tbl:replenishment} yield the same values of $\mathbb Y^K(t)$ and $\tIP^K(t)$ each time. Thus the two steps only need to be performed once at time $-L_K$. 
It is easy to verify that $\mathcal R_j(t)$ $(t \ge -L_{k_j}, 1 \le j \le n)$ prescribed in the table satisfy the feasibility requirements defined in Section \ref{sec:formulation}: integer-valued, non-negative, non-decreasing, right-continuous, and non-anticipating. 

It is worth pointing out that in ATO systems with a single product, this replenishment policy specializes to the one prescribed by  \cite{Rosling1989}.  In this special case, the latter policy is exactly optimal because all components are used according to fixed proportions, so one can always keep inventory positions at their targets (after some initial lapse for these positions to reach ``coordinated levels''). In systems with identical lead times, our policy specializes to the one formulated in \cite{Reiman2015}, which uses base stock policies for replenishment and is asymptotically optimal in the large lead time regime. 

\subsection{Allocation Policy}
\label{subsec:a-policy}

\subsubsection{Overview}

On the RHS of (\ref{eq:bmatch}), $\B^*$ is the optimal solution of the LP
\begin{equation}
\label{eq:miniB}
\min_{\B} \left \{\mathbf c \B~ \big |~ \B \ge 0,~ A^k \B \le A^k \x - \y^{k*},1 \le k \le K\right \},
\end{equation}
which is defined on each sample path of $(\BD^K,\cdots,\BD^1)$, with $\x$ denoting $\BD^K+\cdots +\BD^1$ and $\y^{k*}$ determined recursively by minimizing (\ref{eq:ykdef}) for values of $(\BD^K,\cdots,\BD^{k+1})$ ($1 \le k \le K$) on the same path. 
In the corresponding ATO system, $\BD(t-L_K,t)$ has the same distribution as $\BD^K+\cdots+\BD^1$ and $\mathbb Y^k(t-L_k)$ has the same distributions as $\y^{k*}$ $(t \ge 0, 1 \le k \le K)$. Hence we can mimic $\B^*$ at time $t$ $(t \ge 0)$ by letting $\bB^*(t)$ be the value of $\B$ that minimizes (\ref{eq:miniB}) with  $A^k \x - \y^{k*}$ replaced by
\begin{equation}
\label{eq:iQ}
\mathbb Q^k(t):= A^k \BD(t-L_K,t) - \mathbb Y^k(t-L_k).
\end{equation}
As a result, (\ref{eq:bmatch}) holds for $\bB^*(t)$, i.e.,
\begin{equation}
\label{eq:bmatch-1}
\E[\bB^*(t)]=\E[\mathbf B^*].
\end{equation}

While $\bB^*(t)$ $(t \ge 0)$ gives the desired backlog levels for the inventory cost to reach its lower bound, it is determined by 
\begin{align*}
\begin{split}
\mathbb Q^k(t) &= A^k \BD(t-L_K,t)-(\tIP^k(t-L_k)+A^k \BD(t-L_K,t-L_k))
\\
&=A^k \BD(t-L_k,t) - \tIP^k(t-L_k), \qquad 1 \le k \le K.
 \end{split}
 \end{align*}
Comparing with $\mathbf Q^k(t)$ in Table \ref{tbl:notations}, $\mathbb Q^k(t)$ differs from actual component balances by having inventory position targets $\tIP^k(t-L_k)$ in place of actual positions $\ip^k(t-L_k)$ $(t \ge 0, 1 \le k \le K)$. Since our replenishment policy in general does not keep inventory positions at their targets, $\mathbb Q^k(t)$ does not always capture the exact state of component availability in the system.  
Recognizing this discrepancy, we set the backlog target at $\bB(t)$ $(t \ge 0)$, which, like $\bB^*(t)$, is determined by minimizing (\ref{eq:bmatch}), but using $\mathbf Q^k(t)$ instead of $\mathbb Q^k(t)$ to replace $A^k \x-\y^{k*}$ $(1 \le k \le K)$. 
While (\ref{eq:bmatch}) may not hold for  $\bB(t)$ $(t \ge 0)$ at equality, our analysis below  shows that the difference is asymptotically negligible. 

Actual backlog levels may exceed or fall below their targets. It is generally infeasible to keep all backlog levels at their targets: a below-target backlog level cannot be raised immediately to the target if there is no new demand arrival, and an above-target backlog level cannot be reduced if a required component is not available. Naturally, we prescribe an allocation policy that uses all available components to reduce backlogs that are above their targets and never serves a demand when its backlog level is at or below the target. 

\subsubsection{Specific Policy Procedures:} 
\label{subsubsec:apolicydef}
Table \ref{tbl:allproc} presents the specific policy procedure that implements the above idea on component allocation. 

\begin{table}[ht]
\caption{Allocation Policy Procedure}
\label{tbl:allproc}
\begin{tcolorbox}
%\begin{center}\textbf{Backlog Target Procedure}\end{center}
Initialization: let
\[
\mathcal Z(t)=0 \quad \mbox{for all}~t <0 \quad (\mbox{so}~\B(0^-)=\mD(0^-)).
\]

For $t \ge 0$, 

\smallskip

\begin{itemize}
\item[]
if $t=0$ or $Q_j(t)-Q_j(t^-)\ne 0$ for some $j=1,\cdots,n$, then

\medskip

\begin{enumerate}
\item
let $\bB(t)$ be the value of $\x$ that minimizes:
\begin{equation}
\label{eq:setBtar}
\left \{\mathbf c \x~ \big |\x \ge 0,A^k \x \ge \mathbf Q^k(t),1 \le k \le K\right \}.
\end{equation} 
\item
let
\begin{equation}
\label{eq:Z}
\mathcal Z(t)=\mathcal D(t)-\B(t),
%\quad i=1,\cdots,m,
\end{equation}
where $\B(t)$ must be integer-valued and its choice satisfies:
\smallskip

\begin{enumerate}
\item 
for all $i=1,\cdots,m$,
\begin{equation}
\label{eq:Z-a}
\left[\lfloor B_i(t)-\mathbb B_i(t)\rfloor\right]^+ \wedge \min_{j: a_{ji}>0} \left \{\left(I_j(t)-a_{ji}+1\right)^+\right\}=0,
\end{equation}

and
\begin{equation}
\label{eq:Z-b}
B^-_i(t) \wedge \mathbb B_i(t) \le B_i(t) \le B_i^-(t),
\end{equation}

where
\begin{equation}
\label{eq:Z-}
B_i^-(t)=B_i(t^-)+d_i(t),
\end{equation}

and by definition (\ref{eq:balancedef}), 
\[
I_j(t)=\sum_{i'=1}^m a_{ji'} B_{i'}(t)-Q_j(t).
\]
\item
for all $j=1,\cdots,n$,
\begin{equation}
\label{eq:Z-c}
\sum_{i=1}^m a_{ji} B_i(t)  \ge  Q_j (t).
\end{equation}
\end{enumerate}
\end{enumerate}
\end{itemize}
\end{tcolorbox}
\end{table}

To explain, the allocation decision starts at time $0$ when the first batches of replenishment orders arrive. Before that time, backlogs simply accumulate as new demands arrive. After that time, new allocations are triggered by changes of $Q_j$ ($1 \le  j \le n)$, the balance of some components. (This occurs as demands and replenishments arrive.) As Table \ref{tbl:notations} shows, $\mathbf Q^k(t)$ is the difference between the accumulated demand and supply,  $A^k \mathcal D(t)-\mR^k(t-L_k)$ $(t\ge 0,1 \le k \le K)$. The former is a compound Poisson process and the latter is a pure jump process under the replenishment policy in Table \ref{tbl:replenishment}. Thus allocation actions are taken at discrete points of time, which include using (\ref{eq:setBtar}) to set backlog targets and serving demand under conditions (\ref{eq:Z})-(\ref{eq:Z-c}). Of the latter conditions, (\ref{eq:Z}) and (\ref{eq:Z-}) simply repeat definitions of $\B^-(t)$ and $\B(t)$ in Table \ref{tbl:notations}. Distinct features of the policy are characterized by (\ref{eq:Z-a}), (\ref{eq:Z-b}), and (\ref{eq:Z-c}). 
\begin{itemize}
\item
In (\ref{eq:Z-a}), $B_i(t)$ is specified as an integer, $\bB_i(t)$ is real number, and the floor function is applied to eliminate the rounding error in their difference. The equation defines the key feature of the policy: a backlog level can exceed its target (within the rounding error) only when the system runs out of a needed component to serve the demand, i.e., when
\[
a_{ji} > \sum_{i'=1}^m a_{ji'} B_{i'}(t)-Q_j(t)=I_j(t)\quad \mbox{for some}~j~(1 \le j \le n). 
\]
\item
The left inequality of (\ref{eq:Z-b}) requires that the backlog level of a demand ($B_i(t)$) should be kept at its existing level ($B^-_i(t))$ if the latter level does not exceed its target ($\bB_i(t)$).   
\item
The right inequality of (\ref{eq:Z-b}) encodes that serving a demand  can only reduces its backlog level. 
\item 
Referring to (\ref{eq:balancedef}), (\ref{eq:Z-c}) is equivalent to  
\[
I_j(t) \ge 0,\quad 1 \le j \le n,
\]
i.e., demands cannot be served with a non-existing component.
\end{itemize}

In prescribing this procedure, we have omitted specific processes for determining $\B(t)$ $(t \ge 0)$ to satisfy (\ref{eq:Z-a}), (\ref{eq:Z-b}), and (\ref{eq:Z-c}) because there can be many different ones. For instance, one can select all backlogs that exceed their targets and use available components to reduce them one at a time according to a particular order. Each backlog is reduced to the point where (\ref{eq:Z-a}) applies, i.e., either the backlog reaches the target or there are not enough components to bring it down further. This process obviously also satisfies (\ref{eq:Z-b}) and (\ref{eq:Z-c}). We can formulate many other processes by, for instance, changing the order by which backlogs are processed or not following a fixed order at all. 
%Our procedure in Table \ref{tbl:allproc} implicitly defines our family of allocation policies by the requirement that they satisfy $\B(t)$ $(t \ge 0)$ satisfies (\ref{eq:Z-a}),  (\ref{eq:Z-b}), and (\ref{eq:Z-c}). 
Therefore Table \ref{tbl:allproc} defines not a single policy but a family of eligible policies. 
Such policies always exist, for example, the one that follows the process we just described.

To show that this family of eligible policies satisfy requirements for feasible policies defined in Section \ref{sec:formulation}: notice that once $\B(t)$ $(t \ge 0)$ is determined, the original allocation process $\mathcal Z(t)$ $(t \ge 0)$ is fully specified by (\ref{eq:Z}). By (\ref{eq:Z}) and (\ref{eq:Z-}), we can write 
\[
\mathcal Z(t)= \mathcal Z(t^-)+\B^-(t)- \B(t),\quad t \ge 0,
\]
and use (\ref{eq:Z-a}) to show that $\mathcal Z(t)$ $(t \ge 0)$ is non-decreasing, and also non-negative as $\mathcal Z(t)=0$ for $t<0$. Since the value of $\mathcal Z(t)$ $(t \ge 0)$ changes only at discrete points of time, the process is right-continuous. Apply the definition of $\mathbf Q^k(t)$ in Table \ref{tbl:notations}, its alternative expression in (\ref{eq:balancedef}), and use (\ref{eq:Z}) to write
\[
A^k \mathcal Z(t) = A^k \mathcal D(t) - A^k \B(t)=A^k \mathcal D(t)-(\mathbf Q^k(t)+\mathbf I^k(t))=\mathcal R^k(t-L_k)-\mathbf I^k(t),~1 \le k \le K, t\ge 0,
\]
from which it is easy to see that $\mathcal Z(t)$ $(t \ge 0)$ satisfies all conditions in (\ref{eq:feasibility}). 
 
%\[\mathcal Z(t) \le \mathcal D(t),~(t \ge -L_K), \qquad \mbox{and} \qquad  A^k \mathcal Z(t) \le \mathcal R^k(t-L_k)~(t\ge 0,~ 1 \le k \le K).\]

From the first expression of $\mathbf Q^k(t)$ $(1 \le k \le K, t \ge 0)$ in (\ref{eq:balancedef}) and the constraints in the LP (\ref{eq:setBtar}),
\[
I_j(t)= \sum_{i=1}^m a_{ji} B_i(t)-Q_j(t) \ge \sum_{i=1}^m a_{ji} \left(B_i(t)-\bB_i(t)\right),\quad 1 \le j \le n,~t \ge 0.
%\mathbf I^k(t) = A^k \B(t)-\mathbf Q^k(t) \ge A^k \B(t)-A^k \bB(t)= A^k (\B(t)-\bB(t)),\quad k=1,\cdots K,
\]
so for any given product $i$ $(1 \le i \le m)$,
\[
I_j(t)+ \sum_{i' \ne i} a_{ji'} [\bB_{i'}(t)-B_{i'}(t)]  \ge a_{ji} [B_i(t)-\bB_i(t)], \quad 1 \le j \le n,~t \ge 0.
\]
Apply the inequality to (\ref{eq:Z-a}) leads to the following important property of our allocation policy:

\smallskip
For any $i=1,\cdots,m$, if $B_i(t)-\bB_i(t) \ge 1$, then 
\begin{equation}
\label{eq:twoside}
1+\frac{\bar{a}}{\underline{a}} \sum_{i' \ne i} (\bB_{i'}(t)-B_{i'}(t))^+ \ge B_i(t)-\bB_i(t),\quad t \ge 0.
\end{equation}
(see Table \ref{tbl:minmaxpara} for definitions of $\bar{a}$ and $\underline{a}$)).

\smallskip

By this property, our policy does not allow backlog of any product $i$ to exceed its target by more than a rounding error when no backlog is below its target (i.e., $\bB_{i'}(t) \le B_{i'}(t)$ for all $i'\ne i$). A similar observation was made in \cite{Reiman2015} for developing asymptotically-optimal policies to manage ATO systems with identical lead times. Likewise, we will use this fact in Section \ref{subsubsec:blkres}.

\subsection{Further Comments}
\label{subsec:comments}

Tables \ref{tbl:replenishment} and \ref{tbl:allproc} show that our policies can be fully implemented by controlling inventory positions via (\ref{eq:ipcontrol}) and choosing backlog levels that satisfy (\ref{eq:Z-a}), (\ref{eq:Z-b}), and (\ref{eq:Z-c}). For brevity, from now on, we will only use the latter variables to characterize our policies, instead of invoking the original replenishment process $\mathcal R^k(t)$ $(t \ge -L_k,1 \le k \le K)$ and allocation process $\mathcal Z(t)$ $(t \ge -L_K)$.

The procedures in these tables implicitly assume that there is no ambiguity about $\mathbb Y^k(t)$ $(t \ge -L_k, 1 \le k \le K)$ and $\bB(t)$ $(t \ge 0)$, which is true if the related SPs for choosing these values all have unique optimal solutions.  In Section \ref{sec:SPsolution}, we will discuss how to perturb these SPs to satisfy the uniqueness condition without compromising  asymptotic optimality.  

\section{Asymptotic Optimality}
\label{sec:asymopt}

The rest of the paper will focus on performance evaluation, especially asymptotic optimality, of our policy, and this section provides a critical lead in to this analysis.  We first set up the asymptotic regime in Section \ref{subsec:asymregime}. We then define the asymptotic optimality criteria in Section \ref{subsec:asymcriteria}, followed by the development of a key theorem in Section \ref{subsec:keytheorem} (Theorem \ref{thm:final}) that specifies sufficient conditions for meeting these criteria.

\subsection{Large Lead Time Asymptotic Regime}
\label{subsec:asymregime}

We introduce a family of systems indexed by $L$. 
In the $L^{th}$ system, $L_k^{(L)}$ is the lead time of components $\bar{n}_{k-1}+1,\cdots,\bar{n}_k$ ($1 \le k \le K$).   The longest lead time $\LL_K=L$. We define the large lead time asymptotic regime by letting the longest lead time $L \rightarrow \infty$.

For convenience, define $\LL_0=0$. Let
\begin{equation}
\label{eq:scalelt}
\hLL_k \stackrel{\Delta}{=}\frac{\LL_k}{L}, \quad 0 \le k \le K. 
\end{equation}
We impose no assumptions on other lead times except that for all $L$, 
\[
\LL_0<\LL_1<\LL_2<\cdots<\LL_K, \quad  \mbox{and therefore} \quad 0< \hLL_1 < \cdots \hLL_{K-1} < 1.
\]

\subsubsection{Demand Process}
\label{subsubsec:demand}
The $L^{th}$ system is empty at time $-L$ when the demand process $\sD(t)$ ($t \ge -L$) starts. This process is  the same as $\mathcal D(t)$ $(t \ge -L_K)$ in Section \ref{sec:formulation} except for the starting time.
Similar to definitions of $\BD(t_1,t_2)$ and $\mathbf d(t)$ in Table \ref{tbl:notations}, let $\sBD(t_1,t_2)$ be the increment of $\sD(t)$ over $(t_1,t_2]$ $(-L \le t_1<t_2)$, with $\boldsymbol \mu$ and $\boldsymbol \Sigma$ denoting the mean and covariance matrix of $\sBD(t,t+1)$, and $\sjD(t)$ be the instantaneous change of $\sD(t)$ at $t$ $(t \ge -L)$. Define
\begin{align}
\begin{split}
\label{eq:scaledemand}
&\hBD(t_1,t_2) \stackrel{\Delta}{=}\frac{\sBD(Lt_1,Lt_2)-L(t_2-t_1){\boldsymbol \mu}}{\sqrt{L}}, \quad -1 \le t_1 < t_2,
\\
\mbox{and} \quad
&\hjD(t)\stackrel{\Delta}{=}\frac{\sjD(Lt)}{\sqrt{L}}, \quad t \ge -1.
\end{split}
\end{align}
Let $\BD^{k(L)}$ be a random vector with the same distribution as $\sBD(0,\LL_k-\LL_{k-1})$ $(1 \le k \le K)$. 
Then
\[
\hat{\BD}^{k(L)} \stackrel{\Delta}{=} \frac{\BD^{k(L)}-(\LL_k-\LL_{k-1}) \boldsymbol \mu}{\sqrt{L}}
\]
has the same distribution as $\hBD(t-\hLL_{k},t-\hLL_{k-1})$ $(t \ge -1, 1 \le k \le K)$. 
%Consistent with the notation in Table \ref{tbl:notations}, let $\bar{\BD}^{(L)}$ be a random variable that has the same distribution as $\BD^{(L)}(t-L_K,t)$ $(t \ge 0)$. 

\subsubsection{Replenishment and Inventory}
In the $L^{th}$ system, 
replenishments of components with lead time $\LL_k$ start at time $-\LL_k$ ($1 \le k \le K$). Following the policy description in Section \ref{subsec:r-policy}, 
the inventory position targets are
\[
\tIP^{k(L)}(t)=\mathbb Y^{k(L)}(t)-A^k\BD^{(L)}(t+\LL_k-L,t),\quad t \ge -\LL_k, 1\le k \le K,
\]
where $\mathbb Y^{K(L)}(t)$ is a constant vector (and thus will be denoted by $\mathbb Y^{K(L)}$ below) that minimizes
\begin{equation}
\label{eq:SPL}
\mathbf h^K\y^K + \E\left [\varphi^{K-1}(\y^K,\BD^{K(L)})\right],
\end{equation}
and $\mathbb Y^{k(L)}(t)$ $(1 \le k <K)$  minimizes
\begin{equation}
\tag{\ref{eq:SPL}'}
\mathbf h^k\y^k + \E\left [\varphi^{k-1}(\y^k,\mathbb Y^{k+1(L)}(t+\LL_k-\LL_{k+1}),\cdots,\mathbb Y^{K(L)}, \x+\BD^{k(L)})\right]
\end{equation}
on the sample path of 
\[
\sBD(t+\LL_k-\LL_{k'},t+\LL_k-\LL_{k'-1}), \quad k'=K,\cdots, k+1,
\]
with $\x$ denoting the value of $\sBD(t+\LL_k-L,t)$ on that path.

Let $\ip^{k(L)}(t)$ $(t \ge 0, 1 \le k \le K)$ be the actual inventory position in the $L^{th}$ system. Under the aforementioned replenishment policy,  
for $k=1,\cdots, K$:
\begin{align*}
\begin{split}
 & \ip^{k(L)-}(-\LL_k)=-A^k \mathcal D^{(L)}(-\LL_k), \\
& \ip^{k(L)}(t)=\tIP^{k(L)}(t) \vee \ip^{k(L)-}(t), \quad \qquad t \ge -\LL_k, \\
\mbox{and}  \quad &
\ip^{k(L)-}(t)=\ip^{k(L)}(t^-)-A^k \sjD(t), \qquad t > -\LL_k, 
\end{split}
\end{align*}
where the operator $\vee$ is applied componentwise. Component balances under the targeted and actual inventory positions are, respectively,
\begin{align}
\begin{split}
\label{eq:defQj}
\mathbb {Q}^{k(L)}(t)&=A^k \sBD(t-\LL_k,t)-\tIP^{k(L)}(t-\LL_k),
\\
\mbox{and} \quad
\mathbf{Q}^{k(L)}(t)&= A^k \sBD(t-\LL_k,t)-\mathbf{IP}^{k(L)}(t-\LL_k), \qquad t \ge 0, \quad 1 \le k \le K. 
\end{split}
\end{align}

For $k=1,\cdots, K$, let 
\begin{align}
\label{eq:scaleinv}
\begin{split}
%\hat{\mathbb Y}^{k(L)}(t) \stackrel{\Delta}{=} \frac{\mathbb Y^{k(L)}(Lt)}{\sqrt{L}}, \qquad &
 \hat{\mathbb Y}^{k(L)}(t) \stackrel{\Delta}{=} \frac{\mathbb Y^{k(L)}(Lt)-L A^k {\boldsymbol \mu}}{\sqrt{L}},  & \qquad \quad
t \ge -\hLL_k, \\
\widehat{\mathbb{IP}}^{k(L)}(t) \stackrel{\Delta}{=} \frac{\tIP^{k(L)}(Lt)-A^k \LL_k \boldsymbol \mu}{\sqrt{L}},
\qquad &\widehat{\mathbf{IP}}^{k(L)}(t) \stackrel{\Delta}{=} \frac{\mathbf{IP}^{k(L)}(Lt)-A^k\LL_k \boldsymbol \mu}{\sqrt{L}},\qquad t \ge -\hLL_k,
\\
\hat{\mathbb Q}^{k(L)}(t) \stackrel{\Delta}{=} \frac{\mathbb Q^{k(L)}(Lt)}{\sqrt{L}}, \qquad &\hat{\mathbf Q}^{k(L)}(t) \stackrel{\Delta}{=} \frac{\mathbf Q^{k(L)}(t)}{\sqrt{L}},\qquad \qquad  \qquad \qquad t \ge 0. 
\end{split}
\end{align}

%Consistent with previous notation, denote components of vectors $\mathbb Y^{k(L)}(t)$, $\hat{\mathbb Y}^{k(L)}(t)$,  $\widehat{\mathbb{IP}}^{k(L)}(t)$ and $\widehat{\mathbf{IP}}^{k(L)}(t)$ by $\mathbb Y_j^{k(L)}(t)$, $\hat{\mathbb Y}_j^{k(L)}(t)$,  $\thIP_j(t)$ and $\hIP_j(t)$ respectively ($\bar{n}_{k-1}+1 \le j \le \bar{n}_k$, $t \ge -\hLL_k$, $1 \le k \le K$). 

\subsubsection{Allocation and Backlog}
Referring to our allocation policy in Table \ref{tbl:allproc}, denote backlog targets $\bB(t)$  in the $L^{th}$ system by  $\sbB(t)$ ($t \ge 0$). Denote backlog levels $\mathbf B^-(t)$ and $\mathbf B(t)$ by  $\mathbf B^{(L)-}(t)$ and $\mathbf B^{(L)}(t)$ ($t \ge -L$) respectively. Hence $\sbB(t)$ $(t \ge 0)$ is the optimal solution that minimizes 
\begin{equation}
\label{eq:allspL}
\min_{\x} \left \{\mathbf c \x \large | \x \ge 0, A^k \x \ge \Q^{k(L)}(t), 1 \le k \le K\right\}.
\end{equation}
%Let $\mathbf B^{(L)-}(t)$ and $\mathbf B^{(L)}(t)$ $(t \ge -L)$ be backlog levels under our allocation policy in the $L^{th}$ system. There process are defined the same as $\mathbf B^-(t)$ and $\mathbf B(t)$ $(t \ge -L)$ and satisfies (\ref{eq:Z})-(\ref{eq:Z-c}), with $\mathbf B^{(L)}(t)$, $\sbB(t)$, and $\ corresponding to $\mathbf B(t)$ and $\bB(t)$ respectively. 
%We also define  that correspond to $\mathbf B^-(t)$.   
Corresponding to $\bB^*(t)$ prescribed by (\ref{eq:miniB})-(\ref{eq:iQ}), let $\sbbB(t)$ be the optimal solution to (\ref{eq:allspL}) with $\Q^{k(L)}(t)$ replaced by $\mathbb Q^{k(L)}(t)$ $(t \ge 0, 1 \le k \le K)$. 

Define
\begin{align}
\label{eq:scaleblk}
\begin{split}
\hbsB(t) \stackrel{\Delta}{=} \frac{\sbbB(Lt)}{\sqrt{L}}, \qquad & \hbB(t)\stackrel{\Delta}{=}\frac{\sbB(Lt)}{\sqrt{L}}, \qquad \quad  t \ge 0,\\
\hat{\mathbf B}^{(L)}(t) \stackrel{\Delta}{=} \frac{\mathbf B^{(L)}(Lt)}{\sqrt{L}}, \qquad & \hat{\mathbf B}^{(L)-}(t) \stackrel{\Delta}{=} \frac{\mathbf B^{(L)-}(Lt)}{\sqrt{L}}, \qquad t \ge -1.
\end{split}
\end{align}
Obviously, $\hbsB(t)$ and $\hbB(t)$ $(t \ge 0)$ are optimal solutions to (\ref{eq:allspL}) with $\Q^{k(L)}(t)$ replaced by $\hat{\mathbb Q}^{k(L)}(t)$ and $\hat{\Q}^{k(L)}(t)$ $(1 \le k \le K, t \ge 0)$ respectively. 

\subsubsection{Inventory Cost}

In the $L^{th}$ system, the long-run average expected inventory cost is 
\[
\mathcal C^{(L)}
= \limsup_{T \rightarrow \infty} \frac{1}{T} \int_{0}^{T} \mathcal C^{(L)}(t) dt. 
\]
Using (\ref{eq:costdcom}), the expected inventory cost rate can be written as
\begin{equation}
\label{eq:costrate}
\mathcal C^{(L)}(t)= \sum_{k=1}^K \mathbf h^k \cdot\left( \E[\mathbf{IP}^{k(L)}(t-\LL_k)]-A^k \E[\BD^{(L)}(t-\LL_k,t)] \right)+\mathbf c \cdot \E[\mathbf B^{(L)}(t)],\quad t \ge 0.
\end{equation}
Define
\[
\hat{\mathcal C}^{(L)}\stackrel{\Delta}{=} \frac{\mathcal C^{(L)}}{\sqrt{L}}.
\]
Then
\[
\hat{\mathcal C}^{(L)}=\limsup_{T \rightarrow \infty} \frac{1}{T} \int_{0}^{T} \hat{\mathcal C}^{(L)}(t) dt,
\]
where the scaled inventory cost rate is defined to be 
\[
\hat{\mathcal C}^{(L)}(t) \stackrel{\Delta}{=}\frac{\mathcal C^{(L)}(Lt)}{\sqrt{L}}, \quad t \ge 0,
\]
and by applying (\ref{eq:scaledemand}), (\ref{eq:scaleinv}), and (\ref{eq:scaleblk}) to (\ref{eq:costrate}), $\hat{\mathcal C}^{(L)}(t)$ can be written as 
%\begin{align}
%\begin{split}
%\hat{\mathcal C}^{(L)}(t)&=\frac{\sum_{k=1}^K \mathbf h^k \cdot\left( \E[\mathbf{IP}^{k(L)}(Lt-\LL_k)]-A^k \E[\BD^{(L)}(Lt-\LL_k,Lt)] \right)+\mathbf c \cdot \E[\mathbf B^{(L)}(Lt)]}{\sqrt{L}}\\
%&=\sum_{k=1}^k \mathbf h^k \cdot \left(\E\left[\widehat{\mathbf{IP}}^{k(L)}(t-\hLL_k)\right]-
%A^k \E\left[\hBD(t-\hLL_k,t)\right]\right)+\mathbf c \cdot \E\left[\hat{\mathbf B}^{(L)}(t)\right].
%\end{split}
%\label{eq:invcostscale}
%\end{align}
\begin{equation}
\label{eq:scalecostexp}
\hat{\mathcal C}^{(L)}(t)=\sum_{k=1}^K \mathbf h^k \cdot \left(\E\left[\widehat{\mathbf{IP}}^{k(L)}(t-\hLL_k)\right]-A^k \E\left[\hBD(t-\hLL_k,t)\right]\right)+\mathbf c \cdot \E\left[\hat{\mathbf B}^{(L)}(t)\right],\quad t \ge 0.
\end{equation}

\subsubsection{Summary of Definitions}
For future reference, we summarize the definitions made in this subsection for asymptotic analysis. Some vectors of random processe $\x^{(L)}(t)$ are centered and scaled component-wise by 
\[
\hat{\x}^{(L)} (t)= \frac{\x^{(L)}(Lt)-\bar{\x}^{(L)}}{\sqrt{L}},
\]
where $t$ varies in a proper range of time and the centering $\bar{x}^{(L)}$ is proportional to $L$. Similar operations apply to increments of demand processes over some periods $[t_1,t_2)$ and demand vectors used in the lower bound SP.  These definitions are given in Table \ref{tbl:C-S}.

\begin{table}[ht]
\renewcommand{\arraystretch}{1.5}
\begin{center}
\begin{tabular}{|l|c|c|c|c|c|c|}
\hline
process/vector & $\sBD(t_1,t_2)$ & $\BD^{k(L)}$  & $\mathbb Y^{k(L)}(t)$ & $\mathbb{IP}^{k(L)}(t)$ & $\mathbf{IP}^{k(L)}(t)$ \\ \hline
centering & $L(t_2-t_1) \boldsymbol \mu$ & $(\LL_k-\LL_{k-1}) \boldsymbol \mu$ & $L A^k \boldsymbol \mu$ & $ \LL_k A^k \boldsymbol \mu$ & $ \LL_k A^k \boldsymbol \mu$ \\ \hline
%\begin{tabular}{@{}c@{}}centered-scaled \\ process/vector \end{tabular} 
centered and scaled &  $\hBD(t_1,t_2)$ & $\hat{\BD}^{k(L)}$ & $\hat{\mathbb Y}^{k(L)}(t)$ & $\widehat{\mathbb{IP}}^{k(L)}(t)$   & $\widehat{\mathbf{IP}}^{k(L)}(t)$ \\
\hline
\end{tabular}
\caption{Centered and Scaled Processes and Vectors $(1 \le k \le K)$}
\label{tbl:C-S}
\end{center}
\end{table}

Some other processes and variables are scaled component-wise (without any centering) by 
\[
\hat{\x}^{(L)} (t)= \frac{\x^{(L)}(Lt)}{\sqrt{L}},
\]
and their definitions are summarized in Table \ref{tbl:Sonly}.

\begin{table}[ht]
\renewcommand{\arraystretch}{1.5}
\begin{center}
\begin{tabular}{|l|c|c|c|c|c|c|c|c|}
\hline
process/variable        & $\mathbf d^{(L)}(t)$                  &  $\mathbb Q^{k(L)}(t)$      & $\mathbf Q^{k(L)}(t)$       & $\mathbf B^{(L)}(t)$ & $\sbB(t)$ & $\sbbB(t)$ & $\mathcal C^{(L)}(t)$ & $\mathcal C^{(L)}$ \\ \hline
%scaled process/variable & 
scaled & $\hjD(t)$ & $\hat{\mathbb Q}^{k(L)}(t)$ & $\hat{\mathbf Q}^{k(L)}(t)$ & $\hat{\mathbf B}^{(L)}(t)$ & $\hbB(t)$ & $\hbsB(t)$ & $\hat{\mathcal C}^{(L)}(t)$ & $\hat{\mathcal C}^{(L)}$ \\  \hline
\end{tabular}
\caption{Scaled Processes and Variables $(1 \le k \le K)$}
\label{tbl:Sonly}
\end{center}
\end{table}

\subsection{Asymptotic Optimality Criterion}
\label{subsec:asymcriteria}

Applying Theorem \ref{thm:sameobj} to the $L^{th}$ system, 
%which is an instance of the general model formulated in Section \ref{sec:formulation}, 
the average inventory cost $\mathcal C^{(L)}$ has a lower bound $\underline{C}^{(L)}$, determined by the SP (\ref{eq:MSP2}), using $\BD^{k(L)}$ as inputs $\BD^k$ $(1 \le k \le K)$. 
Define  
\[
\hat{\underline{C}}^{(L)} \stackrel{\Delta}{=}\frac{\underline{C}^{(L)}}{\sqrt{L}}.
\]
%We consider our policy is asymptotically optimal on the diffusion scale

%We also consider another optimality criterion:  
Let $\mincost$ denote the infimum of the average inventory cost over all feasible policies. Since these values are bounded below by $\underline{C}^{(L)}$, this infimum exists. We show that
our policy is asymptotically optimal in the `traditional' sense that
\begin{equation}
\label{eq:asymcost}
\lim_{L \rightarrow \infty} \frac{\mathcal C^{(L)}-\mincost}{\mincost}=0,
\end{equation}
i.e., the percentage difference between the cost and its minimum value converges to zero as the longest lead time becomes large. 

Since $\mincost \ge \underline{C}^{(L)}$, we have 
\[
\frac{\mathcal C^{(L)}-\mincost}{\mincost} \le \frac{\sqrt{L}}{\mincost} \left \{\frac{\mathcal C^{(L)}-\underline{C}^{(L)}}{\sqrt{L}}\right \}.
\]
Hence we can prove (\ref{eq:asymcost}) by showing that 
\begin{equation}
\label{eq:poshold}
\limsup_{L \rightarrow \infty} \left \{\frac{\sqrt{L}}{\mincost}\right \} <\infty.
\end{equation}
and
\begin{equation}
\label{eq:asymobj}
\lim_{L \rightarrow \infty}\frac{\mathcal C^{(L)}-\underline{C}^{(L)}}{\sqrt{L}}=\lim_{L \rightarrow \infty} \left\{\hat{\mathcal C}^{(L)}-\hat{\underline C}^{(L)}\right \} =0.
\end{equation}

To show (\ref{eq:poshold}) is true, apply the SP (\ref{eq:MSP2}) to the $L^{th}$ system with following changes: let $i'$ be a product with $a_{ni'}>0$, where component $n$ has the longest lead time $L$. Keep $h_n$ and $b_{i'}$ intact while resetting $h_j=0$ $(j \ne n$) and $b_i=0$ ($i \ne i'$). The problem then becomes
\[
\min_{y_n} \left\{h_n y_n - (a_{ni'} h_n+b_{i'}) \E \left[\min(y_n/a_{ni'},\bar{D}^{(L)}_{i'})\right] \right\},
\]
where we use $\bar{D}^{(L)}_{i'}$ to denote $D^{K(L)}_{i'}$, which has the same distribution as $D^{(L)}_{i'}(t-L,t)$ $(t \ge 0)$.

Let $\varphi^{K(L)}_n$ be the optimal objective value of the above. Obviously, $\varphi^{K(L)}_n \le \varphi^{K(L)}$, where $\varphi^{K(L)}$ is the objective value of (\ref{eq:MSP2}) specialized to the $L^{th}$ system. Apply $\varphi^{K(L)}_n$ and the new costs to (\ref{eq:sameobj}) and compare the resulting value (denoted by $C^{(L)}_n$) with the lower bound and the minimum cost:
\[
C^{(L)}_n \equiv \varphi^{K(L)}_n+b_{i'} \cdot \E[\bar{D}^{(L)}_{i'}] \le \varphi^{K(L)}+\sum_{i=1}^m b_i \E[\bar{D}^{(L)}_i]= \underline{C}^{(L)} \le \mincost.
\]
Furthermore, simple transformations of $\varphi^{K(L)}_n$ shows that 
\[
\varphi^{K(L)}_n+b_{i'} \cdot \E[\bar{D}^{(L)}_{i'}]=\min_{y_n} \left \{ h_n \E\left[(y_n-a_{ni'} \bar{D}^{(L)}_{i'})^+\right]+\frac{b_{i'}}{a_{ni'}} \E\left[(a_{ni'}\bar{D}^{(L)}_{i'} - y_n)^+\right]\right\}
\]
is a Newsvendor model, so $C^{(L)}_n$ is proportional to the standard deviation of $\bar{D}^{(L)}_{i'}$, which is on the order of $\sqrt{L}$, and (\ref{eq:poshold}) follows as a consequence. 

Hence to show (\ref{eq:asymcost}) is satisfied, we only need to prove (\ref{eq:asymobj}), i.e., our policy is asymptotically optimal on the diffusion scale.

\subsection{Sufficient Conditions for Asymptotic Optimality}
\label{subsec:keytheorem}

The theorem below shows that (\ref{eq:asymobj}) holds if inventory positions and backlog levels converge to their respective targets. Recall that these targets were set at  levels at which the average inventory cost can attain its SP-based lower bound (see (\ref{eq:pmatch}) and (\ref{eq:bmatch})). While meeting these targets exactly is not possible in general, the theorem shows that meeting them on the diffusion scale, which we will prove to be feasible in the next section, is sufficient to guarantee asymptotic optimality.

\begin{theorem}
\label{thm:final}
In a family of systems indexed by their longest lead time $L$, if
\begin{equation}
\label{eq:ipres}
\limsup_{L \rightarrow \infty} \left \{\sup_{t \ge -\hLL_{k_j}} \E\left[\left|\hIP_j(t)-\thIP_j(t)\right|\right]\right \}=0,~\jrange,
\end{equation}
and 
\begin{equation}
\label{eq:blkres}
\limsup_{L \rightarrow \infty} \left\{\sup_{t \ge 0} \E\left[\left|\hB_i(t)-\hbsB_i(t)\right|\right]\right \}= 0,~\irange,
\end{equation}
then (\ref{eq:asymobj}) holds, i.e., our policy is asymptotically optimal on the diffusion scale, and (\ref{eq:asymcost}) follows as a consequence.
\end{theorem}

\smallskip
\noindent
{\bf Proof}: Following (\ref{eq:lbound2}), in the $L^{th}$ system, 
\[
\underline{C}^{(L)}=\sum_{k=1}^K  \mathbf h^k\cdot \left (\E[\y^{k*(L)}]-A^k \E[\BD^{K(L)}] \right)+\mathbf c \cdot \E[\B^{*(L)}],
\]
where $\y^{k*(L)}$ $(1 \le k \le K)$ is obtained from (\ref{eq:MSP2}) and $\B^{*(L)}$ is obtained from (\ref{eq:phi0pi}). 

Specializing (\ref{eq:pmatch-1}) and (\ref{eq:bmatch-1}) in our policy formulation to the $L^{th}$ system, 
 %be the optimal solutions of the SP in (\ref{eq:MSP2}). Then
%Specializing (\ref{eq:pmatch-1}) and (\ref{eq:bmatch-1}) to the $L^{th}$ system, 
%in which case $\mathbb Y^{k(L)}(t-\LL_k)$ $(t \ge 0, 1 \le k \le K)$ and $\bB^{*(L)}(t)$ $(t \ge 0)$ correspond to $\mathbb Y^k(t-L_k)$ $(t \ge 0, 1 \le k \le K)$ and $\bB^(t)$ $(t \ge 0)$ respectively, 
\[
\E[\y^{k*(L)}]=\E[\mathbb Y^{k(L)}(t-\LL_k)],~1 \le k \le K \quad \mbox{and} \quad \E[\B^{*(L)}]=\E[\bB^{*(L)}(t)], \quad t \ge 0.
\]
Moreover, $\BD^{K(L)}$ has the same distribution as $\BD^{(L)}(t-L,t)$ $(t \ge 0)$ by definition. 
Thus 
\begin{align*}
\begin{split}
\underline{C}^{(L)}
&=\sum_{k=1}^K  \mathbf h^k\cdot \left (\E[\mathbb Y^{k(L)}(t-\LL_k)]-A^k \E[\BD^{(L)}(t-L,t)] \right)+\mathbf c \cdot \E[\mathbb B^{*(L)}(t)]\\
&=\sum_{k=1}^K  \mathbf h^k\cdot \left (\E[\mathbb{IP}^{k(L)}(t-\LL_k)]- A^k \E[\BD^{(L)}(t-\LL_k,t)] \right)+\mathbf c \cdot \E[\mathbb B^{*(L)}(t)], \quad t \ge 0,
\end{split} 
\end{align*}
where the second equality results from the determination of inventory position targets in (\ref{eq:target}). 
%Applying the above to (\ref{eq:lbound2}) and using (\ref{eq:target}) to substitute $\mathbb Y^{k(L)}(t-\LL_k)$ in the $L^{th}$ system, 

By the definition of $\hat{\underline{C}}^{(L)}$ and applying the second equality of the above,
\begin{align*}
\begin{split}
\hat{\underline{C}}^{(L)} 
&=\frac{\sum_{k=1}^K  \mathbf h^k\cdot \left (\E[\mathbb{IP}^{k(L)}(Lt-\LL_k)]- A^k \E[\BD^{(L)}(Lt-\LL_k,Lt)] \right)+\mathbf c \cdot \E[\mathbb B^{*(L)}(Lt)]}{\sqrt{L}}\\
 &=\sum_{k=1}^K \mathbf h^k \cdot \left (\E\left[\widehat{\mathbb{IP}}^{k(L)}(t-\hLL_k)\right]- A^k \E\left[\hBD(t-\hLL_k,t)\right]\right)+\mathbf c \cdot \E\left[\hat{\mathbb B}^{*(L)}(t)\right], \quad t \ge 0. 
\end{split}
\end{align*}
Comparing the above with (\ref{eq:scalecostexp}),
\begin{eqnarray}
\nonumber
\hat{\mathcal C}^{(L)}-\hat{\underline{C}}^{(L)} 
&=& \limsup_{T \rightarrow \infty} \frac{1}{T} \int_0^{T} \left(  \hat{\mathcal C}^{(L)} (t) - \underline{\hat{C}}^{(L)}\right) dt
\\
\nonumber
&=& \limsup_{T \rightarrow \infty} \frac{1}{T} \int_0^{T} \sum_{k=1}^K \mathbf h^k \cdot \left(\E\left[\widehat{\mathbf{IP}}^{k(L)}(t-\hLL_k)\right]-\E\left[\widehat{\mathbb{IP}}^{k(L)}(t-\hLL_k)\right]\right) dt
\\
&~& + \limsup_{T \rightarrow \infty} \frac{1}{T} \int_0^{T}\mathbf c \cdot \left(\E\left[\hat{\B}^{(L)}(t)\right]-\E\left[\hat{\mathbb B}^{*(L)}(t)\right]\right)dt,
\label{eq:errorsum}
%\\&\le& \sum_{j=1} h_j\sup_{t \ge 0} \left| \E\left[\hIP_j(t)-\thIP_j(t)\right]\right|+ \sum_{i=1}^m b_i  \sup_{t \ge 2} \left |\E\left [\hB_i(t)-\hbsB_i(t)\right] \right|.
\end{eqnarray}
The theorem follows immediately by applying (\ref{eq:ipres}) and (\ref{eq:blkres}) to (\ref{eq:errorsum}). $\qquad \blacksquare$

\section{Proof of Sufficient Conditions}
\label{sec:sufficond}

%We have shown in Section \ref{sec:policy} that the average cost of the ATO system could have reached the SP-based lower bound if both inventory positions and backlog levels could be kept at their target levels for all time. While the latter conditions are not feasible in general, Theorem \ref{thm:final} shows that we can still achieve asymptotic optimality on the diffusion scale if both levels converge to their targets. 
Continuing from Theorem \ref{thm:final}, we prove in this section that (\ref{eq:ipres})-(\ref{eq:blkres}) are satisfied under Target Stability Conditions to be defined below. To put results obtained in this section in perspective, it was shown in \cite{Reiman2015} that for systems with identical lead times only (\ref{eq:blkres}) was needed since asymptotic optimality can be achieved under a base-stock policy that keeps the inventory positions constant. The convergence (\ref{eq:blkres}) was shown in Theorem 4 in \cite{Reiman2015}, while here it is shown in Corollary \ref{cor:backlog} below. The proof of (\ref{eq:blkres}) in \cite{Reiman2015} does not carry over to this setting because it assumes the use of a base stock policy for replenishment.

%When developing our inventory policy in Section \ref{sec:policy}, we point out that the resulting average inventory cost generally does not reach the SP-based lower bound because it is not possible to keep inventory positions and backlog levels at their target levels for all time. We then show in Theorem \ref{thm:final} that our policy is asymptotically optimal on the diffusion scale if, as are specified by (\ref{eq:ipres}) and (\ref{eq:blkres}), inventory positions and backlog levels converge to their targets on the diffusion scale. In this section, we show that the latter sufficient conditions can indeed be satisfied, provided that the following stability conditions hold as we keep re-solving the aforementioned SP to update these targets:
   
\medskip
\noindent {\bf Target Stability Conditions:}
Under our inventory policy, $\mathbb Y^k(t)$, $(t \ge -L_k, 1 \le k \le K)$, obtained by solving (\ref{eq:Yset}), $\mathbb B^*(t)$ ($t \ge 0$), obtained by solving (\ref{eq:miniB})-(\ref{eq:iQ}), and $\mathbb B(t)$ $(t \ge 0)$, obtained by solving (\ref{eq:setBtar}), have the following properties:  there exists some constant $\kappa$, depending only on $\mathbf h$, $\mathbf b$, and $A$, such that on each sample path of these processes, 
\begin{itemize}
\item
for all $t_2 >t_1  \ge -L_{k_j}$, 
\begin{equation}
\label{eq:iptargetcont}
|\mathbb Y_j(t_2)-\mathbb Y_j(t_1)| \le \kappa \sum_{k>k_j} ||\mathbf D(t_2-(L_k-L_{k_j}),t_2)-\mathbf D(t_1-(L_k-L_{k_j}),t_1)||_1, \quad 1 \le j \le n;
\end{equation}
\item
for all $t \ge 0$, 
\begin{equation}
\label{eq:btarclose}
|\mathbb B_i(t)-\mathbb B^*_i(t)| \le \kappa \sum_{j=1}^n |Q_j(t)-\mathbb Q_j(t)|, \quad 1 \le i \le m,
\end{equation}
and for all $t_2>t_1 \ge 0$,
\begin{equation}
\label{eq:blktargetcont}
|\mathbb B^*_i(t_2)-\mathbb B^*_i(t_1)| \le \kappa \sum_{j=1}^n |\mathbb Q_j(t_2)-\mathbb Q_j(t_1)|,\quad 1 \le i \le m. 
\end{equation}
\end{itemize}

\smallskip
We will develop an SP solution procedure in Section \ref{sec:SPsolution} to satisfy (\ref{eq:iptargetcont})-(\ref{eq:blktargetcont}). Here in this section, we assume these conditions hold to prove (\ref{eq:ipres}) and (\ref{eq:blkres}).  

Let us first give an informal explanation of the intuition of our proof.  Under our replenishment policy, an inventory position can differ from its target only by exceeding it. When this happens, our policy will stop ordering the component, so its inventory position drops at the same rate as its demand arrival rate. The target may also change as new demand arrives. However, condition (\ref{eq:iptargetcont}) ensures that the latter change in the target is ``slow'' in comparison with demand arrivals, so the excess of inventory position over its target can be eliminated fast enough to satisfy (\ref{eq:ipres}). 
We then use a similar argument to show that under conditions (\ref{eq:btarclose})-(\ref{eq:blktargetcont}), (\ref{eq:blkres}) cannot be violated by having a backlog level falling below its target. To show that (\ref{eq:blkres}) can also not be violated by having a backlog level exceeding its target, we make use of the property of our allocation policy given in (\ref{eq:twoside}), i.e., no backlog level will exceed its target when no other backlog level is below its target. 

To reduce redundancy and improve generality,  in Section \ref{subsec:tracking}, we define a general problem, referred to as the ``Stochastic Tracking Model'', based on common features of (\ref{eq:ipres}) and (\ref{eq:blkres}). We develop Theorem \ref{thm:master}, which applies to the general model and contains (\ref{eq:ipres}) and (\ref{eq:blkres}) as special cases.

\subsection{Stochastic Tracking Model}
\label{subsec:tracking}

Consider a family of systems indexed by $L >0$, where the $L^{th}$ system is associated with the following parameters and processes:
\begin{equation}
\label{eq:stm}
\big\{\sD(t), t \ge -L; ~\mA;~s_l^{(L)}~ (l \in \mK);~\tL_0;~\sT(t), t \ge t_0;~\sW(t),t \ge t_0;~\sW_0\big\}.
\end{equation}
Here,
\begin{itemize}
\item
$\sD(t)$ ($t \ge -L)$ is the same compound Poisson demand process defined in Section \ref{subsubsec:demand}  for the $L^{th}$ system, with the arrival rate $\underline{\lambda}$ and order sizes $\mathbf S$ apply to all systems. For the results of this section, it is important to recall that components of $\mathbf S$, $S_i$ $(\irange)$, are assumed to have a finite moment of order $6$. As in Section \ref{subsubsec:demand}, the increment of $\sD(t)$ over interval $(t_1,t_2]$ is denoted by $\sBD(t_1,t_2)$ $(-L \le t_1<t_2)$ and the jump of $\sD(t)$ at time $t$ is denoted by $\sjD(t)$ $(t \ge -L)$. 
\item
$\mA$ is an $m$-dimensional vector of non-negative integers. The value of the vector is the same for all systems and at least one of its elements is strictly positive.
\item
Each system is associated with a set of constants $s_l^{(L)}$ ($l \in \mK$). Values of $s_l^{(L)}$ $(l \in \mK)$ differ between systems, but the index set $\mK$, which is finite, remains the same for all systems. In the $L^{th}$ system,   
\[
0 \le s_l^{(L)} \le L, \quad (l \in \mK).
\]
\item
$\sT(t)$ $(t \ge \tL_0)$, which we refer to as the \emph{target process}, is a pure jump process that starts at time $\tL_0$, where in the $L^{th}$ system,
\[
-L+s^{(L)}_l \le \tL_0 \le 0, \quad \mbox{and thus} \quad \tL_0 -s^{(L)} \ge -L, \quad \mbox{for all}~l \in \mK. 
\]
\item
$\sW(t)$ $(t \ge \tL_0)$, which we refer to as the \emph{tracking process}, is another pure jump process that also starts at time $\tL_0$ from an initial level $\sW_0$. The process is defined by:
\begin{align}
\label{eq:psctp}
\begin{split}
\sWp(\tL_0)&=\sW_0,
\\
\sW(t)&=\sWp(t)\vee \sT(t)=\sT(t)+\left(\sWp(t)-\sT(t)\right)^+, \quad t \ge \tL_0,
\\
\sWp(t)&=\sW(t^-)-\mA \cdot \sjD(t), \quad t>\tL_0.
\end{split}
\end{align}
\end{itemize}
For convenience, we denote the initial difference between the target and tracking processes by 
\[
\mG^{(L)}_0=\sW(\tL_0)-\sT(\tL_0)=\left(\sW_0-\sT(\tL_0)\right)^+
\]
Observe that the tracking process can instantly catch the target that rises higher than its level. When the target is lower, the gap can be closed after a sufficient amount of demand has arrived.    

Let $\hBD(t_1,t_2)$ $(-1 \le t_1<t_2)$ and $\hjD(t)$ $(t \ge -1)$ be defined the same as in (\ref{eq:scaledemand}). 
\iffalse 
\begin{align*}
\begin{split}
&\hBD(t_1,t_2) \stackrel{\Delta}{=}\frac{\sBD(Lt_1,Lt_2)-L(t_2-t_1){\boldsymbol \mu}}{\sqrt{L}}, \quad -1 \le t_1 < t_2,
\\
\mbox{and} \quad&
\hjD(t)\stackrel{\Delta}{=}\frac{\sjD(Lt)}{\sqrt{L}},\quad t \ge -1.
\end{split}
\end{align*}
\fi
Let 
\[
\sL_l\stackrel{\Delta}{=}\frac{s_l^{(L)}}{L}~(l \in \mK) \quad \mbox{and} \quad \sttL\stackrel{\Delta}{=}\frac{t_0^{(L)}}{L}.
\]
Following the conditions imposed on $s^{(L)}_l$ and $\tL_0$ in their definitions,
\[
0 \le \sL_l < 1 \quad \mbox{and}\quad -1+\sL_l \le \sttL < 0,\qquad l \in \mK. 
\]
Let
\begin{align}
\begin{split}
\label{eq:scaling}
 \hT(t)\stackrel{\Delta}{=}\frac{\sT(Lt)}{\sqrt{L}} \quad &\mbox{and} \quad \hW(t) \stackrel{\Delta}{=} \frac{\sW(Lt)}{\sqrt{L}}, \qquad t \ge \sttL,
\\
\mbox{with} \quad \hW_0\stackrel{\Delta}{=}\frac{\sW_0}{\sqrt{L}} \quad &\mbox{and} \quad \hG_0 \stackrel{\Delta}{=} \frac{\mG^{(L)}_0}{\sqrt{L}}=\left(\hW_0-\hT(\sttL)\right)^+.
\end{split}
\end{align}

Finally, the definition of the model also requires the following asymptotic Lipschitz Continuity condition on the target processes in this family of systems:

\medskip
\noindent
{\bf Asymptotic Lipschitz Condition:} there exists some constant $\kappa$ that applies to all systems, such that 
for all $\tL_0 \le t_1 < t_2$,
\begin{equation}
\label{eq:Lcontinuity}
|\mT^{(L)}(t_2)-\mT^{(L)}(t_1)| \le \kappa \sum_{l \in \mK} ||\sBD(t_2-s^{(L)}_l,t_2)-\sBD(t_1-s^{(L)}_l,t_1)||_1+\sE(t_1)+\sE(t_2)
%\quad \tL_0 \le t_1 < t_2,
\end{equation}
where $\sE(t)$ ($t \ge \tL_0$) satisfies the following condition:
let 
\[
\hE(t) \stackrel{\Delta}{=} \frac{\sE(Lt)}{\sqrt{L}}, ~~t \ge \sttL.
\]
Then 
\begin{equation}
\label{eq:alcon}
\lim_{L \rightarrow \infty} \E\left[\sup_{t \ge \sttL} |\hE(t)|\right]= 0. 
\end{equation}
%and in each system, there exists a process $\mT_c^{(L)}(t)$ $(t \ge \tL_0)$ such that 
\iffalse
\begin{align}
\begin{split}
\label{eq:Lcontinuity}
& |\mT^{(L)}_c(t_2)-\mT^{(L)}_c(t_1)| \le \kappa \sum_{l \in \mK} ||\sBD(t_2-s^{(L)}_l,t_2)-\sBD(t_1-s^{(L)}_l,t_1)||_1,\quad \tL_0 \le t_1 < t_2,
\\
\mbox{and} & \qquad
\lim_{L \rightarrow \infty} \E\left [ \sup_{t \ge \sttL} \frac{|\sT(Lt)-\sT_c(Lt)|}{\sqrt{L}}\right]=0.
\end{split}
\end{align}

\medskip
For convenience, we define  
\[
\sE(t)\stackrel{\Delta}{=} \sT(t)-\sT_c(t),~~ t \ge -L \qquad\mbox{and} \qquad \hE(t) \stackrel{\Delta}{=} \frac{\sE(Lt)}{\sqrt{L}}, ~~t \ge \sttL,
\]
and write the latter condition as
\begin{equation}
\label{eq:alcon}
\lim_{L \rightarrow \infty} \E\left[\sup_{t \ge \sttL} |\hE(t)|\right]= 0. 
\end{equation}
\fi

The main conclusion we draw from this model is in Theorem \ref{thm:master} below, which shows that as $L$ increases, the tracking process converges to the target process on the diffusion scale.

\begin{theorem}
\label{thm:master}
Assume that $S_i$ $(\irange)$ has a finite moment of order $6$. If
\begin{equation}
\label{eq:errorini}
%\lim_{L \rightarrow \infty} \E\left[\sup_{t \ge \sttL} |\hE(t)|\right]= 0~~~\mbox{and}~~~ 
\lim_{L \rightarrow \infty} \E\left[\hG_0\right]=0,
\end{equation}
 then 
\begin{equation}
\label{eq:limitcon0}
\lim_{L \rightarrow \infty} \E\left[\sup_{t \ge \sttL} \left\{\hW(t)-\hT(t)\right\}\right]=0.
\end{equation}
\end{theorem}

Before proving the theorem, we provide some context and intuition for the result.
The convergence in (\ref{eq:limitcon0}) is, in essence, an example of ``state-space collapse'', a phenomenon that plays a critical role in the heavy traffic analysis of queueing systems and stochastic processing networks (\cite{Reiman1984},\cite{Harrison1997},\cite{Bramson1998},\cite{Bell2001},\cite{Atar2019}). 
As noted above, the tracking process $\hW$ is able to instantly match the upward jumps of the target process, but not the downward jumps. 
Under the scaling that we have introduced, the family of processes $\hBD$ satisfies a functional central limit theorem (FCLT). This implies that $\hBD$ has ``almost continuous'' paths for large $L$. The assumed asymptotic Lipschitz continuity of $\hT$ implies that it has almost continuous paths as well.
When $\hW(t) > \hT(t)$, it is following total demand downward. Note that this is uncentered demand, so that $\hW$ is decreasing at a ``rate'' that is $O(\sqrt{L})$ when $\hW(t) > \hT(t)$. Thus $\hW $ never gets ``far away'' from $\hT$, and the two dimensional-processes $(\hT(t),\hW(t))$ ``collapse'' into a one-dimensional process with $\hT(t)=\hW(t)$ in the limit. The proof of Theorem \ref{thm:master} makes this intuitive description rigorous.

As noted above, $\hBD$ satisfies an FCLT. Interestingly, we never actually invoke this fact in our proofs.
A key reason is that it would not be enough to obtain the result that we need. FCLTs involve convergence over a finite time interval.
In order to prove asymptotic optimality under the long run average cost criterion that we use, we need uniform convergence over an infinite time interval, which does not follow from the FCLT.

The proof of Theorem \ref{thm:master} is built upon the following four lemmas. The first two are restatements of two lemmas in Reiman and Wang (2015) and the next two are new. 
In all lemmas, order sizes $S_i$ are assumed to have a finite moment of order $2+\delta$ ($\delta>0$), except for Lemma \ref{lem:NBD}, which requires a stronger condition of having a finite moment of order $6$.
%A key result that we use repeatedly in these proofs is Doob's inequality.

\begin{lemma} {(Lemma 2 in Reiman and Wang (2015))}
\label{lem:RW2}
For all $i=1,\cdots,m$, 
\begin{equation}
\E\left[\sup_{0 \le \tau \le 1} \hejD_i(\tau)\right] \le 3 \underline{\lambda}^{1/(2+\delta)} (1+\eta_i)L^{-\delta/(2(2+\delta))},
\end{equation}
where $\underline{\lambda}$ is the demand arrival rate and $\eta_i :=\E[S_i^{2+\delta}]$ $(\delta>0)$. 
\end{lemma}

\begin{lemma} {(Lemma 3 in Reiman and Wang (2015))}
\label{lem:RW3} 
For all $i=1,\cdots,m$,
\begin{equation}
\label{eq:boundonED}
\E \left [\sup_{0 \le \tau \le L^{-1/4}}|\hD_i(0,\tau)|\right]\le (1+\sigma^2_{ii})L^{-1/4}.
\end{equation}
where $\sigma_{ii}$ is the variance of $D_i(0,1)$. 

Let $g$ be any strictly positive constant. Then for all $i=1,\cdots,m$,
\begin{equation}
\label{eq:boundonEDD}
\E\left [\sup_{L^{-1/4} \le \tau \le 1}(|\hD_i(0,\tau)|-\sqrt{L}\tau g )^+\right]\le \frac{\sigma^2_{ii}}{g} L^{-1/4}.
\end{equation}
\end{lemma}

\begin{lemma} {(see Appendix A for the proof)}
\label{lem:NBD}
Assume that $S_i$ has a finite moment of order $6$ $(1 \le i \le m)$. Let $\nu$ be any strictly positive constant. Then for all $i=1,\cdots,m$,
\begin{equation}
\label{eq:dlim}
\lim_{L \rightarrow \infty}
\E\left[\sup_{t \ge 0}\left(|\hD_i(0,t)|-\sqrt{L}t\nu\right)^+\right]=0.
\end{equation} 
\end{lemma}

\begin{lemma}{(see Appendix A for the proof)}
\label{lem:RW4}
Let $\nu$ be any strictly positive constant. Then for all $i=1,\cdots,m$, 
\begin{equation}
\label{eq:increkey}
\lim_{L \rightarrow \infty}
\E\left[\sum_{\tau=0}^\infty \sup_{0\le t <1}\left(\hejD_i(t)-\sqrt{L}\nu\tau\right)^+\right]=0.
\end{equation}
\end{lemma}

The theorem is proved below by showing that there is uniform upper bound on the expected value at the left-hand side on (\ref{eq:limitcon0}) and the bound converges to $0$ as $L$ increases.  

\smallskip
\noindent
{\bf Proof of Theorem \ref{thm:master}}: 
We prove (\ref{eq:limitcon0}) in the theorem by bounding $(\hW(t)-\hT(t))^+$ and prove the bound converges uniformly to $0$ for all $t \ge \sttL$. 
%For any give time $t$ $(t \ge \sttL)$,
By definition, 
\[
\hW(t) \ge \hT(t),\quad t \ge \sttL,
\]
so to bound $(\hW(t)-\hT(t))^+$, we only need to consider the case where $\hW(t)>\hT(t)$, 
and for the latter case, only need to consider changes of $\hW(t)-\hT(t)$ since the last moment before $t$ when $\hW(\cdot)$ jumps from $\hW(\tau^-)=\hT(\tau^-)$ to $\hW(\tau)>\hT(\tau)$. 
That moment is defined by 
\[
\stL=\sup_{\sttL \le \tau \le t} \{\tau:\hW(\tau)=\hT(\tau)\},
\]
or if the set on the right-hand side is empty, let $\stL=\sttL$. Observe that while not explicit in its notation, $\stL$ depends on $t$.

Observe that
\begin{equation}
\label{eq:setdiff}
|\hW(t)-\hT(t)|=(\hW(t)-\hT(t))\mathbf 1\{\stL=\sttL\}+(\hW(t)-\hT(t))\mathbf 1\{\stL>\sttL\},
\end{equation}
so the proof of the theorem is reduced to proving
\begin{eqnarray}
\label{eq:casetl0}
\lim_{L \rightarrow \infty} \E\left[\sup_{t \ge \sttL} \left \{(\hW(t)-\hT(t)) \mathbf 1\{\stL=\sttL\} \right\}\right] &=& 0,
\\
\mbox{and}~~
\label{eq:casetl1}
\lim_{L \rightarrow \infty}\E\left[\sup_{t \ge \sttL} \left \{(\hW(t)-\hT(t))\mathbf 1\{\stL>\sttL\}\right\}\right] &=&0. 
\end{eqnarray}

To prove (\ref{eq:casetl0}), when $\stL=\sttL$, $\hW(\tau)>\hT(\tau)$ for all $\tau \in [\sttL,t]$, so by (\ref{eq:psctp}),
\[
\hW(t)-\hW(\sttL)=-\mA \cdot \frac{\sBD(L\sttL,Lt)}{\sqrt{L}}=-\mA \cdot \left(\hBD(\sttL,t)+\sqrt{L} {\boldsymbol \mu}(t-\sttL)\right). 
\]
By applying the above and condition (\ref{eq:Lcontinuity}), for any time $t$ $(t \ge \sttL)$, 
\begin{align*}
\begin{split}
\nonumber
%replaced two hT_c below with hT_c and commented out \hE
\hW(t)-\hT(t) &= (\hW(t)-\hW(\sttL))+\hG_0-(\hT(t)-\hT(\sttL))
%-\hE(t)+\hE(\sttL)
\\
&\le -\mA \cdot (\hBD(\sttL,t)+\sqrt{L} {\boldsymbol \mu}(t-\sttL))+\hG_0
\\
&~~~~~~
+ \kappa \sum_{l \in \mK} ||\hBD(t-\sL_l,t)-\hBD(\sttL-\sL_l,\sttL)||_1
+|\hE(t)|+|\hE(\sttL)|.
\end{split}
\end{align*}
Therefore  under condition (\ref{eq:errorini}), (\ref{eq:casetl0}) holds if 
\begin{align}
\begin{split}
\label{eq:case0cond}
&~\lim_{L \rightarrow \infty} \E \Big [\sup_{t \ge \sttL} \big \{-\mA \cdot [\hBD(\sttL,t)+\sqrt{L} {\boldsymbol \mu}(t-\sttL)]
\\
&~~~~~~~~~~~~~~~~~~+\kappa \sum_{l \in \mK} ||\hBD(t-\sL_l,t)-\hBD(\sttL-\sL_l,\sttL)||_1 \big \} \Big] = 0, 
\end{split}
\end{align}
which we prove next. 

Let
\[
\zeta=\frac{\mA \cdot {\boldsymbol \mu}}{m(||\mA||_\infty+2\kappa|\mK|)} >0~~\mbox{and}~~~\tilde{\kappa}=||\mA||_\infty+\kappa |\mK| > 0. 
\]
Then
\begin{align}
\begin{split}
\label{eq:boundinside}
&~~-\mA \cdot [\hBD(\sttL,t)+\sqrt{L} {\boldsymbol \mu}(t-\sttL)]+\kappa \sum_{l \in \mK} ||\hBD(t-\sL_l,t)-\hBD(\sttL-\sL_l,\sttL)||_1
\\
\le & ||\mA||_\infty ||\hBD(\sttL,t)||_1-\sqrt{L} \mA \cdot {\boldsymbol \mu}(t-\sttL)
\\
&~~~~~~~~~~+ \kappa \sum_{l \in \mK} (||\hBD(\sttL,t)||_1+||\hBD(\sttL-\sL_l,t-\sL_l)||_1)
\\
\le & \tilde{\kappa} \sum_{i=1}^m(|\hD_i(\sttL,t)|-\sqrt{L} \zeta (t-\sttL))^+
\\
&~~~~~~~~~~+ \kappa  \sum_{l \in \mK} \sum_{i=1}^m(|\hD_i(\sttL-\sL_l,t-\sL_l)|-\sqrt{L} \zeta (t-\sttL))^+.
\end{split}
\end{align}
Since $\sD(t)$ ($t \ge -L$) is a stationary process,
\begin{equation}
\label{eq:equidist}
|\hD_i(\sttL,t)|\stackrel{d}{=}|\hD_i(\sttL-\sL_l,t-\sL_l)|\stackrel{d}{=}|\hD_i(0,t-\sttL)|,~l \in \mK,~\irange. 
\end{equation}
Since $\sttL \in [-1,0]$, for all $i=1,\cdots,m$, 
\begin{align}
\begin{split}
\label{eq:conv0}
~&~~\lim_{L \rightarrow \infty} \E\left[\sup_{t \ge \sttL } \left (|\hD_i(0,t-\sttL)|-\sqrt{L} \zeta (t-\sttL)\right)^+\right] 
\\
&\le \lim_{L \rightarrow \infty} \E\left[\sup_{t \ge t'} \sup_{-1 \le t' \le 0}  \left (|\hD_i(0,t-t')|-\sqrt{L} \zeta (t-t')\right)^+ \right]
\\
&=\lim_{L \rightarrow \infty}\E\left[\sup_{t \ge 0} \left (|\hD_i(0,t)|-\sqrt{L} \zeta t\right)^+\right]
 \\
&= 0~~~~~~\mbox{(Lemma \ref{lem:NBD})},
\end{split}
\end{align}
and (\ref{eq:case0cond}) follows directly from (\ref{eq:boundinside})-(\ref{eq:conv0}).

To prove (\ref{eq:casetl1}), by the definition of $\stL$, in cases where $\stL>\sttL$, 
\[
\hW(\stLp)=\hT(\stLp),
\]
and thus
\begin{eqnarray}
\nonumber
\hW(t)-\hT(t) &=& \hW(t)-\hT(t)-[\hW(\stLp)-\hT(\stLp)]
\\
\nonumber
&=& [\hW(t)-\hW(\stL)]+[\hW(\stL)-\hW(\stLp)]
\\
\label{eq:items}
&~& 
+[\hT(\stLp)-\hT(\stL)]+[\hT(\stL)-\hT(t)].
\end{eqnarray}
Since $\hW(\tau)>\hT(\tau)$ for all $\tau \in (\stL,t]$,
by (\ref{eq:psctp}),
\begin{align}
\begin{split}
\label{eq:bnditemca1}
\hW(t)-\hW(\stL) &=-\mA \cdot (\hBD(\stL,t)+\sqrt{L} {\boldsymbol \mu}(t-\stL)).
\\
& \le ||\mA||_\infty \sum_{i=1}^m |\hD_i(\stL,t)|-\sqrt{L}(t-\stL) \mA \cdot {\boldsymbol \mu},
\end{split}
\end{align}
and 
\begin{equation}
\label{eq:bnditemca2}
\hW(\stL)-\hW(\stLp)=-\mA \cdot \hjD(\stL) \le 0.
\end{equation}
For the last two items in (\ref{eq:items}), by (\ref{eq:Lcontinuity}),
\begin{align}
\begin{split}
\label{eq:bnditemca3}
&~~~~~\left|\hT(\stLp)-\hT(\stL) \right|
% commented this line \\&= \hT_c(\stLp)-\hT_c(\stL)+\hE(\stLp)-\hE(\stL)
\\
& \le
\kappa \sum_{l \in \mK}\left(|| \hjD(\stL-\sL_l)||_1+||\hjD(\stL)||_1\right)+\hE(\stLp)+\hE(\stL)
\\
& =  \kappa  \sum_{i=1}^m \sum_{l \in \mK}\left(|\hejD_i(\stL-\sL_l)|+|\hejD_i(\stL)|\right)+\hE(\stLp)+\hE(\stL),
\end{split}
\end{align}
and 
\begin{align}
\begin{split}
\label{eq:bnditemca4}
&~~~~~~\left|\hT(\stL)-\hT(t) \right|
% commented this line \\&= \hT_c(\stL)-\hT_c(t) +\hE(\stL)-\hE(t)
\\
& \le
\kappa \sum_{l \in \mK}||\hBD(t-\sL_l,t)-\hBD(\stL-\sL_l,\stL)||_1 +\hE(\stL)+\hE(t)
\\
& =\kappa \sum_{l \in \mK}||\hBD(\stL,t)-\hBD(\stL-\sL_l,t-\sL_l)||_1 +\hE(\stL)+\hE(t)
\\
& \le 
\kappa \sum_{i=1}^m \sum_{l \in \mK}\left(|\hD_i(\stL,t)|+|\hD_i(\stL-\sL_l,t-\sL_l)|\right) +\hE(\stL)+\hE(t).
\end{split}
\end{align}
%where (\ref{eq:Lcontinuity}) is used to bound increments of $\mT_c$ in both cases. 

Let
\[
\nu=\frac{\mA \cdot \boldsymbol \mu}{m(||\mA||_\infty+4 \kappa |\mathcal K|)}.
\]
By applying (\ref{eq:bnditemca1}), (\ref{eq:bnditemca2}), (\ref{eq:bnditemca3}), and (\ref{eq:bnditemca4}) to bound each item in the last expression of (\ref{eq:items}) and replacing $\mA \cdot \boldsymbol \mu$ with 
$m(||\mA||_\infty+4 \kappa |\mathcal K|) \nu$, we can prove that (\ref{eq:casetl1}) holds under condition (\ref{eq:errorini}) by showing that for all $i=1,\cdots,m$ and $l \in \mK$,
\begin{eqnarray}
\label{eq:limit1}
\lim_{L \rightarrow \infty} \E\left[\sup_{\sttL \le t' \le t < \infty} \left(|\hD_i(t',t)|-\sqrt{L}(t-t') \nu \right)^+\right] &=& 0,
\\
\label{eq:limit2}
\lim_{L \rightarrow \infty} \E\left[\sup_{\sttL \le t' \le t<\infty} \left(|\hD_i(t'-\sL_l,t-\sL_l)|-\sqrt{L}(t-t') \nu \right)^+\right] &=& 0, 
\\
\label{eq:limit3}
\lim_{L \rightarrow \infty} \E\left[\sup_{\sttL \le t' \le t < \infty}\left(\hejD_i(t'-\sL_l)-\sqrt{L}(t-t')\nu \right)^+\right] &=& 0,
\\
\label{eq:limit4}
\mbox{and}~~~~~~~~~~~~~~~~~~~~\lim_{L \rightarrow \infty} \E\left[\sup_{\sttL<t' \le t< \infty}\left(\hejD_i(t')-\sqrt{L}(t-t')\nu \right)^+\right]  &=& 0.
\end{eqnarray}

For each given $i$ $(\irange)$ and $l$ $(l \in \mK)$,  
\begin{eqnarray*}
\hD_i(t',t) &\stackrel{d}{=}& \hD_i(t'-\sL_l,t-\sL_l)\stackrel{d}{=}\hD_i(0,t-t'),\quad \sttL \le t' \le t,
\\
\hejD_i(t-\sL_l) &\stackrel{d}{=}&\hejD_i(t'), \quad t' \ge \sttL.
\end{eqnarray*}
(Recall that $-1+\sL_l \le \sttL$, so $t'-\sL_l \ge -1$ in the above). 
Therefore, (\ref{eq:limit1}) and (\ref{eq:limit2}) follow from a similar argument that proves (\ref{eq:conv0}): by the use of Lemma \ref{lem:NBD},
\[
\lim_{L \rightarrow \infty} \E\left[\sup_{\sttL \le t' \le t <\infty} \left(|\hD_i(t',t)|-\sqrt{L}(t-t') \nu \right)^+\right]=\lim_{L \rightarrow \infty} \E\left[\sup_{t \ge 0} \left(|\hD_i(0,t)|-\sqrt{L}t \nu \right)^+\right]=0,
\]
and
\begin{align*}
\begin{split}
\lim_{L \rightarrow \infty} \E\left[ \sup_{\sttL \le t' \le t<\infty} \left(|\hD_i(t'-\sL_l,t-\sL_l)|-\sqrt{L}(t-t') \nu \right)^+\right]
%\\=& \lim_{L \rightarrow \infty} \E\left[\sup_{\sttL \le t' \le t} \left(|\hD_i(t',t)|-\sqrt{L}(t-t') \nu \right)^+\right]\\
=& \lim_{L \rightarrow \infty} \E\left[\sup_{t \ge 0} \left(|\hD_i(0,t)|-\sqrt{L}t \nu \right)^+\right]
\\
=& 0. 
\end{split}
\end{align*}

To prove (\ref{eq:limit3}), since $ -1+\sL_l \le \sttL \le t' \le t$ and the demand process is stationary,
\begin{align*}
\begin{split}
&~\E\left[\sup_{\sttL \le t' \le t<\infty}\left(\hejD_i(t'-\sL_l)-\sqrt{L}(t-t')\nu \right)^+\right]
\\
\le &
\E\left[\sup_{-1 \le t' \le t<\infty}\left(\hejD_i(t')-\sqrt{L}(t-t')\nu \right)^+\right]
\\
\le & \E\left[\sup_{t \ge 0} \left (\sum_{\tau=0}^{\lfloor t \rfloor}\left(\hejD_i(t')-\sqrt{L}\tau\nu \right)^+ \mathbf 1(t-\tau-1 \le t' \le t-\tau)+\sup_{-1 \le t' \le 0} \hejD_i(t')\right) \right]
\\
\le & \E\left[\sum_{\tau=0}^\infty \sup_{0 \le t' \le 1} \left(\hejD_i(t')-\sqrt{L}\tau\nu \right)^+\right]+\E\left[\sup_{-1 \le t' \le 0} \hejD_i(t') \right],
\end{split}
\end{align*}
and (\ref{eq:limit3}) follows from Lemmas \ref{lem:RW2} and \ref{lem:RW4}. 
We then use the same argument to prove (\ref{eq:limit4}), and thus complete the proof of (\ref{eq:casetl1}).~$\blacksquare$

\subsection{Convergence to Targets}
\label{subsec:proofcond}

We again consider a family of ATO systems indexed by $L$, with $\sD(t)$ $(t \ge -L)$ as the demand process in the $L^{th}$ system.  For each component and for each product, we specify other parameters in (\ref{eq:stm}) to define an instance of the stochastic tracking model, and prove (\ref{eq:ipres}) and (\ref{eq:blkres}) as corollaries of Theorem \ref{thm:master}. In fact, the theorem shows convergence of $\E[\sup(\cdot)]$ over an infinite time horizon, which is a stronger result than what is needed, which is the convergence of $\sup \E[(\cdot)]$. 

The proofs use the following lemma.

\begin{lemma}{(see Appendix A for the proof)}
	\label{lem:initconv}
	Let $\hbY_j(-\hLL_{k_j})$ be defined in (\ref{eq:scaleinv}) with $t=-\hLL_{k_j}$ $(\jrange)$ and $\nu$ be any strictly positive constant. Then 
	\begin{equation}
	\label{eq:initdon}
	\lim_{L \rightarrow \infty} \E[(|\hbY_j(-\hLL_{k_j})|-\sqrt{L} \nu)^+]=0,~~~\jrange. 
	\end{equation}
\end{lemma}

\iffalse
\begin{table}[ht]
\begin{center}
\begin{tabular}{|c|c|c|c|} 
\hline
& general model& Condition (\ref{eq:ipres}) & Condition (\ref{eq:blkres}) \\ \hline
target process & $\sT(t)$ & $\lceil\tIP^{(L)}_j(t)\rceil$ & $-\sbB_i(t)$ \\ \hline 
tracking process & $\sW(t)$ & $IP_j(t)$ & $-\nsB_i(t)$ \\ \hline
& $\sE(t)$ & $\lceil\tIP^{(L)}_j(t)\rceil- \tIP^{(L)}_j(t)$ & $\sbbB_i(t)-\sbB_i(t)$ \\ \hline
\end{tabular}
\end{center}
\end{table}
\fi

\subsubsection{Proof of Condition (\ref{eq:ipres}):}
\label{subsubsec:ipres}

For each component $j$ $(1 \le j \le n)$,  define
\begin{align}
\begin{split}
\label{eq:IPtracking}
&\mA=\aj, \quad \mK=\{k_j+1,\cdots,K\}, \quad s_l^{(L)}=\LL_l-\LL_{k_j}~(l \in \mK), \quad \tL_0=-\LL_{k_j},
\\	
&\sT(t)=\lceil\tIP^{(L)}_j(t)\rceil,~~t \ge \tL_0,
 \\
& \sW(t)=\sIP_j(t),~~~t \ge \tL_0, \qquad \mbox{and} \qquad \sW_0=-\aj \cdot \sBD(-L,-\LL_{k_j}). 
\end{split}
\end{align}

To verify this is an instance of the stochastic tracking model, by our replenishment policy in Table \ref{tbl:replenishment}, $\lceil\tIP^{(L)}_j(t)\rceil$ ($t \ge 0$) is a pure jump process.
Apply (\ref{eq:target}) to the $L^{th}$ system, 
\[
\tIP^{(L)}_j(t)=\sbY_j(t)-\mathbf A_j \cdot \sBD(t+\LL_{k_j}-L,t),~t \ge -\LL_{k_j}.
\]
Thus under (\ref{eq:iptargetcont}), the asymptotic Lipschitz condition (\ref{eq:Lcontinuity})-(\ref{eq:alcon}) is satisfied with $\sE(t)=1$. 
%If component $j$ has the longest lead time $L$, then $\tIP^{(L)}_j(t)=\sbY_j(t)$ $(t \ge -L)$ is constant. Otherwise, by induction, $\sbY_j(t)$ and thus $\tIP^{(L)}_j(t)$ $(t \ge -\LL_{k_j})$ are pure jump processes (NOTE:following the replenishment procedure in Table \ref{tbl:replenishment})
%given that the processes are determined by $\sbY_{j'}(t+\tL_{k_j}-\tL_{k_{j'}})$ for all $j'$ such that $\tL_{k_{j'}}>\tL_{k_j}$, and $\sBD(t+\LL_{k_j}-L,t)$, an increment of a compound Poisson demand process. 
%Let
%\[\sT_c(t)=\sT(t)=\tIP^{(L)}_j(t),\quad t \ge -\LL_{k_j}.\]
%Then $\sT_c(t)$ and $\sT(t)$ $(t \ge -\LL_{k_j})$ satisfy (\ref{eq:Lcontinuity}) given that (\ref{eq:iptargetcont}) holds. 

By the definition of $IP^-_j(t)$ in Table \ref{tbl:notations} and refer to (\ref{eq:ipcontrol}), in the $L^{th}$ system,
\begin{align*}
\begin{split}
& \sIP_j(t)=IP^{(L)-}_j(t) \vee \lceil \tIP^{(L)}_j(t)\rceil, \quad t  \ge -\LL_{k_j}, \\
\mbox{and} \quad &IP^{(L)-}_j(t)=\sIP_j(t^-)- \aj \cdot \sjD(t), \quad t > -\LL_{k_j},
\end{split}
\end{align*}
which fits the definition of the tracking process $\sW(t)$ $(t \ge \tL_0)$. The initial value $\sW_0$ is the inventory position of component $j$ before the placement of the first order at time $-\LL_{k_j}$.

\begin{corollary}
\label{cor:invtra}
For all $j=1,\cdots,n$, 
\begin{equation}
\label{eq:ipresco}
\lim_{L \rightarrow \infty} \E\left[\sup_{t \ge -\hLL_{k_j}} \left|\hIP_j(t)-\thIP_j(t)\right|\right]=\lim_{L \rightarrow \infty} \E\left[\sup_{t \ge -\hLL_{k_j}} \left|\frac{\sIP_j(Lt)-\lceil\tIP^{(L)}_j(Lt)\rceil}{\sqrt{L}}\right|\right]
=0.
\end{equation}
\end{corollary}

\smallskip
\noindent{\bf Proof:} The first equality follows from definitions in (\ref{eq:scaleinv}). 
%The second equality holds for all $j$ where $L_{k_j}=L$, in which case $\sIP_j(t)=\tIP^{(L)}_j(t)$ for all $t$. For all other components, 
Theorem \ref{thm:master} shows that the second equality holds if
%for $j$ ($\jrange$) if
\[
\lim_{L \rightarrow \infty} \E\left[\hG_0\right]\stackrel{\Delta}{=} \lim_{L \rightarrow \infty}\frac{\left(\sW_0-\sT(\tL_0)\right)^+}{\sqrt{L}} = 0,
\]
where by (\ref{eq:IPtracking}), for given $j$ ($\jrange$), 
\[
\hG_0 
%\stackrel{\Delta}{=} \frac{\left(\sW_0-\sT(\tL_0)\right)^+}{\sqrt{L}}
= \left(\frac{-\aj \cdot \sBD(-L,-\LL_{k_j})-\lceil\tIP^{(L)}_j(-\LL_{k_j})\rceil}{\sqrt{L}}\right)^+.
\]
Applying (\ref{eq:target}) to the $L^{th}$ system, 
\[
\tIP^{(L)}_j(-\LL_{k_j}) =\sbY_j(-\LL_{k_j})-\aj \cdot \sBD(-L,-\LL_{k_j}).
\]
Therefore by Lemma \ref{lem:initconv} with $\nu=\aj \cdot \boldsymbol \mu$, 
\[
\lim_{L \rightarrow \infty} \hG_0=\lim_{L \rightarrow \infty} \E\left[\frac{\left(-\sbY_j(-\LL_{k_j})\right)^+}{\sqrt{L}}\right] \le \lim_{L \rightarrow \infty} \E\left[\left(|\hbY_j(-\hLL_{k_j})|-\sqrt{L} \aj \cdot \boldsymbol \mu \right)^+\right]=0.~~~~\blacksquare
\]

\subsubsection{Proof of Condition (\ref{eq:blkres}):}
\label{subsubsec:blkres}

It is helpful to first review the processes that will be involved in the analysis.  In the $L^{th}$ system and for product $i$ $(1 \le i \le m$), $\sbB_i(t)$ $(t \ge 0)$ is the backlog target and $\sbbB_i(t)$ $(t \ge 0)$ is the ideal backlog level for reaching our SP-based lower bound. Both are determined by the LP in (\ref{eq:allspL}), using actual inventory positions and inventory position targets as their respectively inputs. The backlog level under our allocation policy is $\sB_i(t)$ $(t \ge 0)$. To help with the analysis, we will define below an auxiliary process $\nsB_i(t)$ $(t \ge 0)$, which mimics $\sB_i(t)$ $(t \ge 0)$, except that it is not allowed to rise above the target $\sbB_i(t)$ $(t \ge 0)$.

For each product $i$ $(1 \le i \le m)$, we construct an instance of the stochastic tracking model by letting
\begin{align}
\begin{split}
\label{eq:blktrack}
&\mA=\mathbf e_i,\quad \mK=\{1,\cdots,K\},\quad s^{(L)}_l = \LL_l~(l \in \mK), \quad \tL_0=0,
\\	
&\sT(t)=-\sbB_i(t),~~~t \ge 0,
\\
&\sW(t)=-\nsB_i(t),~~~t \ge 0, \qquad \mbox{and} \qquad \sW_0=-D^{(L)}_i(-L,0),
\end{split}
\end{align}
where $\nsB_i(t)$ $(t \ge 0)$ is defined by 
\begin{align*}
\begin{split}
& \underline{B}^{(L)-}_i(0)= -\sW_0=D^{(L)}_i(-L,0),
\\
& \nsB_i(t)=\underline{B}^{(L)-}_i(t) \wedge \sbB_i(t),\qquad t \ge 0,
\\
\mbox{and} \qquad & \underline{B}^{(L)-}_i(t) = \nsB_i(t^-)+\mA \sjD=\nsB_i(t^-)+d^{(L)}_i(t), \qquad t>0.
\end{split}
\end{align*}
Below we will show that $-\sbB_i(t)$ $(t \ge 0)$ is a target process, which, following the above definition, also implies that $-\nsB_i(t)$ $(t \ge 0)$ is a tracking process.

When defining our allocation policy in Section \ref{subsubsec:apolicydef}, we have shown that $\sbB_i(t)$ ($t \ge 0$) is a pure jump process. 
%Recall that we use (\ref{eq:setBtar}) to set backlog targets, so the value of $\sbB_i(t)$ is updated only at the occurrence of a change of the following balance level of some component $j$,
%\[Q^{(L)}_j(t)= \mathbf A_j \sBD(t-\LL_{k_j},t)-IP^{(L)}_j(t-\LL_{k_j}), \qquad 1 \le j \le n.\]
%Under our replenishment policy, $IP^{(L)}_j(t)$ $(t \ge -\LL_{k_j})$ is determined by its inventory position target $\tIP^{(L)}_j(t)$ $(t \ge -\LL_{k_j})$ and the demand process $\sD(t)$ $(t \ge -L)$, both of which are pure jump processes. Therefore $-\sbB_i(t)$ $(t \ge 0)$ is also a pure jump process. 
To show that $-\sbB_i(t)$ $(t \ge 0)$ is asymptotically Lipschitz continuous,
let 
\begin{equation}
\label{eq:deE}
\sE(t)=|\sbbB_i(t)-\sbB_i(t)|\quad t \ge 0.
\end{equation}
Then
\[
|\sbB_i(t_1)-\sbB_i(t_2)| \le |\sbbB_i(t_1)-\sbbB_i(t_2)|+ \sE(t_1)+\sE(t_2),
\]
which satisfies conditions (\ref{eq:Lcontinuity})-(\ref{eq:alcon}) because:
\begin{enumerate}
\item
By applying (\ref{eq:iQ}) to (\ref{eq:blktargetcont}),
there exists some constant $\kappa$ such that 
\begin{align*}
\begin{split}
 |\sbbB_i(t_1)-\sbbB_i(t_2)| \le \kappa & \sum_{j=1}^n \huge(|\mathbf A_j \sBD(t_2-L,t_2)-\mathbf A_j \sBD(t_1-L,t_1)|
\\
~~~~~&+|\mathbb Y^{(L)}_j(t_2-\LL_{k_j})-\mathbb Y^{(L)}_j(t_1-\LL_{k_j})|\huge).
\end{split}
\end{align*}
where under (\ref{eq:iptargetcont}), there also exists some constant $\kappa$ such that
\begin{align*}
\begin{split}
~&|\mathbb Y^{(L)}_j(t_2-\LL_{k_j})-\mathbb Y^{(L)}_j(t_1-\LL_{k_j})| 
\\
\le & \kappa \sum_{k > k_j} ||\sBD(t_2-\LL_k,t_2-\LL_{k_j})-\sBD(t_1-\LL_k,t_1-\LL_{k_j})||_1 \\
\le & \kappa \sum_{k > k_j} \huge(||\sBD(t_2-\LL_k,t_2)-\sBD(t_1-\LL_k,t_1)||_1
+||\sBD(t_2-\LL_{k_j},t_2)-\sBD(t_1-\LL_{k_j},t_1)||_1\huge)
\end{split}
\end{align*}
\item
By (\ref{eq:btarclose}), there exists some constant $\kappa$ such that
\[
|\sE(t)| \le \kappa \sum_{j=1}^n |\mathbb Q^{(L)}_j(t)- Q^{(L)}_j(t)|=\kappa \sum_{j=1}^n |IP^{(L)}_j(t-\LL_{k_j})-\tIP^{(L)}_j(t-\LL_{k_j})|, \quad t \ge 0.
\]
Thus (\ref{eq:alcon}) follows directly from Corollary \ref{cor:invtra}.
\end{enumerate}

Having shown that $(-\sbB_i(t),-\nsB_i(t))$ $(t \ge 0)$ is an instance of the target and tracking processes, we now can prove (\ref{eq:blkres}) as the following corollary to Theorem \ref{thm:master}. 
\begin{corollary}
	\label{cor:backlog}
	For $i=1,\cdots, m$, let
	\[
	\hnsB_i(t)=\frac{\nsB_i(Lt)}{\sqrt{L}}.
	\]
	Then
	\begin{equation}
	\label{eq:blkresco}
	\lim_{L \rightarrow \infty}\E\left[\sup_{t \ge 0}\left|\hbB_i(t)-\hnsB_i(t)\right|\right]= 0,
\end{equation}
which implies Condition (\ref{eq:blkres}). 
\end{corollary}

\smallskip
\noindent{\bf Proof}:
By Theorem \ref{thm:master}, (\ref{eq:blkresco}) holds if 
\begin{equation}
\label{eq:initblk}
\lim_{L \rightarrow \infty} \E\left [\left|\hbB_i(0)-\hnsB_i(0)\right|\right]=0.
\end{equation}

To prove this initial condition, observe that
\[
\sbB_i(0)-\nsB_i(0)=\sbB_i(0)-\sbB_i(0) \wedge B^{(L)-}_i(0)=\left(\sbB_i(0)-B^{(L)-}_i(0)\right)^+= \left(\sbB_i(0)-D^{(L)}_i(-L,0) \right)^+. 
\] 
Since $\sbB_i(0)$ minimizes (\ref{eq:setBtar}) with $\mathbf Q(t)=\mathbf Q^{(L)}(0)$, there exists some constant $\kappa$,  which depends only on matrix $A$, such that 
\[
\sbB_i(0) \le \kappa \sum_{j=1}^n \left[Q^{(L)}_j(0)\right]^+=\kappa \sum_{j=1}^n \left[\aj \cdot\sBD(-L^{(L)}_{k_j},0)-\sIP_j(-L^{(L)}_{k_j})\right]^+.
\]
Since $\tIP^{(L)}_j(-L^{(L)}_{k_j}) \le \sIP_j(-L^{(L)}_{k_j})$ $(\jrange)$ under our replenishment policy (see (\ref{eq:ipcontrol})),  
\begin{align*}
\begin{split}
\sbB_i(0)-D^{(L)}_i(-L,0) 
& \le \kappa \sum_{j=1}^n \left[\aj \cdot\sBD(-L^{(L)}_{k_j},0)-\tIP^{(L)}_j(-L^{(L)}_{k_j})\right]^+-D^{(L)}_i(-L,0)
\\
&= \kappa \sum_{j=1}^n \left[\aj \cdot\sBD(-L,0)-\sbY_j(-L^{(L)}_{k_j})\right]^+-D^{(L)}_i(-L,0).
\end{split}
\end{align*}

Let $\nu=\mu_i/(\kappa n(m\bar{a}+1)+1)$. By its definition, $\nsB_i(0)=\sbB_i(0) \wedge D^{(L)}_i(-L,0)$. Therefore, 
\begin{align*}
\begin{split}
&~\left|\hbB_i(0)-\hnsB_i(0)\right|
\\
&\le \frac{1}{\sqrt{L}} \left(\kappa\sum_{j=1}^n \left[\aj\cdot \sBD(-L,0)-\mathbb Y_j^{(L)}(-L^{(L)}_{k_j}) \right]^+-D_i^{(L)}(-L,0)\right)^+
\\
& \le \left(\kappa \sum_{j=1}^n \aj\cdot |\hBD(-1,0)|+ \kappa \sum_{j=1}^n |\hat{\mathbb Y}_j^{(L)}(-\hat{L}_{k_j})| -\sqrt{L} \mu_i-\hat{D}_i^{(L)}(-1,0)\right)^+
\\
& \le  n \kappa \bar{a}\sum_{i'=1}^m \left(|\hat{D}_{i'}(-1,0)|-\sqrt{L}\nu\right)^++ \kappa \sum_{j=1}^n \left(|\hat{\mathbb Y}_j^{(L)}(-\hat{L}_{k_j})|-\sqrt{L}\nu\right)^+
+\left(\hat{D}_i^{(L)}(-1,0)-\nu\right)^+.
\end{split}
\end{align*}
Since $\mD^{(L)}(t)$ $(t \ge -L$) is stationary, Lemmas \ref{lem:NBD} and \ref{lem:initconv} imply (\ref{eq:initblk}), which proves (\ref{eq:blkresco}).

To show that (\ref{eq:blkresco}) implies (\ref{eq:blkres}),
\[
|\hB_i(t)-\hbsB_i(t)| \le |\hB_i(t)-\hbB_i(t)|+|\hbB_i(t)-\hbsB_i(t)|.
%=|\hB_i(t)-\hbB_i(t)|+|\hat{\mE}^{(L)}(t)|.
\]
By (\ref{eq:deE}), the last term is simply $\hE(t)$ that satisfies (\ref{eq:alcon}), thus we only need to prove
% $\sbB_i(t)$ $(t \ge 0)$ is asymptotically Lipschitz continuous with $\sT_c(t)=\sbbB_i(t)$ $(t\ge 0)$, the last term is negligible, and we only need to prove that
\begin{equation}
\label{eq:laststep}
\lim_{L\rightarrow \infty} \E\left[\sup_{t \ge 0} \left|\hB_i(t)-\hbB_i(t)\right|\right]=0.
\end{equation}
As is shown in (\ref{eq:Z-b}), under our allocation policy,
\[
\sB_{i} (t) \ge B^{(L)-}_{i} (t) \wedge \sbB_{i} (t) =\nsB_{i}(t),\qquad t \ge 0,
\]
and thus 
\[
\left(\hbB_{i}(t)-\hB_{i}(t)\right)^+ \le \left(\hbB_{i}(t)-\hnsB_{i}(t)\right)^+=\left|\hbB_{i}(t)-\hnsB_{i}(t)\right|, \quad t \ge 0. 
\]
Moreover, since property (\ref{eq:twoside}) applies when $\sB_i(t) - \sbB_i(t) \ge 1$, in all cases,
\[
\hB_i(t)-\hbB_i(t) \le \frac{1}{\sqrt{L}}+\frac{\bar{a}}{\underline{a}} \sum_{i'\ne i} \left(\hbB_{i'} (t)-\hB_{i'} (t)\right)^+
\le \frac{1}{\sqrt{L}}+\frac{\bar{a}}{\underline{a}} \sum_{i'\ne i} \left|\hbB_{i'} (t)-\hnsB_{i'} (t)\right|,\quad t\ge 0.
\]
Combine the above two inequalities,
\[
\left|\hB_i(t)-\hbB_i(t)\right| \le \frac{1}{\sqrt{L}}+ \frac{\bar{a}}{\underline{a}} \sum_{i'=1}^m \left|\hbB_{i'} (t)-\hnsB_{i'} (t)\right|,\quad t\ge 0,
\]
so (\ref{eq:laststep}) follows immediately from (\ref{eq:blkresco}). $\qquad \qquad \qquad \qquad \qquad \blacksquare$

\section{Stability of SP Optimal Solution}
\label{sec:SPsolution}

In this section, we finalize our analysis by showing that with proper treatments of the SP (\ref{eq:MSP2}), conditions (\ref{eq:iptargetcont}), (\ref{eq:blktargetcont}), and (\ref{eq:btarclose}) can be satisfied. To this end, we need to address two issues: uniqueness of the optimal solutions and values of Lipschitz constants. 

Recall that $\mathbb B_i(t)$ and $\mathbb B_i^*(t)$ $(t \ge 0$) are optimal solutions to the same LP (\ref{eq:miniB}) with different RHS coefficients in constraints. Hence by Hoffman's Lemma, for each $t$ ($t \ge 0$), there exist $\mathbb B_i(t)$ and $\mathbb B^*_i(t)$ that satisfy (\ref{eq:btarclose}), and for each $t_1$ and $t_2$ $(t_2>t_1 \ge 0$), there exist $\mathbb B_i^*(t_1)$ and $\mathbb B_i^*(t_2)$ that satisfy (\ref{eq:blktargetcont}). Nevertheless, the lemma does not exclude the possibility that if the optimal solution of (\ref{eq:miniB}) is not unique, then to satisfy (\ref{eq:blktargetcont}) at the same $t_1$, we may need different $\mathbb B_i^*(t_1)$ for different $t_2$.  In this case, a well-defined process $\mathbb B_i^*(t)$ $(t \ge 0)$ that satisfies (\ref{eq:blktargetcont}) for all $t_1$ and $t_2$ ($0 \le t_1 <t_2$) is not guaranteed. To avoid this situation, we use perturbation to keep the optimal solution of (\ref{eq:miniB}) unique. As a result,  both $\mathbb B_i(t)$  and $\mathbb B_i^*(t)$ $(t \ge 0)$ are uniquely-defined processes that satisfy (\ref{eq:btarclose}) and (\ref{eq:blktargetcont}).  
 
It takes significantly more effort to address (\ref{eq:iptargetcont}). Because the supports of $\BD^k$  $(1 \le k \le K)$ are unbounded, $\varphi^k(\y^{k+1},\cdots,\y^K,\x)$ in (\ref{eq:MSP2}) are infinite-dimensional problems, so their optimal solutions may not satisfy the continuity condition specified by (\ref{eq:iptargetcont}). We develop finite-dimensional LPs to approximate the latter problem and prove that we can always keep the approximation error negligible by keeping the dimension of the LP sufficiently high. We also use perturbation to maintain uniqueness of the optimal solution. Both the approximation and perturbation (which also addresses uniqueness of the optimal solution to (\ref{eq:miniB})) are developed in Section \ref{subsec:finitedef}. 

In the asymptotic analysis, when lead times increase, probabilities of having larger sample values of $\BD^k(t)$ $(1 \le k \le K)$ increase, so the dimension of the approximating LP needs to increase to keep the approximation sufficiently accurate. To sustain (\ref{eq:iptargetcont}), we need to rule out the possibility that the Lipschitz constant $\kappa$ has to grow unboundedly with the problem dimension. In Section \ref{subsec:FCBM}, we show that $\kappa$ can be kept at a finite value regardless how large the dimension of the LP becomes.

\subsection{Approximation and Perturbation}
\label{subsec:finitedef}
We develop finite-dimension approximations to SP (\ref{eq:MSP2}) by taking following steps: let $\mathbf M$ be a $m$-dimensional vector with all entries equal to an integer $M>0$. 
Denote $\BD^k \wedge \mathbf M$ by $\BD_M^k$ ($1 \le k \le K$).
Replace $\BD^k$ in (\ref{eq:MSP2}) with $\BD^k_M$ in $\varphi^k(\y^{k+1},\cdots,\y^K,\x)$ ($1 \le k \le K$) to define
\begin{align}
\label{eq:MSPapp}
\begin{split}
\varphi^K_M &=\min_{\y^K \in \mathbb R^{n_{\scaleto{K}{2.5pt}}}} \left \{\mathbf h^K \cdot \y^K+\E\left[\varphi^{k-1}_M(\y^K,\BD^K_M)\right]\right \},
\\
\varphi^k_M(\y^{k+1},\cdots,\y^K,\x) &=\min_{\y^k \in \mathbb R^{n_{\scaleto{k}{2.5pt}}}} \left \{\mathbf h^k \cdot \y^k+\E\left[\varphi^{k-1}_M(\y^k,\cdots,\y^K,\x+\BD^k_M)\right]\right \},~1 \le k <K,
\\
\varphi_M^0(\y^1,\cdots,\y^K,\x) &=\varphi^0(\y^1,\cdots,\y^K,\x) = -\max_{\z \in \mathbb R^m} \left \{\mathbf c \cdot \z|\z \le \x, A^k \z \le \y^k,1\le k \le K\right \}.
\end{split}
\end{align}

Since $\BD^k$ are integers, for given $M$, $\BD^k_M$ have a finite support ($ 1 \le k \le K$), so $\varphi^k_M(\y^{k+1},\cdots,\y^K,\x)$ ($0 \le k \le K$) are finite-dimension problems. There is no approximation error when $k=0$, in which case the original problem is a finite-dimension LP.  For $k>0$, the theorem below shows that the approximation error can be kept negligible if $M$ is sufficiently large.  

\begin{theorem}{(see Appendix A for the proof)}
\label{thm:unicom}
As $M \rightarrow \infty$, $\varphi^k_M(\y^{k+1},\cdots, \y^{K},\x)$ converges uniformly to $\varphi^k(\y^{k+1},\cdots,\y^K,\x)$ 
over $\y^{k'} \in \mathbb R^{k'}$ ($k+1 \le k' \le K$) and $\x \in \mathbb R^m$. 
\end{theorem}

For each $k$ $(1 \le k \le K)$, $\varphi_M^k(\y^{k+1},\cdots,\y^K,\x)$ in (\ref{eq:MSPapp}) can be equivalently formulated as an LP. The reformulation, which involves characterizing a multi-stage decision problem by a scenario tree, is thoroughly elaborated on in Section 3.1.3 of \cite{Shapiro2009}. Appendix C describes the detailed steps of specializing this standard process to (\ref{eq:MSPapp}). The resulting LP formulation is: 
\begin{equation}
\label{eq:finitLPobj}
\varphi_M^k(\y^{k+1},\cdots,\y^K,\x)=\min_{\y^k,\cdots,\y^1, \mathbf z} \left \{\mathbf h^k \cdot \y^k + \sum_{k'=1}^{k-1}  \sum_{\bar{\omega} \in \Omega_k^{k'}}  \mathbb P (\bar{\omega}) \mathbf h^{k'} \cdot \y^{k'}(\bar{\omega})
 -\sum_{\bar{\omega} \in \Omega_k^0} \mathbb P(\bar{\omega}) \mathbf c \cdot \z(\bar{\omega}) \right \}
\end{equation}
subject to: 
\begin{eqnarray}
\label{eq:mapping1}
\z(\bar{\omega}) ~~~~~~~~~~~~~~~~~~~~~~~~~~~~~~&\le&~~~ \x+\underline{\BD}^k_M (\bar{\omega}),~~~\bar{\omega} \in \Omega_k^0,
\\
\label{eq:mapping2}
A^{k'} \mathbf z(\bar{\omega}) ~~~- \y^{k'}(\bar{\omega}')~~~~~~~~~~~~~~&\le&~~~ \mathbf 0,~~~~~~~~~~~~~~~\bar{\omega} \in \Omega_k^0,~\bar{\omega}' \in \Omega_k^{k'},~\bar{\omega}' \sqsubset \bar{\omega}, ~1 \le k' <k,
\\
\label{eq:mapping3}
A^k \mathbf z(\bar{\omega}) ~~~~~~~~~~~~~~~~~~~~~- \y^k~~~ &\le&~~~ \mathbf 0,~~~~~~~~~~~~~~~\bar{\omega} \in \Omega_k^0,
\\
\label{eq:mapping4}
A^{k'} \mathbf z(\bar{\omega}) ~~~~~~~~~~~~~~~~~~~~~~~~~~~~~~ &\le&~~~ \y^{k'},~~~~~~~~~~~~~\bar{\omega} \in \Omega_k^0, ~k < k' \le K.
\end{eqnarray}
Here $\Omega_k^{k'}$ is the set of strings that encode sample paths $(\BD_M^{k'+1},\cdots,\BD_M^k)$, $0 \le k' < k \le K$. For each sample path $\bar{\omega} \in \Omega_k^0$,  
\[
\underline{\BD}_M^k(\bar{\omega})=\BD^k_M(\bar{\omega})+...+\BD^1_M(\bar{\omega}), \qquad 1 \le k \le K. 
\]
The probability attached to sample path $\bar{\omega}$ is denoted by $\mathbb P(\bar{\omega})$. For $\omega' \in \Omega_{k}^{k'}$ and $\bar{\omega} \in \Omega_k^0$ ($0<k' <k, 1< k \le K$), we write $\bar{\omega}' \sqsubset \bar{\omega}$ if the sample path $(\BD^{k'+1}_M,\cdots,\BD^k_M)$ encoded by $\bar{\omega}'$ is on the same sample path $(\BD^1_M,\cdots,\BD^k_M)$ encoded by $\bar{\omega}$. 

We perturb parameter values in the objective function (\ref{eq:finitLPobj}) to keep the optimal solution of the LP unique. Because $A^k$ and $\BD^k_M$ $(1 \le k \le K)$ are integer-valued, the optimal solution must be rational numbers when $k=K$, and by induction, so are the optimal solutions when $k=K-1,\cdots,1$. By adding tiny different irrational values to $\mathbb P(\bar{\omega})$ ($\bar{\omega} \in \Omega_k^l, 1 \le l \le k \le K)$), $\mathbf h^k$ $(1 \le k \le K)$, and $\mathbf c$, we can exclude cases in which two different solutions yield the same objective value. 

Again, we refer readers to Appendix C for the details of the perturbation. To incorporate the above developments into our policy: in the allocation policy defined in Table \ref{tbl:allproc}, (\ref{eq:setBtar}) is equivalent to $\varphi^0_M(\cdot)$ in (\ref{eq:MSPapp}) for all $M$, so we can solve it with the same perturbed values of $\mathbf c$ to ensure the optimal solution is unique, in which case Hoffman's lemma immediately leads to:

\begin{lemma}
\label{lem:Bcont}
Let $\mathbb B(t)$ be the (unique) optimal solution of (\ref{eq:setBtar}) and $\mathbb B^*(t)$ be the (unique) optimal solution of the same LP with $\mathbf Q^k(t)$ replaced by $\mathbb Q^k(t)$ $(1 \le k \le K, t\ge 0$). Then $\mathbb B(t)$ $(t \ge 0)$ and $\mathbb B^*(t)$ $(t \ge 0)$ satisfy (\ref{eq:btarclose}) and $\mathbb B^*(t)$ $(t \ge 0)$ satisfies (\ref{eq:blktargetcont}).
 \end{lemma}

In the replenishment policy defined in Table \ref{tbl:replenishment}, instead of (\ref{eq:Yset}), we solve its approximating LP (\ref{eq:finitLPobj})-(\ref{eq:mapping4}) at Step 2 (a) to set inventory position targets, using 
\begin{equation}
\label{eq:LPinputs}
\y^{k'}=\mathbb Y^{k'}(t+L_k-L_{k'}) ~(k<k' \le K) \quad \mbox{and} \quad \x=\BD(t+L_k-L_K,t).
\end{equation}
as inputs, a sufficiently large $M$ to keep the approximation error negligible, and perturbed parameter values in the objective function to keep the optimal solution unique. 
%What remains to be shown, as we will do next, is that there does exist a Lipschitz constant $\kappa$ that is independent of $M$ to satisfy (\ref{eq:iptargetcont}). 

To ensure that the above approximation and perturbation will not introduce errors that can compromise asymptotic optimality, we can fix some constant $\epsilon$. In the asymptotic analysis, for the $L^{th}$ system, Theorem \ref{thm:unicom} indicates that we can always choose a sufficiently large $M$ to keep 
\[
|\varphi^{K(L)}-\varphi_M^{K(L)}| \le \epsilon,
\]
where $\varphi^{K(L)}$ and $\varphi_M^{K(L)}$ are respectively optimal objective values of (\ref{eq:MSP2}) and (\ref{eq:MSPapp}) specialized to the $L^{th}$ system. 
As $M$ increases with $L$, $\Omega_k^l$ in (\ref{eq:finitLPobj})-(\ref{eq:mapping4}) will have more elements $\bar{\omega}$ ($0 \le l \le K$, $1 \le k \le K$), but there is an infinite number of irrational values we can use to perturb $\tilde{\mathbb P}(\bar{\omega})$ (see Appendix C for details of choosing these values). 
By increasing the integer divisor (again, see Appendix C) as $L$ and $M$ grow,
we can  keep the latter values small enough  so the perturbation will not change the optimal objective value of (\ref{eq:MSPapp}) by more than $\epsilon$. Under these arrangements, we can set inventory positions by solving the LP (\ref{eq:finitLPobj})-(\ref{eq:mapping4}) and backlog targets by solving the LP (\ref{eq:setBtar}) with perturbed parameter values. In comparison with targets set by solving the lower-bound SP (\ref{eq:MSP2}), the difference between the average inventory costs is no more than $2\epsilon$. This error is invisible on the diffusion scale. 

\subsection{Values of Lipschitz Constants}
\label{subsec:FCBM}

The approximation and perturbation scheme in Section \ref{subsec:finitedef} is developed to set inventory position targets that satisfy (\ref{eq:iptargetcont}) and backlog targets that satisfy (\ref{eq:btarclose})-(\ref{eq:blktargetcont}). In the asymptotic analysis, the LP (\ref{eq:setBtar}) stays the same when $L$ changes, so (\ref{eq:btarclose})-(\ref{eq:blktargetcont}) are satisfied as long as the optimal solution is kept unique by the aforementioned perturbation.  On the other hand, as $L$ increases, larger $M$ is needed in (\ref{eq:finitLPobj})-(\ref{eq:mapping4}) for an accurate approximation of the demand distribution. For (\ref{eq:iptargetcont}) to hold, the value of $\kappa$ needs to be independent of the change of $M$. Theorem \ref{thm:SPconti} and its corollary below show that such value always exists.
%In our asymptotic analysis, the same value of $\kappa$ is used in (\ref{eq:iptargetcont}) for all systems to bound the change of the optimal solution relative to the change of RHS of constraints. As lead times increase, larger $M$ should be used in (\ref{eq:finitLPobj})-(\ref{eq:mapping4}) to keep the approximation error negligible. Hence it is necessary to show, as in Theorem \ref{thm:SPconti} below, that the change of the optimal solution of (\ref{eq:finitLPobj})-(\ref{eq:mapping4}) is bounded by a constant that is independent of $M$. As a corollary, the constant can be used as $\kappa$ in  (\ref{eq:iptargetcont}). 

\begin{theorem}
\label{thm:SPconti}
	For $k=1,\cdots, K$, let $\y^{k*}_{M,a}$ and $\y^{k*}_{M,b}$ be (unique) optimal solutions to (\ref{eq:finitLPobj})-(\ref{eq:mapping4}) under inputs $(\y_a^{k+1},\cdots, \y_a^{K},\x_a)$ and $(\y_b^{k+1},\cdots, \y^{K}_b,\x_b)$ respectively. Then 
		\begin{equation}
		\label{eq:LPconti}
		||\y_{M,a}^{k*}-\y_{M,b}^{k*}||_\infty \le \kappa \left(||\y_a^{k+1}-\y_b^{k+1}||_\infty+\cdots + ||\y_a^K-\y_b^K||_\infty+||\x_a-\x_b||_\infty\right),
		\end{equation}
		where $\kappa$ depends only on $(m,n_1,\cdots,n_K)$ and  $A^k$ $(1\le k \le K)$. 
\end{theorem}

Before proving the theorem, we first show its implication by deriving a corollary: as is specified in Section \ref{subsec:finitedef}, we solve (\ref{eq:finitLPobj})-(\ref{eq:mapping4}) to obtain values of $\mathbb Y^k(t)$ ($t \ge -L_k, 1 \le k \le K$) for implementing our replenishment policy. For this purpose, $\mathbb Y^k(t)$ $(t \ge -L_k)$ needs to satisfy (\ref{eq:iptargetcont}), which is clearly satisfied when $k=K$, in which case these values are kept constant. Theorem  \ref{thm:SPconti} allows us to use induction to extend (\ref{eq:iptargetcont}) to $k<K$, which is stated as the following corollary:

\begin{corollary}
	\label{cor:iptargcont}
	For $t \ge 0$, let $\mathbb Y^k(t)$ $(t \ge -L_k, 1 \le k \le K)$ be the (unique) optimal solution of $\varphi^k_M(\y^{k+1},\cdots,\y^K,\x)$ defined in (\ref{eq:finitLPobj})-(\ref{eq:mapping4}), with inputs given by (\ref{eq:LPinputs}).  Then (\ref{eq:iptargetcont}) holds for a constant $\kappa$ that depends only on $(m,n_1,\cdots,n_K)$ and  $A^k$ $(1\le k \le K)$. 
\end{corollary}

\smallskip
\noindent{\bf Proof of Corollary \ref{cor:iptargcont}}: For given $k$, suppose that (\ref{eq:iptargetcont}) holds for all $l>k$ (which is the case for $k=K-1$). This means that there exists some $\kappa$, depending only on $(m,n_1,\cdots,n_K)$ and  $A^k$ $(1\le k \le K)$, such that for all $-L_l \le t_1<t_2$, 
\begin{align}
\label{eq:inducassum}
\begin{split}
&||\mathbb Y^l(t_2-(L_l-L_k))-\mathbb Y^l(t_1-(L_l-L_k))||_1 
\\
\le & \kappa \sum_{k'=l+1}^K ||\BD(t_2-(L_{k'}-L_k),t_2-(L_l-L_k))-\BD(t_1-(L_{k'}-L_k),t_1-(L_l-L_k))||_1,\\
=& \kappa \sum_{k'=l+1}^K ||\left\{\BD(t_2-(L_{k'}-L_k),t_2)-\BD(t_2-(L_l-L_k),t_2)\right\}
\\
& \qquad \qquad -\left\{\BD(t_1-(L_{k'}-L_k),t_1)-\BD(t_1-(L_l-L_k),t_1)\right\}||_1.
\end{split}
\end{align}
Since $\mathbb Y^k(t)$ ($t \ge -L_k$) is obtained by solving (\ref{eq:finitLPobj})-(\ref{eq:mapping4}) using (\ref{eq:LPinputs}) for inputs $\y^l$ ($k<l \le K$) and $\x$, by Theorem \ref{thm:SPconti}, 
\begin{align}
\begin{split}
\label{eq:SPLip}
||\mathbb Y^k(t_2)-\mathbb Y^k(t_1)||_1 ~\le~ & \kappa \sum_{l=k+1}^K ||\mathbb Y^l(t_2-(L_l-L_k))-\mathbb Y^l(t_1-(L_l-L_k))||_1 
\\ 
&~+ \kappa ||\mathbf D(t_2-(L_K-L_k),t_2)-\mathbf D(t_1-(L_K-L_k),t_1)||_1
\\
\le ~& \kappa (K-k) \sum_{l=k+1}^K ||\BD(t_2-(L_l-L_k),t_2)-\BD(t_1-(L_l-L_k),t_1)||_1,
\end{split}
\end{align}
which shows that (\ref{eq:iptargetcont}) holds for $k$. $\quad \blacksquare$

The remainder of this section is dedicated to proving Theorem \ref{thm:SPconti}. We start from the following Theorem 2.4 in \cite{Mangasarian1987} (for convenience, we denote $\mathbf u$ and $\mathbf v$ in their theorem by $\bar{\mathbf u}$ and $\bar{\mathbf v}$ respectively here): 

\emph{Let $\beta \ge 1$ and $\beta^*=1/(1-1/\beta)$ (so that $||\x||_{\beta^*}$ is the dual norm of $||\x||_\beta$).
% be the $l_\beta$ norm of a vector $\x$. Let  and thus   $\x$. 	
Let the linear program
\[
\max_{\mathbf x} \{\mathbf p \cdot \mathbf x ~~~~\mbox{s.t.}~~ A \mathbf x \le \mathbf b,~C \mathbf x = \mathbf d \}
\]
have non-empty solution sets $\mathcal S^1$ and $\mathcal S^2$ for right-hand sides $(\mathbf b^1,\mathbf d^1)$ and $(\mathbf b^2,\mathbf d^2)$, respectively. 
For each $\mathbf x^1 \in \mathcal S^1$, there exists $\mathbf x^2 \in \mathcal S^2$ such that
\[
||\mathbf x^1-\mathbf x^2||_\infty  \le \nu_{\beta}(A;C) \left | \left |
\begin{array}{c}
\mathbf b^1-\mathbf b^2
\\
\mathbf d^1 -\mathbf d^2
\end{array} 
\right| \right |_\beta
\]
where 
\[
\nu_\beta(A;C):=\sup_{\bar{\mathbf u},\bar{\mathbf v}} \left \{ \left |\left | \begin{array}{c}
\bar{\mathbf u}
\\
\bar{\mathbf v}
\end{array} 
\right| \right |_{\beta^*}
\left |
\begin{array}{l}
||\bar{\mathbf u} A+\bar{\mathbf v} C||_1=1
\\
\mbox{Rows of } \left (\begin{array}{c} A \\ C\end{array}\right) \mbox{ corresponding to nonzero}
\\
\mbox{elements of } \left (\begin{array}{c} \bar{\mathbf u} \\ \bar{\mathbf v} \end{array}\right)\mbox{ are linearly independent}
\end{array}
\right .
\right\}.
\]
}

To apply the theorem to $\varphi_M^k(\y^{k+1},\cdots,\y^K,\x)$ defined in (\ref{eq:finitLPobj})-(\ref{eq:mapping4}), let $\mathbb A^k$ be the LHS constraint matrix of the latter LP. Since there is no equality constraint, we can ignore $\bar{\mathbf v}$ and $C$. Let $\beta=\infty$, and thus $\beta^*=1$. Then $\nu_\beta(A;C)$ and $\bar{\mathbf u}$ in the theorem specialize to
\[
\nu_\infty (\mathbb A^k):=\sup_{\bar{\mathbf u}} ||\bar{\mathbf u}||_1,
\]
where 
\[
||\bar{\mathbf u} \mathbb A^k ||_1=1,
\]
and rows of $\mathbb A^k$ corresponding to nonzero elements of $\bar{\mathbf u}$ are linearly independent.

The lemma below builds a critical connection between the formulation of the above quantities and the conclusion of Theorem \ref{thm:SPconti}.

\begin{lemma}(see Appendix A for the proof)
\label{lem:bondU}
	For $k=1,\cdots, K$, let $\hat{\mathbb A}^k$ be a maximal nonsingular submatrix of $\mathbb A^k$. Let $\mathbf u$ be a vector that has the same number of components as the number of rows in $\hat{\mathbb A}^k$. Then there exists a constant $\kappa$, depending only on $(m,n_1,\cdots,n_K)$ and the values of the components of $A^k$ $(1\le k \le K)$, such that 
	\begin{equation}
	\label{eq:limitnorm}
	||\mathbf u||_1 \le \kappa ||\mathbf u \hat{\mathbb A}^k||_1.
	\end{equation}
\end{lemma}

Proving this lemma is quite involved because we need to exploit the particular structure of $\mathbb A^k$ to show that the result, which does not apply to an arbitrary matrix, is true here. As may be seen from the proof of this lemma in Appendix A, it takes some effort to define this structure, especially for cases where $k>1$.  To help with the understanding of this result, we  explain the intuition behind the lemma's conclusion, after we first use the lemma to give the following proof of Theorem \ref{thm:SPconti}:

\smallskip
\noindent{\bf Proof of Theorem \ref{thm:SPconti} }: Let $\bar{\mathbf u}$ be any vector such that $||\bar{\mathbf u} \mathbb A^k||_1=1$ and rows of $\mathbb A^k$ corresponding to nonzero elements of $\bar{\mathbf u}$ are linearly independent.  Let $\hat{\mathbb A}^k$ be a maximal nonsingular submatrix of $\mathbb A^k$ that contains all these independent rows. Let $\mathbf u$ be a sub-vector of $\bar{\mathbf u}$ that is composed of components corresponding to rows in $\hat{\mathbb A}^k$ in $\bar{\mathbf u} \mathbb A^k$, and thus includes all nonzero elements of $\bar{\mathbf u}$. 

By Lemma \ref{lem:bondU}, there exists a constant $\kappa$, depending only on $(m,n_1,\cdots,n_K)$ and component values of $A^{k'}$ $(1 \le k' \le K)$, such that 
\[
||\bar{\mathbf u}||_1=||\mathbf u||_1 \le  \kappa ||\mathbf u \hat{\mathbb A}^k||_1=  \kappa ||\bar{\mathbf u} \mathbb A^k||_1=\kappa,
\]
and thus 
\[
\nu_\infty(\mathbb A^k)= \sup_{\bar{\mathbf u}} ||\bar{\mathbf u}||_1 \le \kappa.
\] 
By Theorem 2.4 in \cite{Mangasarian1987}, $\nu_\infty(\mathbb A^k)$ satisfies (\ref{eq:LPconti}), so does $\kappa$. $\qquad \blacksquare$

To explain the insights on proving Lemma \ref{lem:bondU}, we start from the identity
\begin{equation}
\label{eq:recover}
\mathbf u = \mathbf u \hat{\mathbb A}^k (\hat{\mathbb A}^k)^{-1}, \quad \mbox{i.e.,} \quad
\mathbf u^T=[(\hat{\mathbb A}^k)^{-1}]^T [\mathbf u \hat{\mathbb A}^k]^T,
\end{equation}
where $T$ denotes transpose. 
Following standard definitions \cite{Meyer2000},
the $L_1$ norm of a matrix $A$ is 
\[
||A||_1=\sup_{\x \ne 0}\frac{||A\x||_1}{||\x||_1},
\]
and its value is the maximum absolute column sum of $A$. This is also the value of $||A^T||_\infty$, 
%since the infinity norm of a matrix is simply 
which is the maximum absolute row sum of that matrix. 
Applying definitions of these norms to (\ref{eq:recover}) with $[(\hat{\mathbb A}^k)^{-1}]^T$ in place of $A$,
\[
||(\hat{\mathbb A}^k)^{-1}||_\infty=||[(\hat{\mathbb A}^k)^{-1}]^T||_1
=\sup_{\x \ne 0} \frac{||[(\hat{\mathbb A}^k)^{-1}]^T \x||_1}{||\x||_1} \ge 
\frac{||[(\hat{\mathbb A}^k)^{-1}]^T [\mathbf u \hat{\mathbb A}^k]^T||_1}{||[\mathbf u \hat{\mathbb A}^k]^T||_1}=\frac{||\mathbf u||_1}{||\mathbf u \hat{\mathbb A}^k||_1}.
\]
Thus we can let $\kappa$ in (\ref{eq:limitnorm}) be the largest $||(\hat{\FCBA}^k)^{-1}||_\infty$ over all $\hat{\FCBA}^k$ (nonsingular submatrices of $\FCBA^k$). What needs to be explained then is why this value of $\kappa$ can stay bounded as $M$ increases, which leads to more elements in $\Omega_k^0$ and $\Omega_k^{k'}$ ($k<k' \le K$) and thus a larger $\FCBA^k$.

\iffalse 
the identity 
\begin{equation}
    \label{eq:recover}
\mathbf u = \mathbf u \hat{\mathbb A}^k (\hat{\mathbb A}^k)^{-1}
\end{equation}
implies that 
\[
||\mathbf u||_1 \le || \mathbf u \hat{\mathbb A}^k||_1 || (\hat{\mathbb A}^k)^{-1}||_\infty,
\]
where the infinity norm $|| (\hat{\mathbb A}^k)^{-1}||_\infty$ is obtained by adding absolute values of entries on each row of $(\hat{\mathbb A}^k)^{-1}$ and taking the maximum value of these sums.  We let $\kappa$ in (\ref{eq:limitnorm}) be the largest $||(\hat{\FCBA}^k)^{-1}||_\infty$ over all $\hat{\FCBA}^k$ (nonsingular submatrices of $\FCBA^k$). What needs to be explained then is why the value of $\kappa$ can stay bounded as $M$ increases, which leads to more elements in $\Omega_k^0$ and $\Omega_k^{k'}$ ($k<k' \le K$) and thus a larger $\FCBA^k$.
\fi

Our explanation focuses on $k=1$, in which case (\ref{eq:mapping2}) is irrelevant and
\[
\FCBA^1=\left(\begin{array}{ccccc} \FCBH^1 & & & & \\  & \ddots & & & \FCBE^1 \\ & & & \FCBH^1 & \end{array} \right),
\]
where $\FCBH^1$ is the coefficient matrix of $\z(\bar{\omega})$ for a given $\bar{\omega}$ and $\FCBE^1$ is the coefficient matrix of $\y^1$. Submatrix $\FCBH^1$ has $m$ columns, corresponding to the number of components in $\z(\bar{\omega})$, and $m+n_1+\cdots+n_K$ rows, corresponding to the number of constraints in (\ref{eq:mapping1}), (\ref{eq:mapping3}), and (\ref{eq:mapping4}). Entries of $\FCBH^1$ are $0$, $1$, or values of components in $A^k$ $(1 \le k \le K)$. 

Any nonsingular submatrix of $\FCBA^1$ can be written as
\begin{equation}
\label{eq:AM}
\hat{\FCBA}^1=\left(\begin{array}{ccccc} \FCBH'_1 & & & & \\  & \ddots & & & \FCBE' \\ & & & \FCBH'_N & \end{array} \right),
\end{equation}
where $\FCBH'_i$ is a submatrix of $\FCBH^1$, $N$ is the number of such submatrices contained in $\hat{\FCBA}^1$, and $\FCBE'$ is a submatrix of $\FCBE^1$. As $M$ increases, the number of block matrices $\FCBH^1$ in $\FCBA^1$ increases, as does $N$ and thus the dimension of many submatrices $\hat{\FCBA}^1$, especially the maximal ones of $\FCBA^1$.       

As a thought experiment, suppose that for every aforementioned nonsingular matrix $\hat{\FCBA}^1$, we can find a matrix $\mathcal B$ to make the following linear transformation
\[
\mathcal B \hat{\FCBA}^1=\hat{\FCBA}_d.
\]
Here $\hat{\FCBA}_d$ is a nonsingular matrix in the form of 
\[
\left(\begin{array}{ccc} \hat{\FCBH}_1 & & \\  & \ddots &  \\ & &  \hat{\FCBH}_{N'}  \\ \end{array}
 \right),
\]
where $\hat{\FCBH}_i$ $(1 \le i \le N')$ are nonsingular submatrices with their dimensions bounded by an integer function of $m$ and $n_1$ and their entries depending only on $A^k$ ($1 \le k \le K$). Then
\[
(\hat{\FCBA}^1)^{-1}=\hat{\FCBA}_d^{-1} \mathcal B = \left(\begin{array}{ccc} \hat{\FCBH}_1^{-1} & & \\  & \ddots &  \\ & &  \hat{\FCBH}_{N'}^{-1}  \\\end{array} \right) \mathcal B.
\]
Suppose further that for all $\hat{\FCBA}^1$, a) entries of the corresponding matrix $\mathcal B$ are values in a finite set 
%that depends only on $A^k$ $(1 \le k \le K)$, 
and b) the number of non-zero entries on each row of $\mathcal B$ is bounded by an integer function of $m$ and $n_1$. Then the expression above implies that regardless how large $M$ is, the maximum absolute value of entries of $(\hat{\FCBA}^1)^{-1}$ stays bounded, and the number of non-zero entries on each row of $(\hat{\FCBA}^1)^{-1}$ remains finite.  Thus $||(\hat{\FCBA}^1)^{-1}||_\infty$ and thus the value of the aforementioned  $\kappa$ is bounded by some constant that is independent of $M$. 
 
The actual proof of Lemma \ref{lem:bondU} is based on a similar idea to the above even though we do not need to specify $\mathcal B$ and carry out the aforementioned transformation explicitly. Roughly speaking, we observe that in (\ref{eq:AM}), all $\FCBH'$s have finite dimensions and $\FCBE'$ has a finite number of columns ($\le n_1$) with $0$ or $1$ as their entries. This structure allows all non-zero entries of $(\hat{\FCBA}^1)^{-1}$, which recovers $||\mathbf u||_1$ from $||\mathbf u \hat{\FCBA}^1||_1$ (in the sense of (\ref{eq:recover})), to be obtained by inverting submatrices of $\hat{\FCBA}^1$. 
%Similar to $\hat{\FCBH}_i$ in $\hat{\FCBA}_d$,
These submatrices have finite dimensions and values of their entries do not depend on $M$. Consequently, in $(\hat{\FCBA}^1)^{-1}$, both the number of non-zero entries on each row and their values are finite, so $||(\hat{\FCBA^1})^{-1}||_\infty$ can remain bounded as $M$ increases.

Similar insights are used to prove Lemma \ref{lem:bondU} for cases with $k>1$, but the procedure gets considerably more complicated. For instance, when $k=2$, the LHS matrix of (\ref{eq:mapping1})-(\ref{eq:mapping4}) is  
\[
\FCBA^2=\left(\begin{array}{ccccc} \FCBH^2 & & & & \\  & \ddots & & & \FCBE^2 \\ & & & \FCBH^2 & \end{array} \right),
\]
where $\FCBH^2$ is the coefficient matrix of $\z(\bar{\omega})$ and $\y^1(\bar{\omega}')$, and $\FCBE^2$ is the coefficient matrix of $\y^2$.  Each $\FCBH^2$ corresponds to a particular $\bar{\omega}'$ and contains all $\bar{\omega}$ such that $\bar{\omega}' \sqsubset \bar{\omega}$ (see (\ref{eq:mapping1})-(\ref{eq:mapping2}) with $k=2$ and $k'=1$). Like $\FCBE^1$, $\FCBE^2$ has a fixed number $(n_2)$ of columns, and like $\FCBA^1$, entries of $\FCBA^2$ come from a finite set of values that do not depend on $M$. Different from $\FCBH^1$, the dimension of $\FCBH^2$ is not fixed but grows with $M$. Because of this difference, we need to introduce additional constructs to define the matrix structure recursively and prove the lemma by induction. 
\section{Numerical Results}
\label{sec:simulation}
We evaluate the performance of our policy by simulating its application to the examples of ATO systems shown in Figure \ref{fig:spesys}.  Each system features two distinct lead times. We vary cost parameters and lead times to generate various cases. 
For all cases the demand for products consists of independent Poisson processes,
with the arrival rate of demand for product $i$ denoted by $\lambda_i$.
In each case, $\underline{C}$ is the SP lower bound, $C_s $ is the average inventory cost determined by the simulation,
and the following optimality gap
\[
\Delta = \frac{C_s - \underline{C}}{\underline{C}}
\]
serves as the performance metric of our policy. 

For each case, we carry out 30 simulation runs. Depending on the lead times, the length of each run ranges from $150,000$ to $600,000$ time units, ensuring that it is at least 1250 times of the longer lead time. The first one-tenth of the simulation time is for warm-up. Values presented below are summaries of outputs from the remaining periods, averaged over the 30 simulation runs.

We first consider the $W$ system, using $27$ parameter sets given in \cite{Dogru2010} (Section 4.1).  The inventory holding cost of the common component, $h_0$, is normalized to unity and demand arrival rates are kept at $\lambda_1=\lambda_2=25$.  Values of $h_1$, $h_2$, $b_1$, and $b_2$, which are shown in the tables below, are chosen to cover a wide range of cost relationships.  Components $1$ and $2$ have the same lead time, which differs from that of component $0$. 

\begin{table}
	\footnotesize
	\begin{center}
			\caption{Optimality gaps: $W$ system, component 0 has the shorter lead time,$h_0=1$,$\lambda_1=\lambda_2=25$}
		\begin{tabular}{c|c|c|c|c|c|c|c|c|c|c|c}
			\hline
			case no.  & $h_1$  &$h_2$    & $b_1$ & $b_2$   &  \multicolumn{7}{c}{Optimality Gap ($L_1$: component $0$, $L_2$: components 1 and 2)}  \\  \cline{6-12}
			         &             &              &             &        & {$L_1 = 1$} &{$L_1 = 5$}&{$L_1 = 10$}&{$L_1 = 20$} &{$L_1=40$} & {$L_1=80$} & {$L_1=160$} \\
			&&&&& {$L_2 = 1.5$} &{$L_2 = 7.5$}&{$L_2 = 15$}&{$L_2 = 30$} & {$L_2=60$} & {$L_2=120$} & {$L_2=240$} \\
			\hline
			1 & 1 & 1 & 4 &4  & $0.03\%$ * & $0.14\%$ *  &  $0.04\%$ * & - & - & - &-  \\
			2 &	0.2& 0.2 & 2.4 &2.4 & $0.01\%$ * & $0.04\%$ * &  $0.08\%$ * & - &- & -& - \\
			3 & 1 & 5& 10 &6  & $0.10\%$ * & $0.07\%$ * & $0.19\%$ *  &- & - & - & - \\
			4 & 5& 5 & 12 &12 & $0.05\%$  & $0.01\%$ * & $0.00\%$ * & - & - & - & - \\
			5 & 0.2& 1 & 6 &4 & $ 0.56\%$ & $0.12\%$ * & $0.24\%$ * &- & - & - & - \\
			6 & 0.2& 0.2 & 2.4 &1.2 & $3.51 \%$ & $1.92\%$ & $1.34\%$ & $0.83\%$ &- & - & - \\
			7 & 1& 1 &4 &2 & $1.04\%$ &$0.68\%$ & $0.49\%$ &- & - &- &- \\
			8 & 5& 5 &12 &6 & $0.13\%$ & $0.08\%$ & $0.04\%$ * & - & - &- &-\\
			9 & 1& 0.2 &4 &2.4 & $2.20\%$ & $1.03\%$ & $0.79\%$ & -& - &- & - \\
			10& 0.2& 0.2 &6 &2 & $3.41\%$ & $1.63\%$ & $1.22\%$ & $0.87\%$ &- & - & -\\
			11& 1& 1 &10 &4 & $2.07\%$ & $1.24\%$ & $0.76\%$ &- &- &- &- \\
			12& 5& 5 &30 &12 & $0.40\%$ & $0.12\%$ * & $0.04\%$ * &- &- &- &- \\
			13& 0.2& 0.2 &6 &2.4 & $5.47\%$ & $2.89\%$ & $2.03\%$ & $1.56\%$ & $1.03\%$ &- & - \\
			14& 1& 0.2 &4 &1.2 & $5.85\%$ & $2.98\%$ & $2.07\%$ &1.45\% & $0.98\%$& - & - \\
			15& 0.2& 0.2 &6 &1.2 & $14.36\%$ & $7.17\%$ & $5.34\%$ & $3.75\%$ & $2.48\%$ & $1.92\%$ & $1.37\%$ \\
			16& 1& 1 &10 &2 & $5.20 \%$ & $2.55\%$ & $2.12\%$ & $1.45\%$ & $0.86\%$ & - & - \\
			17& 5& 1 &12 &4 & $1.45\%$ & $1.09\%$ & $0.72\%$ &- & - & - &-  \\
			18& 5& 5 &30 &6 & $0.84 \%$ & $0.32\%$ & $0.23\%$ & - & - & - & - \\
			19& 1& 0.2 &10 &2.4 & $7.10 \%$ & $3.73 \%$ & $2.79 \%$& $2.04\%$ & $1.52\%$ & $1.01\%$ & -  \\
			20& 5& 1 &12 &2 & $3.12\%$ & $1.63\%$ & $1.15\%$ & $0.78\%$ &- & - & - \\
			21& 1& 0.2 &10 &1.2 & $14.85\%$ & $7.98\%$ & $5.47\%$ & $4.08\%$ & $2.88\%$ & $2.17\%$ & $1.36\%$ \\
			22& 5& 0.2 &12 &2.4 & $4.65\%$ & $2.45\%$ & $1.83\%$ & $1.35\%$ & $0.87\%$ & - & -\\
			23& 5& 1 &30 &4 & $4.21 \%$ & $2.21\%$ & $1.47\%$& $1.30\%$  & $0.75\%$ & - & - \\
			24& 5& 0.2 &12 &1.2 & $9.03\%$ & $4.48 \%$ & $3.53\%$ & $2.51\%$ & $1.78\%$ & $1.26\%$ & -  \\
			25& 5& 1 &30 &2 & $7.79 \%$ & $3.81 \%$ & $3.08\%$ &$2.06\%$ & $1.59\%$ & $1.05\%$ & - \\
			26& 5& 0.2 &30 &2.4 & $9.65\%$ & $5.24\%$ & $3.91 \%$ & $2.76\%$ & $1.89\%$ & $1.38\%$ & - \\
			27& 5& 0.2 &30 &1.2 & $17.01\%$ & $8.65\%$ & $6.47\%$ & $4.53\%$ & $2.92\%$ & $2.30\%$ & $1.55\%$ \\
			\hline
		\end{tabular}
		\label{tbl:wcomshort}
	\end{center}
\end{table}

\begin{table}
	\footnotesize
	\begin{center}
		\caption{Optimality gaps: $W$ system, component 0 has the longer lead time, $h_0=1$,$\lambda_1=\lambda_2=25$}
		\begin{tabular}{c|c|c|c|c|c|c|c|c|c|c|c}
			\hline
			case no.  &   $h_1$  &$h_2$    & $b_1$ & $b_2$   & \multicolumn{7}{c}{ Optimality Gap ($L_1$: components 1 and 2, $L_2$: component 0)}   \\  \cline{6-12}
			&&&&& {$L_1 = 1$} &{$L_1 = 5$}&{$L_1 = 10$}&{$L_1 = 20$} &{$L_1=40$} & {$L_1=80$} & {$L_1=160$}  \\
			&&&&& {$L_2 = 1.5$} &{$L_2 = 7.5$}&{$L_2 = 15$}&{$L_2 = 30$} & {$L_2=60$} & {$L_2=120$} & {$L_2=240$} \\
			\hline
			1&1 &  1 & 4 &4 & $0.16\%$ & $0.08\%$ * & $0.08\%$ *  & - &- &- & -  \\
			2 &	 0.2& 0.2 & 2.4 &2.4 & $0.11\%$ * & $0.05\%$ * & $0.11\%$ * & -& - & - &- \\
			3&  1& 5 & 10 &6  & $0.30\%$ & $0.10\%$ * & $0.15\%$ *  &- &- &- &- \\
			4& 5& 5 & 12 &12 & $0.34\%$ & $0.11\%$ * & $0.14\%$ &- & - & - & - \\
			5&  0.2& 1 & 6 &4 &$0.78\%$ & $0.31\%$ & $0.20\%$ & - &- & - & - \\
			6& 0.2& 0.2 & 2.4 &1.2 & $3.46\%$ & $2.04\%$ & $1.28\%$ & $0.94\%$ & - & - & - \\
			7&  1& 1 &4 &2 & $1.65\%$ & $0.72\%$ & $0.43\%$ &- & - & - &- \\
			8&  5& 5 &12 &6 & $0.94\%$ & $0.42\%$ & $0.13\%$ &-& - & - & - \\
			9& 1& 0.2 &4 &2.4 & $2.72\%$ &$1.25\%$ &$0.95\%$ &- & - & - & - \\
			10&  0.2& 0.2 &6 &2 & $3.97\%$ & $1.88\%$ & $1.39\%$ & $0.92\%$ & - & - & -  \\
			11&  1& 1 &10 &4 & $2.63\%$ & $1.38\%$ & $0.97\%$ &-& - & - & - \\
			12& 5& 5 &30 &12 & $0.96\%$ & $0.37\%$ & $0.33\%$ &- & - & - & - \\
			13& 0.2& 0.2 &6 &2.4 & $5.76\%$ &  $2.77\%$ & $1.97\%$ &$1.52\%$ & $1.00\%$ & $0.53\%$ & -  \\
			14& 1& 0.2 &4 &1.2 & $7.74\%$ & $3.61\%$ & $2.74\%$ & $2.04\%$ & $1.33\%$ & $0.97\%$ &-  \\
			15& 0.2& 0.2 &6 &1.2 & $14.47\%$ & $7.29\%$ & $5.07\%$ & $3.67\%$ & $2.67\%$ & $2.04\%$ & $1.50\%$  \\
			16& 1& 1 &10 &2 & $6.56\%$ &$3.19\%$&$2.55\%$& $1.81\%$ & $1.41\%$ & $0.95\%$ & -   \\
			17& 5& 1 &12 &4 & $2.79\%$ & $1.55\%$ & $1.17\%$ & $0.86\%$ & - & - &- \\
			18& 5& 5 &30 &6 & $2.23\%$ &$1.11\%$ & $0.69\%$&-& - & - & -  \\
			19& 1& 0.2 &10 &2.4 & $8.78\%$ & $4.54\%$ & $3.33\%$ & $2.46\%$ & $1.60\%$ & $1.10\%$ & $0.77\%$ \\
			20& 5& 1 &12 &2 & $5.73\%$ &$3.19\%$ & $2.32\%$ & $1.58\%$  & $1.15\%$ & $0.83\%$ & - \\
			21& 1& 0.2 &10 &1.2 &$17.04\%$ &$8.88\%$ & $6.77\%$ &$4.87\%$ & $3.35\%$ & $2.79\%$ & $1.95\%$  \\
			22& 5& 0.2 &12 &2.4 & $7.43\%$ & $3.71\%$ & $2.72\%$ & $1.97\%$ & $1.33\%$ & $0.87\%$ & -\\
			23& 5& 1 &30 &4  & $6.32\%$ & $3.49\%$ & $2.44\%$ & $1.57\%$ & $1.19\%$ & $0.77\%$ & - \\
			24& 5& 0.2 &12 &1.2 & $13.02\%$ & $6.64\%$ & $4.99\%$ & $3.60\%$ & $2.62\%$ & $1.75\%$ & $1.31\%$  \\
			25& 5& 1 &30 &2 &  $10.83\%$ & $5.87\%$ & $4.17\%$ & $2.98\%$ & $2.21\%$ & $1.64\%$ & $0.86\%$  \\
			26& 5& 0.2 &30 &2.4& $13.29\%$ & $6.75\%$ & $5.02\%$ & $3.67\%$ & $2.59\%$ & $1.97\%$ & $1.03\%$ \\
			27& 5& 0.2 &30 &1.2  & $24.68\%$ & $12.13\%$ & $8.52\%$ & $6.00\%$ & $4.39\%$ & $3.09\%$ & $1.98\%$  \\
			\hline
		\end{tabular}
				\label{tbl:wcomlong}
	\end{center}
\end{table}

Table \ref{tbl:wcomshort} shows optimality gaps when component $0$ has the shorter lead time. 
Entries marked by * correspond to cases where the SP solution is within the $95\%$ confidence interval (CI) of the average cost estimated by $30$ simulation runs.
In these cases, the difference between the inventory cost under our policy and its lower bound is not statistically significant.
In all other cases, the optimality gap decreases as the lead times increase, consistent with the trend predicted by Theorem \ref{thm:final} above.
The optimality gaps fall to a level that is close to or below $1\%$ in many cases when $(L_1,L_2)=(20,30)$, and in a majority of cases when $(L_1,L_2)=(40,60)$. We stop running simulations for these cases with longer lead times (as indicated by ``-'' in the table). 
In a few cases, the optimality gap stays significantly above $1\%$ in all simulations, but nevertheless, still follows a clear pattern of converging to $0$. These cases (e.g., 15, 21, 27) are normally associated with a high $c_1/c_2$ ratio. Discussions in Section  4.2 in \cite{Dogru2010} explain why the gap tends to be larger under these parameter values, for systems with identical lead times. The same intuition applies here for systems with non-identical lead times. 

In Table \ref{tbl:wcomlong}, we use the same parameter set but let the common component have the longer lead time and the other two components have the same shorter lead time.  Like results in Table \ref{tbl:wcomshort}, it is evident that in each case, the optimality gap converges to zero as the lead time increases.  Nevertheless, these gaps are generally larger here. For instance, when $(L_1,L_2)=(160,240)$, the gap is close to $2\%$ in cases 21 and 27. We have additional runs for these two cases with $(L_1,L_2)=(320,480)$ and found the gap drops to $1.03\%$ and $1.26\%$ respectively.  

It is interesting to focus on the first four cases where $c_1=c_2$.  In Table \ref{tbl:wcomshort}, 
 the SP lower bound is within the $95\%$ CI except for case 4 when $(L_1,L_2)=(1,1.5)$.  A further calculation shows that this bound is within the (wider) $99.9\%$ CI.  In Table \ref{tbl:wcomlong}, the SP lower bound is outside the $95\%$ CI in cases 1, 3, and 4 when $(L_1,L_2)=(1,1.5)$.  Further calculations show that in cases 3 and 4, the SP lower bound is also outside the $99.9\%$ CI. Comparing the sample mean of the average cost obtained from the simulation with the lower bound, the $t$-value is $4.69$ for case 3 and $6.48$ for case 4. Given that the sample size is $30$, these values suggest that the inventory cost in both cases is significantly higher than the lower bound (with $\alpha<0.0005$). 

The observation is consistent with Theorem 4 in \cite{Reiman2012}, which concludes that in $W$ systems with $c_1=c_2$, the SP lower bound is reachable when the common component has the shorter lead time. However, the conclusion does not extend to $W$ systems in which the common component has the longer lead time.  Below  we use a special case of the $W$ system,  the $N$ system shown in Figure \ref{fig:spesys}, to explain this subtlety. 

When $c_1=c_2$,  serving a unit of either product $1$ or $2$ removes the same amount of inventory cost from the system. Thus the optimal allocation outcome prescribed by the SP solution is trivially attainable. However, when the common component has the longer lead time, it may not be possible to meet the SP-based inventory position target of the other component for all time.
To see this point, consider a time $t$ when a large batch of demand for product $1$ has just arrived, exhausting all inventory of component 0, both on-hand and in-transit at the moment. 

When component $0$ has the longer lead time, any new order of it will not arrive until $t+L_2$, which means that its net inventory level at $t+L_1$ is non-positive.  Correspondingly at  time $t$, the inventory position target of component $1$, constrained by the availability of component $0$ at $t+L_1$,  will be set at a non-positive level.  The system may not be able to meet this target because the actual inventory position is affected by previous ordering decisions, and therefore can be positive at $t$ even without ordering any new unit. Reducing it to a non-positive level requires removing existing units, which is not feasible. 

This situation does not happen when component $0$ has the shorter lead time. In this case, the SP sets a constant inventory position target for component $1$, which is always met (excluding the initial period) under a constant base stock policy. 
Any usage of component $1$ is accompanied by the usage of component $0$ by the same amount. 
Therefore, when the inventory position target of component $0$ needs to be reduced because of the lack of component $1$ within the next (shorter) lead time, its actual inventory position is also at a lowered level,  which allows the target to be met without removing any unit from the system. 

The impact of this difference on the optimality gap is shown by a few examples in Table \ref{tbl:N_system}. The first three columns give cost parameters. For each example, we compare the SP lower bound with the average cost from 100 simulation runs.  Results in columns 4, 5, 6 are from examples in which component $0$ has the longer lead time. In each example,  the SP lower bound is strictly below the lower end of the $99.9\%$ CI of the average cost. Comparing the sample mean of the latter cost with the lower bound, the $t$-values are $25.08$, $16.81$, and $33.30$ in these three cases, so the optimality gap is strictly positive at exceedingly high significance levels. Results in columns 7, 8, 9 are from examples in which component $0$ has the shorter lead time. In all examples, the SP lower bound is inside the $95\%$ CI, and thus certainly inside the (wider) $99.9\%$ CI. The $t$-values are $1.28$, $0.97$, and $1.34$ respectively, none of them is significant at $\alpha=0.05$. 

\begin{table}
	\footnotesize
	\begin{center}
	\caption{Optimality gaps: $N$-system with symmetric costs ($h_0=1$, $L_1=1$, $L_2=1.5$, $\lambda_1=\lambda_2=5$)}
		\begin{tabular}{c|c|c|c|c|c|c|c|c}
			\hline & & 
			 & \multicolumn{3}{c|}{$L_1$: component 1} & \multicolumn{3}{c}{$L_1$: component 0}\\ 
			$h_1$   & $b_1$ & $b_2$  &  \multicolumn{3}{c|}{$L_2$: component 0} & \multicolumn{3}{c}{$L_2$: component 1}\\ \cline{4-9}
			& & & lower bound & average & 99.9 \%CI & lower bound & average & $95\%$ CI \\ \hline
			0.1 & 0.4 & 0.5 & 21.38 & 21.56 &[21.54,21.59] & 18.95& 18.95 & [18.94,18.97] \\
			0.1 & 0.7 & 0. 8 & 28.82 & 29.00 & [28.96,29.04] & 25.26 & 25.27 & [25.25,25.29] \\
			1    &   1  & 2  & 51.38 & 51.98 & [51.91,52.04] & 49.24 & 49.25 & [49.23,49.28]\\
			\hline   
		\end{tabular}
		\label{tbl:N_system}
	\end{center}
\end{table}

We have also conducted simulations on the $M$ system shown in Figure \ref{fig:spesys}. We use the same cost parameters as these in \cite{Dogru2017}. The inventory holding costs of both components, $h_1$ and $h_2$, are set to unity. Backlog costs, $b_0$, $b_1$, and $b_2$, are varied (see the second row in Tables \ref{tbl:M1short} and \ref{tbl:M1long}) to generate four regions of different cost relationships:  $c_1 + c_2 < c_0$ (region A); (2) $c_2 \le c_1 < c_0 \le c_1 + c_2$ (region B); (3) $c_2 < c_0 \le c_1$ (region C); and (4) $c_0 \le c_2 \le c_1$ (region D). As discussed in \cite{Dogru2017}, our allocation policy specializes to allocation rules that are qualitatively different between the regions. 
Tables \ref{tbl:M1short} ande \ref{tbl:M1long} show optimality gaps for cases when component $1$ has the shorter and longer lead times respectively.  In both cases, the gaps are above zero when $(L_1,L_2)=(1,1.5)$, and significantly so in regions A and D.  Nevertheless, in each case, there is also a clear trend that the gap converges towards zero as the lead times increase, as is predicted by Theorem \ref{thm:final}.

\setlength{\tabcolsep}{20pt}
\begin{table}
\footnotesize
	\caption[caption]{Optimality gaps: $M$ system, component 1 has the shorter lead time, $h_1=h_2=1$,$\lambda_0=25$,$\lambda_1=\lambda_2=50$}
\begin{center}
\begin{tabular}{l l|c|c|c|c}
\hline
\multicolumn{2}{c|}{lead time} & A & B & C & D
\\ \cline{3-6}
&  &\multicolumn{1}{l|}{$b_0=8$}  &\multicolumn{1}{l|}{$b_0=3$}  & \multicolumn{1}{l|}{$b_0=2$}  &\multicolumn{1}{l}{$b_0=1$} \\ 
\multicolumn{2}{c|}{$L_1$:component 1}& \multicolumn{1}{l|}{$b_1=3.5$}  & \multicolumn{1}{l|}{$b_1=2.5$} & \multicolumn{1}{l|}{$b_1=3.5$} & \multicolumn{1}{l}{$b_1=8$} \\
\multicolumn{2}{c|}{$L_2$:component 2} &  \multicolumn{1}{l|}{$b_2=1$}  & \multicolumn{1}{l|}{$b_2=1$} & \multicolumn{1}{l|}{$b_2=1$} & \multicolumn{1}{l}{$b_2=3$} \\ \hline
$L_1=1$ & $L_2=1.5$ & $10.03\%$ & $4.60\%$ & $5.72\%$ & $23.20\%$
\\ 
$L_1=5$ & $L_2=7.5$ & $4.82\%$ & $2.30\%$ & $2.87 \%$ & $11.27\%$
\\ 
$L_1=10$ &$L_2=15$ & $3.29\%$ & $1.72\%$ & $2.07\%$ & $8.18\%$ 
\\
$L_1=20$ & $L_2=30$ & $2.32\%$ & 1.17\% & $1.44\%$ & $6.13\%$ 
\\
$L_1=40$ & $L_2=60$ & $1.75\%$ & - & 1.00\% & $4.45\%$
\\ 
$L_1=80$ & $L_2=120$ & $1.37\%$ & - & - & $3.20\%$ 
\\
\hline
\end{tabular}
\end{center}
\label{tbl:M1short}
\end{table}

\setlength{\tabcolsep}{20pt}

\begin{table}
\footnotesize
	\caption[caption]{Optimality gaps: $M$ system, component 2 has the shorter lead time, $h_1=h_2=1$, $\lambda_0=25$, $\lambda_1=\lambda_2=50$}
\begin{center}
\begin{tabular}{l l|c|c|c|c}
\hline
\multicolumn{2}{c|}{lead time} & A & B & C & D
\\ \cline{3-6}
&  &\multicolumn{1}{l|}{$b_0=8$}  &\multicolumn{1}{l|}{$b_0=3$}  & \multicolumn{1}{l|}{$b_0=2$}  &\multicolumn{1}{l}{$b_0=1$} \\ 
\multicolumn{2}{c|}{$L_1$:component 2}& \multicolumn{1}{l|}{$b_1=3.5$}  & \multicolumn{1}{l|}{$b_1=2.5$} & \multicolumn{1}{l|}{$b_1=3.5$} & \multicolumn{1}{l}{$b_1=8$} \\
\multicolumn{2}{c|}{$L_2$:component 1} &  \multicolumn{1}{l|}{$b_2=1$}  & \multicolumn{1}{l|}{$b_2=1$} & \multicolumn{1}{l|}{$b_2=1$} & \multicolumn{1}{l}{$b_2=3$} \\ \hline
$L_1=1$ & $L_2=1.5$ & $9.07\%$ & $4.20\%$ & $5.20\%$ & $22.48\%$
\\
$L_1=5$ & $L_2=7.5$ & $4.42\%$ & $2.25\%$ & $2.67 \%$ & $11.16\%$
\\ 
$L_1=10$& $L_2=15$ & $3.17\%$ & $1.40\%$ & $1.87\%$ & $8.12\%$ 
\\ 
$L_1=20$ & $L_2=30$ & $2.25\%$ & 0.99\% & $1.29\%$ & $5.74\%$ 
\\ 
$L_1=40$ & $L_2=60$ & $1.46\%$ & - & $0.88\%$ & $3.99\%$
\\ 
$L_1=80$ & $L_2=120$ & $0.94\%$ & - & - & $2.69\%$ 
\\ \hline
\end{tabular}
\end{center}
\label{tbl:M1long}
\end{table}

\iffalse
\begin{figure}
	\begin{center}
		\includegraphics[width=5in,height=3in]{M_system_result.png}
		\caption{Asymptotic Optimality in four regions of $M$ system when ${ {{L_2} \over {L_1}} = 1.5}$.}\label{M_system}
	\end{center}
\end{figure}
\fi
\section{Conclusion}
\label{sec:conclusion}

Optimal inventory control of ATO systems with multiple products is a long-standing problem in the literature. In principle, this paper has settled the issue of developing an asymptotically optimal control policy for systems with general BOMs and deterministic lead times. We have ``collapsed'' the design of a dynamic control policy for minimizing the average cost over an infinite time horizon to the solution of certain multi-stage stochastic programs. The approach leads to drastically simplified analyses and feasible inventory policies with provable optimality properties. 

This paper is a continuation of a stream of past work on ATO inventory systems. The multi-stage SP in Section \ref{sec:SPdisc} generalizes the two-stage SP in \cite{Reiman2015}, which applies only to systems with identical lead times. In comparison with the formulation in \cite{Reiman2012}, the SP here is a better alternative as it directly sets the same lower bound on the average cost, eliminating the need to take the infimum of the objective values, and also yields optimal solutions that have finite values and hence can be used as parameters of inventory control policies. While the allocation policy follows the same principle as that in \cite{Reiman2015}, we break new ground in developing a replenishment policy. The conventional base stock policy is justified by asymptotic analysis for  use in systems with identical or near identical lead times (\cite{Reiman2015}, \cite{Reiman2016}) while a different policy has been shown to be optimal for one-product systems with non-identical lead times (\cite{Rosling1989}). Subsuming both as special cases, our policy is supported by the proof of asymptotic optimality and numerical results for its use in systems with general BOMs and lead times.

As a result of the new policy development, the asymptotic analysis is more involved than for systems with identical lead times in \cite{Reiman2015}. Instead of setting constant targets for inventory positions, the new policy updates them repeatedly over time based on changes of relevant system states,  giving rise to the question about the stability of these targets.  In general, it is not possible to keep actual inventory positions at their targets for all times, and in theory, the influence of the discrepancies can last into the indefinite future. We address these issues with new approaches, such as the formulation and analysis of the ``stochastic tracking model'' and the characterization and exploitation of the matrix structures of SP constraints.  We develop these technical machineries for general settings, making it possible to apply them to models other than ATO systems, e.g.,  to prove  convergence of other stochastic processes or address Lipschitz continuity of certain types of infinite LPs.

This paper is, in a sense, the culmination of the stream of work mentioned above. On the other hand, this  is not the end of the story.
In particular, making our policy implementable in practical ATO systems requires new approaches to overcome the computational complexity of multi-stage SPs. In this paper, we consider the use of finite-dimensional LPs to approximate the original problems to generate stable targets for inventory positions and backlog levels. For the approximation to be sufficiently accurate, the dimensions of these LPs may need to be exceedingly high, especially when the system has many different and very long lead times. Moreover, the SPs need to be re-optimized repeatedly over time to update policy parameters, which makes the computation even more expensive.

Recent studies on two-stage SPs show that exploring structural properties of ATO problems implied by their BOMs can facilitate solution procedures  (\cite{DeValve2018}, \cite{Dogru2017},  \cite{Zipkin2016}). How to extend this strategy to multistage SPs is an issue yet to be explored. It has also been shown that when the percentage difference between the longest and shortest lead times becomes insignificant, it can be asymptotically optimal to follow a base stock policy \cite{Reiman2016}, reducing the SP to a two-stage problem and eliminating the need to update inventory position targets. Similarly, it also seems possible to preserve asymptotic optimality while reducing computational efforts by applying simple replenishment rules without resorting to SP solutions on components whose lead times are insignificant fractions of the longest lead time. In general, there can be opportunities to systematically explore exploiting small lead time differences to reduce the number of stages of the SPs. Furthermore, one may also find ways to reduce the computational burden by re-solving multi-stage SPs less frequently.  While interesting and important, these topics are certainly beyond the scope of one paper, and we leave them for future research.

As a broader message, our work highlights the power of asymptotic analysis to tackle difficult inventory models that defy exact analysis. There is no shortage of such problems in the inventory theory literature, giving rise to ample opportunities for innovative research. On this subject, we refer to a recent survey \cite{Goldberg2019} for detailed discussions.

\section*{Appendix A: Proof of Theorems}
\label{sec:appendix}

\medskip

\subsection*{Proof of Theorem \ref{thm:solbound}:}
The proof will use the  following primal-dual transformation of the last stage LP in (\ref{eq:MSP2}):
\begin{align}
\begin{split}
\label{eq:defvp0}
\varphi^0(\mathbf y^1,\cdots \mathbf y^K,\mathbf x)
&=-\max_{\mathbf z} \{\mathbf c \cdot \z|\z \le \x, A \z \le \y \}
\\
&=-\min_{\bnu, \bk \ge 0} \left \{\x \cdot \bnu+\y \cdot \bk | \bnu+A' \bk =\mathbf c \right\}
\\
&= -\mathbf c \cdot \x -\min_{\bk \ge 0} \left \{(\y-A \x) \cdot \bk | A' \bk \le \mathbf c \right\}
\\
&= \max_{\bk \ge 0} \left \{(A \x - \y) \cdot \bk | A' \bk \le \mathbf c \right\}- \mathbf c \cdot \x.
\end{split}
\end{align}
Given that all components of $A$ are nonnegative, it is always optimal to set $\kappa_j=0$ for any $j$ such that $\mathbf A_j \cdot \x-y_j \le 0$ $(1 \le j \le n)$. Thus without the loss of generality, 
we replace $(A \x - \y)$ with $(A \x - \y)^+$ in (\ref{eq:defvp0}), in which case
\begin{align}
\begin{split}
\label{eq:dualv0}
\varphi^0(\mathbf y^1,\cdots \mathbf y^K,\mathbf x) &= \varphi^0_d(\mathbf y^1,\cdots \mathbf y^K, \x)-\mathbf c \cdot \x
\\
~~~\mbox{where}~~
\varphi^0_d(\mathbf y^1,\cdots \mathbf y^K, \x) &:=\max_{\bk \ge 0} \left \{(A \x -\y)^+\cdot \bk | A' \bk \le \mathbf c \right\}
\end{split}
\end{align} 

In (\ref{eq:dualv0}), for any feasible solution of $\varphi^0_d(\mathbf y^1,\cdots \mathbf y^K, \x)$, 
\[
\kappa_j \le \bar{\kappa} ~(1 \le j \le n), ~~\mbox{where}~~\bar{\kappa}=\frac{\bar{c}}{\underline{a}},
\]
where $\bar{c}$ and $\underline{a}$ are defined in Table \ref{tbl:minmaxpara}. 
Let $\y^{l*}$ $(1 \le l \le k)$ be an optimal solution on the sample path $(\BD^k,\cdots,\BD^1$).
Then
\begin{equation}
\label{eq:expobj}
\varphi^k(\y^{k+1},\cdots,\y^K,\x)= \sum_{l=1}^{k} \mathbf h^l \cdot \E[\mathbf y^{l*}]+\E[\varphi^0_d(\y^{1*},\cdots, \y^{k*}, \y^{k+1},\cdots,\y^K, \x+\underline{\BD}^k)] -\mathbf c \cdot \E[\x+\underline{\BD}^k].
\end{equation}
  
Let $\varphi^{k',j} (\y^{k+1},\cdots \y^K,\x)$ be a deviation from the optimal objective value $\varphi^k(\y^{k+1},\cdots,\y^K,\x)$, obtained by changing $y_j^*$ to $y_j^*-1$ while keeping all other values at the optimal levels. 
Then
\[
\varphi^{k',j}(\y^{k+1},\cdots \y^K,\x)-\varphi^k(\y^{k+1},\cdots,\y^K,\x) \ge 0,
\]
i.e., 
\begin{eqnarray*}
\nonumber
&~&-h_j+\E[\kappa_j (\mathbf A_j \cdot (\x+\underline{\BD}^k)-(y^*_j-1))^+]-\E[\kappa_j (\mathbf A_j \cdot (\x+\underline{\BD}^k)-y^*_j)^+]
\\
\nonumber
&=&-h_j+\E\left [\kappa_j \mathbf 1 \{\mathbf A_j \cdot (\x+\underline{\BD}^k) \ge y^*_j\}\right]
\\
&\ge& 0.
\end{eqnarray*}
For the above condition to hold, it is necessary that
\[
\bar{\kappa} \Pr\{\bar{a}(||\x||_1+||\underline{\BD}^k||_1) \ge  y^*_j \} \ge \underline{h}.
\]
When $y^*_j>0$, apply Markov's Inequality to the above,
\[
\frac{\E\left[\bar{a}(||\x||_1+||\underline{\BD}^k||_1)\right]}{y^*_j} \ge \underline{h}/\bar{\kappa},
\]
which yields an upper bound 
\[
y^*_j \le \bar{\beta} \left(||\x||_1+\E\left[||\underline{\BD}^k||_1\right]\right),
\]
where $\bar{\beta}>0$ is determined by $\bar{\kappa}$, $\bar{a}$, and $\underline{h}$, and hence ultimately depends only on $A$, $\mathbf b$, and $\mathbf h$.
This upper bound obviously applies when $y^*_j \le 0$.

To prove the lower bound, fix a component $j$ of $\mathbf y^k$ (i.e., $k_j=k$) and consider the following feasible solution to 
(\ref{eq:dualv0})
\begin{equation}
\kappa_{j'} =  \left \{
\begin{array}{lr}
 0 & k_j>k,
\\
 h_j+\underline{b} /\bar{a}& j'=j,
\\
 h_{j'} &~~~\mbox{otherwise},
\end{array}
\right .
\end{equation}
which yields a weakly lower objective value than the optimal one, i.e., 
\[
 \sum_{l=1}^{k}  \mathbf h^l \cdot (A^l \x-\y^{l*}) + \underline{b}/\bar{a} (\mathbf A_j \cdot \x-y^*_j) \le \varphi^0_d(\mathbf y^{1*},\cdots,\y^{k*},\y^{k+1},\cdots,\y^K, \x).
\]
Apply the inequality to (\ref{eq:expobj}) (notice that $\x$ in $\varphi_d^0()$ above corresponds to $\x+\underline{\BD}^k$ in $\varphi_d^0()$ in (\ref{eq:expobj})),
\begin{align}
\begin{split}
\label{eq:ineq1}
~\varphi^k(\mathbf y^{k+1},\cdots \mathbf y^K, \x) 
\ge & \mathbf h^k \cdot \y^{k*} +\mathbf h^k \cdot (A^k(\x+\E[\underline{\BD}^k])-\y^{k*}) 
+\underline{b}/\bar{a} [\mathbf A_j \cdot (\x+\E[\underline{\BD}^k])-y^*_j]
\\
&~+\sum_{l=1}^{k-1} \mathbf h^l \cdot \E[\y^{l*}]+ \sum_{l=1}^{k-1} \mathbf h^l \cdot  (A^l(\x+\E[\underline{\BD}^k]) - \E[\y^{l*}])
 -\mathbf c \cdot (\x+\E[\underline{\BD}^k])
\\
\ge & - (\underline{b} /\bar{a})y^*_j-\mathbf c' \cdot (\x+\E[\underline{\BD}^k]),
\end{split}
\end{align}
where $\mathbf c'=\mathbf c - \sum_{l=1}^k (A^k)' \mathbf h^k-(\underline{b}/\bar{a})\mathbf A_j$.

Let $\bk^0$ be an optimal solution of $\varphi^0_d(\mathbf 0,\cdots,\mathbf 0, \y^{k+1}, \cdots \y^K,\x)$. 
Then $\{\y^l=0~(1 \le l \le k),\bk^0\}$ is a feasible solution of $\varphi^k(\mathbf y^{k+1},\cdots \mathbf y^K, \x)$ and yields the following objective value:
\[
\sum_{j':k_{j'} \le k} \E \left[\kappa_{j'}^0 \mathbf A_{j'} \cdot \left(\x+\underline{\BD}^k\right)\right]+\sum_{j':k_{j'}>k}  \E \left [\kappa_{j'}^0 \left (\mathbf A_{j'} \cdot \left(\x+\underline{\BD}^k\right)-y_{j'}\right) \right]- \mathbf c \cdot (\x+\E[\underline{\BD}^k]). 
\]
Therefore
\begin{equation}
\varphi^k(\mathbf y^{k+1},\cdots \mathbf y^K, \x)  \le\sum_{j'=1}^{n} \mathbf  A_{j'} \cdot\E[\kappa_{j'}^0 (\x+\underline{\BD}^k)]+\sum_{j':k_{j'}>k} |y_{j'}|\E[\kappa_{j'}^0]-\mathbf c \cdot (\x+\E[\underline{\BD}^k])
\label{eq:ineq2}
\end{equation}
Using (\ref{eq:ineq1}) and (\ref{eq:ineq2}), and applying that $\kappa_j \le \bar{\kappa}$:
\begin{equation*}
%\label{eq:lbdy}
(\underline{b}/\bar{a}) y^*_j 
\ge
 -\bar{\kappa} \left(\sum_{j'=1}^n \mathbf A_{j'}\cdot (\x+\E[\underline{\BD}^k])+\sum_{j':k_{j'}>k} |y_{j'}|\right)+(\mathbf c-\mathbf c')\cdot(\x+\E[\underline{\BD}^k]).
\end{equation*}
%Since\[\left(\sum_{j'=1}^{n_k} h_{j'} \mathbf A_{j'}+(\underline{b}/\bar{a})\mathbf A_j\right)\cdot(\x+\E[\underline{\BD}^k]) \ge - (\bar{a} n\bar{h}+m \underline{b})\left(||\x||_1+\E[||\underline{\BD}^k||_1]\right),\]
Therefore, there exists a constant $\underline{\beta}>0$, depending only on $A$, $\mathbf b$, and $\mathbf h$, 
such that:
\[
y^*_j \ge - \underline{\beta} \left(||\x||_1+\E[||\underline{\BD}^k||_1]+\sum_{l=k+1}^{K}||\y^l||_1\right), ~~\jrange.~~~\blacksquare
\]

\subsection*{Proof of Theorem \ref{thm:sameobj}:}

By substituting $\y^k$ with $\y^k-A^k \boldsymbol \alpha$ ($1 \le k \le K$) and $\z$ with $\z-\boldsymbol \alpha$ in (\ref{eq:MSP1}):
\[
\Phi^K(\bal)=\Psi^K_{\bal}-\mathbf b \cdot \bal
\]
where
\begin{align}
 \label{eq:uMSP1}
 \begin{split}
 \Psi^K_{\bal} &= \inf_{\mathbf{y}^K \ge -A^K \boldsymbol \alpha} \{\mathbf{h}^K \cdot \mathbf{y}^K + \E [\Psi^{K-1}_{\bal}(\y^K, \BD^K) ]\},
 \\
 \Psi^k_{\bal}(\y^{k+1},\cdots, \y^K,\x) &= \inf_{\mathbf{y}^k \ge -A^k \boldsymbol \alpha} \{\mathbf{h}^k \cdot \mathbf{y}^k + \E [  \Psi^{k-1}_{\bal}(\y^k,\cdots, \y^K,\x + \BD^k)]\},~~k=K-1,\cdots,1,
 \\
 \Psi^0_{\bal}(\y^1,\cdots, \y^K,\x) &= - \max_{\mathbf{z} \ge -\boldsymbol \alpha} \{ \mathbf{c}\cdot \mathbf{z}  |  \mathbf{z} \le \mathbf{x}, A^k \mathbf{z} \le \y^k,~1 \le k \le K\}.
 \end{split}
 \end{align}
Using (\ref{eq:lbound}), below we prove (\ref{eq:sameobj}) by showing that
\begin{equation}
\label{eq:balconv}
\inf_{\bal \ge \mathbf 0} \{\Psi^K_{\bal}\}= \varphi^K.
\end{equation}

Let $M$ be a positive integer and $\mathbf M$ be a $m$-dimensional vector with all components equal to $M$.
Let $\varphi_M^K$ be the optimal objective value of the SP defined in (\ref{eq:MSP2}) with inputs $\BD^k$ replaced by 
\[
\BD^k_M:=\BD^k \wedge \mathbf M, \quad 1 \le k \le K. 
\]
Likewise, let $\Psi^K_{M,\alpha}$ be the optimal objective value of the SP defined in (\ref{eq:uMSP1}) with the replacement of $\BD^k$ by $\BD^k_M$ $(1 \le k \le K)$.
Following Theorem \ref{thm:unicom}, as $M \rightarrow \infty$, $\varphi_M^K$ converges to $\varphi^K$ and $\Psi^K_{M,\bal}$ converges to $\Psi^K_{\bal}$ \emph{uniformly} over all $\bal$ (see (\ref{eq:Psiconv}) and the discussion below in the proof of Theorem \ref{thm:unicom}. While the latter theorem and its proof appear later in the paper, they do not rely on the conclusion that we are proving here). Thus to prove (\ref{eq:balconv}), we only need to show that for any given $M$, there exists some $\bal$, such that
\begin{equation}
\label{eq:sametruncateSP}
\Psi^K_{M,\bal}=\varphi^K_M.
\end{equation}

Let $\y^{k*}$ $(1 \le k \le K)$ be an optimal solution of $\varphi_M^K$. 
By Theorem \ref{thm:solbound},
 $\y^{K*}_M$ is finite.
For $k=K-1,\cdots,1$, 
since every entry of $\x$ in $\varphi_M^k(\y^{k+1*}_M,\cdots,\y^{K*}_M,\x)$, which is 
\[
\BD^K_M+\cdots \BD^{k+1}_M
\]
is bounded by $M(K-k)$,
by a simple induction on (\ref{eq:solbound}), $\y^{k*}_M$ has a finite upper bound that applies to all sample paths. 
This means that when $\bal$ is sufficiently large, 
$\y_M^{k*}$ ($1 \le k \le K$) are feasible solutions to $\Psi^k_{M,\bal} (\y^{k+1*}_M,\cdots,\y^{K*}_M,\x)$ ($1 \le k \le K$).
On the other hand, any optimal values of $\y^k$ for $\Psi^k_{M,\bal}(\y^{k+1},\cdots,\y^{K},\x)$ are obviously feasible for  $\varphi^k_M(\y^{k+1},\cdots,\y^{K},\x)$ $(1 \le k \le K)$. 
Thus to prove (\ref{eq:sametruncateSP}), 
we only need to show that if $\bal$ is sufficiently large, 
\begin{equation}
\label{eq:lpeq}
\Psi^0_{M,\bal} (\y^1,\cdots, \y^K,\x)= \varphi^0_M(\y^1,\cdots,\y^K,\x)
\end{equation}
for any  $\x  \le K \mathbf M$
and $\y^k$ $(1 \le k \le K)$ bounded by a constant that depends on $M$ but not $\bal$. By definition, except for restrictions on inputs,  $\varphi^0_M(\y^1,\cdots,\y^K,\x)$ and $\Psi^0_{M,\bal} (\y^1,\cdots, \y^K,\x)$ are the same LPs as $\varphi^0(\y^1,\cdots,\y^K,\x)$ and $\Psi^0_{\bal} (\y^1,\cdots, \y^K,\x)$ respectively.

Obviously, for any $\y^k$ $(1 \le k \le K)$ and $\x$
\begin{equation}
\label{eq:inverseinq}
\Psi^0_{M,\bal} (\y^1,\cdots, \y^K,\x) \ge  \varphi^0_M(\y^1,\cdots,\y^K,\x),
\end{equation}
since any feasible solution for the LP on the LHS is also feasible for the LP on the RHS. To prove the inequality holds in reverse, 
observe that if $\mathbf z^*$ optimizes
\begin{eqnarray*}
\varphi^0_M (\y^1,\cdots, \y^K,\x) &=& -\max_{\z} \left \{\mathbf c \cdot \z|
A^k \z \le \y^k~(1 \le k \le K),~ \z \le \x \right\}
\\
&=&  -\max_{\z} \left \{\sum_{i=1}^{m} c_i z_i|
\sum_{i=1}^m a_{ji} z_i \le y_j~(1  \le j \le n),~z_i \le x_i~(1 \le i \le m) \right \}.
\end{eqnarray*}
Then
\[
z^*_i \ge - \underline{\alpha}~~~\mbox{where}~~\underline{\alpha}:=\max_{1 \le j \le n} \left \{\frac{|y_j|+\sum_{i=1}^m a_{ji} x_i}{\min_{a_{ji}>0} \{a_{ji}\}}\right\},~1 \le i \le m.
\]
Otherwise, one can improve the objective value by increasing $z_i$ without violating any constraint. 
Since $\y^k$ ($1 \le k \le K$) and $\x$ are bounded,
$\underline{\alpha}$ is also bounded by some constant $\underline{\alpha}_{\max}$.
When $\alpha \ge \underline{\alpha}_{\max}$, any optimal solution to $\varphi^0_M (\y^1,\cdots, \y^K,\x)$  is a feasible solution to $\Psi^0_{M,\bal} (\y^1,\cdots, \y^K,\x)$ 
so (\ref{eq:inverseinq}) holds in reverse. This completes the proof of the theorem. 
$~~~~\blacksquare$

\subsection*{Proof of Lemma \ref{lem:NBD}:}

Observe that
\begin{eqnarray}
\nonumber
\E\left[\sup_{t \ge 0}\left(|\hD_i(0,t)|-\sqrt{L}t\nu\right)^+\right]
&\le&
\sum_{\tau=1}^{\infty} \E\left[\sup_{\tau  \le t <\tau+1}\left(|\hD_i(0,t)|-\sqrt{L}t\nu\right)^+\right]
\\
\nonumber
&~&+\E\left[\sup_{L^{-1/4} \le t <1}\left(|\hD_i(0,t)|-\sqrt{L}t\nu\right)^+\right]
\\
\label{eq:spt}
&~&+\E\left[\sup_{0 \le t <L^{-1/4}}\left(|\hD_i(0,t)|-\sqrt{L}t\nu\right)^+\right].
\end{eqnarray}

Notice that $\tsD(t)$ $(t \ge -1)$ is stationary process and $|\hD_i(0,t)|$ ($t \ge -1$) is a sub-martingale.
 \begin{eqnarray}
\nonumber
\sum_{\tau=1}^{\infty}\E\left[\left(\sup_{\tau \le t \le \tau+1}|\hD_i(0,t)|-\sqrt{L} t \nu\right)^+\right]
&\le&
\sum_{\tau=1}^{\infty}\E\left[\left(\sup_{\tau \le t \le \tau+1}|\hD_i(0,t)|-\sqrt{L} \tau \nu\right)^+\right]
\\
\nonumber
&=&\sum_{\tau=1}^{\infty}\E\left[\left(\sup_{0 \le t \le 1}|\hD_i(0,t+\tau)|-\sqrt{L} \tau \nu\right)^+\right]
\\
\nonumber
 &=& \sum_{\tau=1}^{\infty}\int_{\sqrt{L} \tau \nu}^\infty \Pr\left\{\sup_{0 \le t \le 1} |\hD_i(0,t+\tau)| \ge x \right \} dx
 \\
 \nonumber
 &\le& \sum_{\tau=1}^{\infty} \int_{\sqrt{L} \tau \nu}^\infty \frac{\E[|\hD_i(0,\tau+1)|^p]}{x^p} dx
\\
\label{eq:doobequ}
&=& \frac{(\sqrt{L} \nu)^{1-p}}{p-1} \sum_{\tau=1}^{\infty} \frac{\E[|\hD_i(0,\tau+1)|^p]}{\tau^{p-1}},
 \end{eqnarray}
where the second inequality is Doob's inequality and $p>0$. 

To bound the sum in the above, observe that
\[
\E[(\hD_i(0,\tau+1))^6]=\E\left[\left(\sum_{s=1}^{\tau+1}\hD_i(s-1,s)\right)^6\right],
\]
where $\E[\hD_i(s-1,s)]=0$.
Since $\hD_i(s-1,s)$ ($1 \le s \le \tau+1$) is an i.i.d sequence,
\[
\E[(\hD_i(s-1,s))^k]=\E[(\hD_i(0,1))^k],~~1 \le k \le 6,~s=1,\cdots,\tau+1.
\]
By eliminating terms that contain $\E[\hD_i(s-1,s)]=0$
and thus equal zero
in the expansion,
\begin{eqnarray}
\nonumber
\E[(\hD_i(0,\tau+1)^6] &=& C^{\tau+1}_1 \times \E[(\hD_i(0,1))^6]
\\
\nonumber
&~&+ 
C_2^{\tau+1} \left(2C_2^6\times\E[(\hD_i(0,1))^4]\times\E[(\hD_i(0,1))^2]+C_3^6\times(\E[(\hD_i(0,1))^3])^2\right)
\\
\label{eq:expansion}
&~&
+C^{\tau+1}_3 \times C_2^6 \times C_2^4 \times (\E[(\hD_i(0,1))^2])^3,
\end{eqnarray}
where $C_r^q$ denotes $q$ choose $r$ ($q$, $r$ integers and $q \ge r$).
Let 
\[
\upsilon=\sum_{\tau=1}^{\infty}\frac{\E[(\hD_i(0,\tau+1))^6]}{\tau^5},
\]
which is a finite constant: since the jump size of the Compound Poisson process is assumed to have a finite moment of order 6,   $\E[(\hD_i(0,1))^k]$ ($k=2,3,4,6$) are all finite. Moreover $C^{\tau+1}_k$ ($k=1,2,3$) are on the order of $\tau^3$ or smaller.
By (\ref{eq:doobequ}) with $p=6$,
\begin{equation}
\label{eq:expbound}
\sum_{\tau=1}^{\infty}\E\left[\left(\sup_{\tau \le t \le \tau+1}|\hD_i(0,t)|-\sqrt{L} t \nu\right)^+\right] \le \frac{\upsilon}{5\nu^5} L^{-5/2}.
\end{equation}
Applying the above and bounds in (\ref{eq:boundonED}) and (\ref{eq:boundonEDD}) of Lemma \ref{lem:RW3} to (\ref{eq:spt}) proves (\ref{eq:dlim}).\label{key}
$\blacksquare$

\subsection*{Proof of Lemma \ref{lem:RW4}:}

For each $i=1,\cdots,m$,
\begin{eqnarray*}
\E\left[\sum_{\tau=0}^{\infty}\sup_{0 \le t<1}\left(\hejD_i(t)-\sqrt{L} \nu \tau\right)^+\right]
&=& 
\E\left[\sup_{0 \le t <1} \hat{d}^{(L)}_i(t)\right]+
\sum_{\tau=1}^{\infty}\E\left[\sup_{0 \le t <1}\left(\hejD_i(t)-\sqrt{L} \nu \tau\right)^+\right]
\\
&\le&
3 \lambda^{\frac{1}{2+\delta}}(1+\eta_i)L^{-\frac{\delta}{2(2+\delta)}}
+\sum_{\tau=1}^{\infty}\E\left[\sup_{0 \le t <1}\left(\hejD_i(t)-\sqrt{L} \nu \tau\right)^+\right],
\end{eqnarray*}
where the first equality follows from Tonelli's Theorem, and the next inequality is obtained by applying Lemma 2 in \cite{Reiman2015} to bound the first expected value ($\tau=0$).
Observe that
\begin{eqnarray*}
\sum_{\tau=1}^{\infty}\E\left[\sup_{0 \le t <1}\left(\hejD_i(t)-\sqrt{L} \nu \tau\right)^+\right]
&=&\sum_{\tau=1}^{\infty}\int_{\sqrt{L}\nu \tau}^{\infty}\Pr\left\{\sup_{0 \le t <1} \hejD_i(t) \ge x \right\}dx
\\
&\le&\sum_{\tau=1}^{\infty}\int_{\sqrt{L}\nu \tau}^{\infty}\frac{\E\left[\left(\sup_{0 \le t<1} \hejD_i(t) \right)^{2+\delta}\right]}{x^{2+\delta}}dx
\\
&=&\sum_{\tau=1}^{\infty}\frac{\E\left[\left(\sup_{0 \le t <1} \hejD_i(t) \right)^{2+\delta}\right]}{(1+\delta)(\sqrt{L}\nu)^{1+\delta}}\frac{1}{\tau^{1+\delta}}
\\
&\le&\sum_{\tau=1}^{\infty}\frac{L \lambda \E\left[(S_i/\sqrt{L}) ^{2+\delta}\right]}{(1+\delta)(\sqrt{L}\nu)^{1+\delta}}\frac{1}{\tau^{1+\delta}}
\\
&=&
L^{-(1/2+\delta)}\frac{\lambda \eta_i }{(1+\delta)\nu^{1+\delta}} \sum_{\tau=1}^{\infty}\frac{1}{\tau^{1+\delta}},
\end{eqnarray*}
where the last inequality follows from
\[
\left(\sup_{0 \le t <1} \hat{d}^{(L)}_i(t) \right)^{2+\delta} \le \sum_{k=0}^{\Lambda^{(L)}(1)} \left(\frac{S_i}{\sqrt{L}}\right)^{2+\delta},
\] 
where $\Lambda^{(L)}(1)$ is the number of arrivals in $[0,L]$ and $S_i$ is the order size, which are independent of each other. 
The analysis above shows that
\[
\E\left[\sum_{\tau=0}^{\infty}\sup_{0 \le t <1}\left(\hat{d}^{(L)}_i(t)-\sqrt{L} \nu \tau\right)^+\right]
\le
\xi_1 L^{-\frac{\delta}{2(2+\delta)}}+\xi_2 L^{-(1/2+\delta)},
\]
where
\[
\xi_1=3 \lambda^{1/(2+\delta)}(1+\eta_i)~\mbox{and}~\xi_2=\frac{\lambda\eta_i}{(1+\delta)\nu^{1+\delta}}\sum_{\tau=1}^{\infty}\frac{1}{\tau^{1+\delta}}
\]
are both finite values that are independent of $L$. 
The lemma follows as a result. 
$\blacksquare$

\subsection*{Proof of Lemma \ref{lem:initconv}:}
It is easy to verify from the expressions in (\ref{eq:MSP2}) that because $\mathbb Y^{k(L)}(t)$ $(t \ge -L_k,1 \le k \le K)$ minimize (\ref{eq:SPL})-(\ref{eq:SPL}'), $\hat{\mathbb Y}^{k(L)}(t)$ $(t \ge -\hLL_k, 1 \le k \le K)$ minimize the same functions with the same linear transformation of inputs as follows: at each stage $k$ $(1 \le k \le K)$, let $\hat{\BD}^{k(L)}$ replace $\BD^{k(L)}$, $\x$ take the value of $\hBD(t+\hLL_{k+1}-1,t)$ instead of $\sBD(t+\LL_{k+1}-L,t)$, and $\hat{\mathbb Y}^{k'(L)}(t+\hLL_k-\hLL_{k'})$ replace $\mathbb Y^{k'(L)}(t+\LL_k-\LL_{k'})$ $(k<k'\le K)$.  

For each $j=1,\cdots,n$, applying  (\ref{eq:solbound}), there exists a constant $\tilde{\beta}$ such that
\begin{equation*}
%\label{eq:scalebound}
\hbY_j(-\hLL_{k_j}) \le \tilde{\beta} \sum_{i=1}^m \sum_{k=k_j}^K  \left(|\hat{D}^{(L)}_i(-1,-\hLL_k)|+\E[|\hat{D}^{(L)}_i(-\hLL_{k_j},0)|] \right).
\end{equation*}

For any constant $\nu>0$, let
\[
\tilde{\nu}=\frac{\nu/\tilde{\beta}}{2m(\bar{n}_k+\cdots+\bar{n}_K)}
\]
Since $\sD(t)$ ($t \ge -L$) is a stationary process and $0 \le \hLL_{k_j}\le 1$ $(\jrange)$,
\begin{eqnarray*}
	&~&\E\left[(|\hbY_j(-\hLL_{k_j})|-\sqrt{L} \nu)^+\right]
	\\
	&\le& 
 	\tilde{\beta} \sum_{i=1}^m \sum_{k=k_j}^K 
	\E\left[\left(|\hat{D}^{(L)}_i(-1,-\hLL_k)|-\sqrt{L} \tilde{\nu}\right)^++\left(\E[|\hat{D}^{(L)}_i(-\hLL_{k_j},0)|]-\sqrt{L}\tilde{\nu}\right)^+\right]
	\\
	&\le& 2 \tilde{\beta} (K-k_j) \sum_{i=1}^m \E\left [\sup_{t \ge 0} \left(|\hat{D}^{(L)}_i(0,t)|-\sqrt{L} \tilde{\nu}\right)^+\right].
\end{eqnarray*}
Applying Lemma \ref{lem:NBD} to the last expression in the above concludes the proof. 
$\blacksquare$

\subsection*{Proof of Theorem \ref{thm:unicom}:}

We use induction to prove that for $k=0,\cdots,K$ and any $\y^{k+1},\cdots,\y^K$ and $\x$, there exists a $\eta$, depending only on $\mathbf c$ and $A^{k'}$ $(1 \le k' \le K)$, such that
\begin{equation}
\label{eq:sandwidth}
0 \le \varphi^k_M(\y^{k+1},\cdots,\y^K,\x) -\varphi^k(\y^{k+1},\cdots,\y^K,\x) \le \eta \sum_{k'=1}^k \E[||\BD^{k'}-\BD^{k'}_M||_1],
\end{equation}
and the theorem follows immediately.

Since $\varphi^0_M(\y^1,\cdots,\y^K,\x)$ in (\ref{eq:MSP2}) and $\varphi^0(\y^1,\cdots,\y^K,\x)$ in (\ref{eq:MSPapp}) are the same LP, (\ref{eq:sandwidth}) holds when $k=0$. For the need of induction, observe that for for any $\x^a$ and $\x^b$ such that $\x^a \le \x^b$, 
\begin{equation}
\label{eq:LPsensitivity}
0 \le \varphi^0(\y^1,\cdots,\y^K,\x^a)-\varphi^0(\y^1,\cdots,\y^K,\x^b) \le \eta ||\x^a-\x^b||_1,
\end{equation}
where $\eta$ is a constant that depends only on $\mathbf c$ and $A^k$ $(1 \le k \le K)$.  The left inequality is directly implied by the formulation of the LP and the right inequality and the specification of $\eta$ are standard results of LP sensitivity analysis (e.g., Section 10.4 of \cite{Schrijver1986}). 

As our induction assumption, for $k=l$ $(0 \le l <K)$,   
\begin{equation}
\label{eq:induct1}
0 \le \varphi^l_M(\y^{l+1},\cdots,\y^K,\x) -\varphi^l(\y^{l+1},\cdots,\y^K,\x) \le \eta \sum_{k'=1}^l \E[||\BD^{k'}-\BD^{k'}_M||_1],
\end{equation}
and for any $\x^a$ and $\x^b$ where $\x^a \le \x^b$, 
\begin{equation}
\label{eq:induct2}
0 \le \varphi^l(\y^{l+1},\cdots,\y^K,\x^a)-\varphi^l(\y^{l+1},\cdots,\y^K,\x^b) \le \eta ||\x^a-\x^b||_1. 
\end{equation}
In both cases, $\eta$ is a constant that depends only on $\mathbf c$ and $A^k$ $(1 \le k \le K)$. 

We show below that both conditions are also satisfied when $k=l+1$. Let $\y_M^{(l+1)*}$ be an optimal solution of $\varphi^{l+1}_M(\y^{l+2},\cdots,\y^K,\x)$. 
\begin{eqnarray*}
	\varphi^{l+1}(\y^{l+2},\cdots,\y^K,\x)&=& \min_{\y^{l+1}}\{\mathbf h^{l+1}\cdot \y^{l+1}+\E[\varphi^{l}(\y^{l+1},\y^{l+2},\cdots,\y^K,\x+\BD^{l+1})]\}
	\\
	&\le& \mathbf h^{l+1}\cdot \y_M^{(l+1)*}+\E[\varphi^l(\y_M^{(l+1)*},\y^{l+2}, \cdots,\y^K,\x+\BD^{l+1})]
	\\
	&\le& \mathbf h^{l+1}\cdot \y^{(l+1)*}+\E[\varphi^l(\y_M^{(l+1)*},\y^{l+2},\cdots,\y^K,\x+\BD_M^{l+1})]
	\\
	&\le& \mathbf h^{l+1}\cdot \y^{(l+1)*}+\E[\varphi^l_M(\y_M^{(l+1)*},\y^{l+2}, \cdots,\y^K,\x+\BD^{l+1}_M)]
	\\
	&=&\varphi_M^{l+1}(\y^{l+2},\cdots,\y^K,\x),
\end{eqnarray*}
where the second inequality follows from induction assumption (\ref{eq:induct2}) with $\x^a=\BD^{l+1}_M$ and $\x^b=\BD^{l+1}$ and the third inequality follows from (\ref{eq:induct1}).

Let $\y^{(l+1)*}$ be an optimal solution to $\varphi^{l+1}(\y^{l+2},\cdots,\y^K,\x)$. Then
\begin{eqnarray*}
	\varphi^{l+1}_M(\y^{k+2},\cdots,\y^K,\x)&=& \min_{\y^{l+1}}\{\mathbf h^{l+1}\cdot \y^{l+1}+\E[\varphi^l_M(\y^{l+1},\y^{k+2},\cdots,\y^K,\x+\BD^{l+1}_M)]\}
	\\
	&\le& \mathbf h^{l+1}\cdot \y^{(l+1)*}+\E[\varphi^k_M(\y^{(l+1)*},\y^{k+2}, \cdots,\y^K,\x+\BD^{l+1}_M)]
	\\
	&\le& \mathbf h^{l+1}\cdot \y^{(l+1)*}+\E[\varphi^k(\y^{(l+1)*},\y^{k+2},\cdots,\y^K,\x+\BD^{l+1}_M)]+\eta\sum_{k'=1}^l \E[||\BD^{k'}-\BD^{k'}_M||_1]
	\\
	&\le& \mathbf h^{l+1}\cdot \y^{(l+1)*}+\E[\varphi^k(\y^{(l+1)*},\y^{k+2}, \cdots,\y^K,\x+\BD^{l+1})]+\eta\sum_{k'=1}^{l+1} \E[||\BD^{k'}-\BD^{k'}_M||_1]
	\\
	&=&\varphi^{l+1}(\y^{l+2},\cdots,\y^K,\x)+\eta\sum_{k'=1}^{l+1} \E[||\BD^{k'}-\BD^{k'}_M||_1],
\end{eqnarray*}
where the second inequality follows from (\ref{eq:induct1}) and third inequality follows from (\ref{eq:induct2})
with $\x^a=\BD_M^{l+1}$ and $\x^b=\BD^{l+1}$. Hence we have proved (\ref{eq:induct1}) holds when $k=l+1$. 

To prove (\ref{eq:induct2}) holds for $k=l+1$, for any $\x^a$ and $\x^b$ where $\x^a \le \x^b$,
\begin{eqnarray*}
	\varphi^{l+1}(\y^{l+2},\cdots,\y^K,\x^a) &=& \min_{\y^{l+1}}\{\mathbf h^{l+1}\cdot \y^{l+1}+\E[\varphi^l(\y^{l+1},\y^{l+2},\cdots,\y^K,\x^a+\BD^{l+1})]\}
	\\
	&\le& \min_{\y^{l+1}}\{\mathbf h^{l+1}\cdot \y^{l+1}+\E[\varphi^l(\y^{l+1},\y^{k+2},\cdots,\y^K,\x^b+\BD^{l+1})+\eta ||\x^a-\x^b||_1]\}
	\\
	&=&\varphi^{l+1}(\y^{l+2},\cdots,\y^K,\x^b)+\eta ||\x^a-\x^b||_1,
\end{eqnarray*}
where the inequality follows from the right inequality of (\ref{eq:induct2}). Following the same step and use the left inequality of (\ref{eq:induct2}),
\[
\varphi^{l+1}(\y^{l+2},\cdots,\y^K,\x^a) \ge \varphi^{l+1}(\y^{l+2},\cdots,\y^K,\x^b).
\]
This concludes the proof of the theorem statement.  

As a slight extension, in proving Theorem \ref{thm:sameobj}, we refer to this theorem to show that
\begin{equation}
\label{eq:Psiconv}
\lim_{M \rightarrow \infty} \Psi_{M,\bal}^K=\Psi_{\bal}
\end{equation}
where $\Psi_{\bal}$ is defined in (\ref{eq:uMSP1}) and $\Psi_{M,\bal}^K$ is defined by the same equation with $\BD^k$ replaced by $\BD^k_M$ $(1 \le k \le K)$.  It is straightforward to follow the same process as above,
with $\varphi^k_M()$ and $\varphi^k()$ replaced by $\Psi_{M,\bal}^k()$ and $\Psi^k_{\bal}()$ respectively $(1 \le k \le K)$,
to prove (\ref{eq:Psiconv}). The only difference is that the determination of $\eta$ in (\ref{eq:LPsensitivity}) needs to include the LHS of the constraint $\mathbf z \ge - \bal$. Since the value of $\eta$ is independent of $\bal$, which is on the RHS, the convergence is uniform over all $\bal$.  
$\quad \blacksquare$

\subsection*{Proof of Lemma \ref{lem:bondU}:}

The proof requires a substantial amount of articulation. To facilitate understanding, we will first prove the special case with $k=1$ before proving the general case in which $k$ can vary from 1 to $K$.

\subsubsection*{Special Case ($\mathbf{k=1}$):}
\label{subsubsec:kone}

Consider $\varphi^1_M(\mathbf y^2,\cdots,\mathbf y^K,\mathbf x)$, in which (\ref{eq:mapping1})-(\ref{eq:mapping4}) specialize to
\begin{align}
\begin{split}
\label{eq:k1const}
	\mathbf z(\bar{\omega}) ~~~~~~~~&\le \mathbf x+\underline{\mathbf D}^1_M(\bar{\omega}), ~~~~~~\bar{\omega} \in \Omega_1^0,
	\\
	A^1 \mathbf z(\bar{\omega})~-~\mathbf y^1 & \le  \mathbf 0,  ~~~~~~~~~~~~~~~~~~~\bar{\omega}  \in \Omega_1^0,
	\\
	A^k \mathbf z(\bar{\omega}) ~~~~~~~~& \le \mathbf y^k,~~~~~~~~~~~~~~~~~~\bar{\omega} \in \Omega_1^0,~1 < k \le K.
\end{split}
\end{align}
The coefficient matrix on the LHS of constraints can be expressed as
\begin{equation}
\label{eq:stage1mat}
\FCBA=(\FCBH,\FCBE)~~~\mbox{where}~~~~\FCBH=\left(\begin{array}{cccc} H_1 & & &\\  & \ddots & & \\ & & H_N & \end{array} \right).
\end{equation}
All components of submatrix $\FCBH$ are zero except for those in block matrices on the diagonal (here and below we will deviate slightly from the standard definition by referring to a submatrix on the diagonal as a block even if it is not a square one).
The number of blocks, $N$, is the number of elements in $\Omega_1^0$. Each  block $H_i$ $(1 \le i \le N)$ is composed of coefficients associated with $\mathbf z(\bar{\omega})$ of a particular $\bar{\omega}$, and thus has $m+n_1+...+n_K$ rows and $m$ columns. 
In particular,
\begin{equation}
\label{eq:defH}
H_1=...=H_N=\left (\begin{array}{c}
I
\\
A^1
\\
.....
\\
A^K
\end{array}
\right).
\end{equation}
The matrix $\FCBE$ is composed of coefficients of $\mathbf y^1$. It has $n_1$ columns (the number of components in $\mathbf y^1$) and its components are either $0$ or $-1$.  Since every $H_i$ $(1 \le i \le N)$ has an identity submatrix, $\FCBE$ has a negative identity submatrix, and these submatrices are on different rows, $\FCBA$ in (\ref{eq:stage1mat}) has full column rank.

In general, we define a matrix $\FCBA$ as a Finite Coupling Block Matrix (FCBM) with characterization numbers $(m,n_1)$ if:
\begin{enumerate}
\item
With necessary permutations of rows and columns, $\FCBA$ can fit into the pattern defined by (\ref{eq:stage1mat});
\item 
$\FCBA$ is of full column rank;
\item
$\FCBH$ has a finite number of blocks (denoted by $N$); all blocks, $H_i$ $(1 \le i \le N)$, have the same finite number of columns (denoted by $m$), but can have different number of rows;
\item
$\FCBE$ has a finite number of columns (denoted by $n_1$), and its components are either $0$ or $-1$.
\end{enumerate}

Let $\hat{\FCBA}$ be a maximal nonsingular submatrix of a FCBM $\FCBA$. Since $\FCBA$ is of full column rank, $\hat{\FCBA}$ has the same number of columns as $\FCBA$.  Referring to (\ref{eq:stage1mat}), since each block $H_i$ $(1 \le i \le N)$ in $\FCBA$ is of full column rank (or otherwise $\FCBA$ will not be so), $\hat{\FCBA}$ must contain at least $m$ rows from each block. To be a square matrix, $\hat{\FCBA}$ must draw more rows than columns from some blocks, and the total number of these extra rows is exactly $n_1$ (the number of columns of $\FCBE$). This means that the number of blocks with more than $m$ rows in $\hat{\FCBA}$ is somewhere between $1$ and $n_1$. Correspondingly we can divide $H_i$ $(1 \le i \le N)$ into two groups: group 1 are those with exactly $m$ rows in $\hat{\FCBA}$, each forms a $m$-dimensional nonsingular submatrix. Group 2  are those with more than $m$ rows in $\hat{\FCBA}$. 
Hence by permuting rows and columns, we can write $\hat{\FCBA}$ as
\begin{equation}
\label{eq:nsFCBM}
\hat{\FCBA}= \left (\begin{array}{lll} \FCBH_1 & ~~\FCBE_1 & ~~\mathbf 0 \\
\mathbf 0 & ~~\FCBE_2 & ~~\FCBH_2
\end{array}
\right),
\end{equation}
where $\FCBH_1$ is a diagonal block matrix, where each block is a $m$-dimensional nonsingular matrix, $\FCBH_2$ is a also a diagonal block matrix with no more than $n_1$ blocks, where each block has more rows than columns, and $\left (\begin{array}{c} \FCBE_1 \\ \FCBE_2 \end{array}\right)$ is a submatrix of $\FCBE$ with rearranged rows. Clearly $\hat{\FCBA}$ is a FCBM itself. 

Under the above specifications, Lemma \ref{lem:bondU} holds when $k=1$ if the following statement is true:

%\begin{lemma}
%\label{lem:stage1FCBM}
\emph{For any FCBM $\FCBA$ with characterization numbers $(m,n_1)$, where $m$ and $n_1$ can be any positive integer, let $\hat{\FCBA}$ be a maximal nonsingular submatrix of $\FCBA$. Let $\mathcal V$ be the set of  values of elements of $\FCBA$. Let $\mathbf u$ be any vector with the same number of components as the number of rows in $\hat{\FCBA}$.
Then there exists a constant $\kappa$, depending only on $\mathcal V$ and $(m,n_1)$, such that
\begin{equation}
\label{eq:CBMbound}
||\mathbf u||_1 \le \kappa ||\mathbf u \hat{\FCBA}||_1.
\end{equation}
}
%\end{lemma}
%\subsection*{Proof of Lemma \ref{lem:stage1FCBM}:}

To prove the above statement, let $\mathcal M_1$ be the set of indexes of nonsingular blocks in $\FCBH_1$.  We partition $\mathbf u$ into $\mathbf u^1$ and $\mathbf u^2$ such that when multiplying $\mathbf u$ with $\hat{\FCBA}$, components of $\mathbf u^1$ are multiplied with $(\FCBH_1~\FCBE_1)$ and components of $\mathbf u^2$ are multiplied with $(\FCBH_2~\FCBE_2)$.
Similarly, we partition $\mathbf u^1$ into $\mathbf u^1_i$ ($i \in \mathcal M_1$) where components of $\mathbf u^1_i$ are multiplied with block $H_i$ in $\FCBH_1$. 
For each $i \in \mathcal M_1$, $H_i$ is nonsingular, thus
\begin{equation}
\label{eq:u1con}
||\mathbf u^1||_1 = \sum_{i \in \mathcal M_1} ||\mathbf u^1_i||_1=\sum_{i \in \mathcal M_1} ||\mathbf u^1_i H_i H_i^{-1}||_1 \le \kappa_1 \sum_{i \in \mathcal M_1} ||\mathbf u^1_i H_i||_1.
\end{equation}
Here $H_i$ $(i \in \mathcal M_1)$ are $m$-dimensional nonsingular matrices, with their components taking values from the finite set $\mathcal V$, and $\kappa_1$ is the maximum component value of all possible $H_i^{-1}$ $(i \in \mathcal M_1$). 

Define $\FCBG_2 \doteq (\FCBE_2~ \FCBH_2)$. Since both  $\hat{\FCBA}$ and $\FCBH_1$ are nonsingular matrices, $\FCBG_2$ is also nonsingular.  Each block in $\FCBH_2$ has $m$ columns. For $\FCBG_2$ to be a square matrix, $\FCBH_2$ has $n_1$ more rows than columns, so the dimension of $\FCBG_2$ is between $n_1+m$ if $\FCBH_2$ has one block with $n_1+m$ rows and $(m+1) n_1$ if $\FCBH_2$ has $n_1$ blocks, each with $m+1$ rows. 
Let $\kappa_2$ be the maximum component of all possible formation of $\FCBG_2^{-1}$, where $\FCBG_2$ is a nonsingular matrix with dimension between $n_1+m$ and $(m+1)n_1$ and component values drawn from the finite set $\mathcal V$. Then
\begin{equation}
\label{eq:u2con}
||\mathbf u^2||_1 = ||\mathbf u^2\FCBG_2 \FCBG_2^{-1} ||_1 \le \kappa_2 ||\mathbf u^2 \FCBG_2 ||_1 =\kappa_2 \left ( ||\mathbf u^2 \FCBH_2||_1+||\mathbf u^2 \FCBE_2||_1 \right).
\end{equation}
It follows that
\begin{align}
\label{eq:uG}
\begin{split}
||\mathbf u \hat{\FCBA}||_1&= \sum_{i \in \mathcal M_1} ||\mathbf u^1_i H_i ||_1+ ||\mathbf u^2 \FCBH_2||_1+|| \mathbf u^1 \FCBE_1+\mathbf u^2 \FCBE_2||_1
\\
& \ge
\sum_{i \in \mathcal M_1} ||\mathbf u^1_i H_i ||_1+||\mathbf u^2 \FCBH_2||_1+||\mathbf u^2 \FCBE_2||_1-||\mathbf u^1 \FCBE_1||_1  
\\
&\ge  \sum_{i \in \mathcal M_1} ||\mathbf u^1_i H_i ||_1+ ||\mathbf u^2 \FCBH_2||_1+
||\mathbf u^2 \FCBE_2||_1- n_1 ||\mathbf u^1||_1
\\
&\ge  \frac{||\mathbf u^1||_1}{\kappa_1}+ \frac{||\mathbf u^2||_1}{\kappa_2}- n_1 ||\mathbf u^1||_1,
\end{split}
\end{align}
where the second inequality holds because components of $\mathbb E_1$ are either $0$ or $-1$, and the last inequality follows from (\ref{eq:u1con}) and (\ref{eq:u2con}). 
Also by (\ref{eq:u1con}),
\[
||\mathbf u^1||_1 \le \kappa_1 \sum_{i \in \mathcal M_1} ||\mathbf u^1_i H_i||_1 \le \kappa_1  ||\mathbf u \hat{\FCBA}||_1,
\]
where the second inequality follows from the structure of $\hat{\FCBA}$ shown in (\ref{eq:nsFCBM}). Thus (\ref{eq:uG}) implies that 
\begin{equation}
\label{eq:conk1}
||\mathbf u||_1 =  ||\mathbf u^1||_1 +||\mathbf u^2||_1 \le (\kappa_1 \vee \kappa_2) \left (||\mathbf u \hat{\FCBA}||_1+n_1 ||\mathbf u^1 ||_1\right) \le (\kappa_1 \vee \kappa_2) (1+n_1 \kappa_1) ||\mathbf u \hat{\FCBA}||_1,
\end{equation}
and (\ref{eq:CBMbound}) is satisfied with
\[
\kappa=(\kappa_1 \vee \kappa_2) (1+n_1 \kappa_1).
\]

\subsubsection*{General Case ($\mathbf {1 \le k \le K}$):}
\label{subsubsec:kgeneral}

We can extend the proof of (\ref{eq:limitnorm}) to cases with $k=1,\cdots, K$ by induction.  To do that, we first generalize the definition of FCBM with the following recursive characterizations:
\begin{itemize}
\item
a matrix $\FCBA^0$ is stage-0 FCBM with characterization number $m$ if it is a finite-dimensional matrix with $m$ columns and full column rank;
\item
a matrix $\FCBA^k$ is a stage-$k$ FCBM $(k \ge 1)$ with characterization numbers $(m,n_1,\cdots,n_k)$ if it is of full column rank, and with necessary permutations of rows and columns, can be written as
\begin{equation}
\label{eq:stagekmat}
\FCBA^k=(\FCBH^k, \FCBE^k),
\end{equation}
where
\begin{enumerate}
\item 
$\FCBH^{k}$ has non-zero entries only in diagonal blocks, the number of which, denoted by $N_k$, can take any finite value. Specifically,
\begin{equation}
\label{eq:FCBHk}
\FCBH^k=\left(\begin{array}{cccc} \FCBA^{k-1}_1 & & &\\  & \ddots & & \\ & & \FCBA^{k-1}_{N_k} & \end{array} \right);
\end{equation}
where each block $\FCBA^{k-1}_i$ $(1 \le i \le N_k$) is a stage-$(k-1)$ FCBM with characterization numbers $(m,n_1,\cdots,n_{k-1})$;
blocks might not be identical, but they all have the same number of columns;
\item 
$\FCBE^k$ has a finite number of columns (denoted by $n_k$) and its components are either 0 or -1.
\end{enumerate} 
\end{itemize}

For $k=1,\cdots,K$, the LHS constraint matrix of $\varphi_M^k(\y^{k+1},\cdots,\y^K,\x)$ is a FCBM. To show this, when $k=2$, (\ref{eq:mapping1})-(\ref{eq:mapping4}) specialize to
\begin{eqnarray}
\label{eq:constz1}
	\mathbf z(\bar{\omega})~~~~~~~~~~~~~~~~~~~~&\le& \mathbf x+\underline{\mathbf D}^2_M(\bar{\omega}), ~~\bar{\omega} \in \Omega_2^0,
	\\
	\label{eq:constz2}
	A^1 \mathbf z(\bar{\omega})~-\mathbf y^1 (\bar{\omega}') ~~~~~~~~& \le & \mathbf 0, ~~~~~~~~~~~~~~~\bar{\omega}' \in \Omega_2^1,~\bar{\omega} \in \Omega_2^0, ~\bar{\omega}' \sqsubset \bar{\omega},
	\\
	\label{eq:constz3}
	A^2 \mathbf z(\bar{\omega}) ~~~~~~~~~~~~~-\mathbf y^2~& \le & \mathbf 0,~~~~~~~~~~~~~~~\bar{\omega} \in \Omega_2^0,
	\\
	\label{eq:constz4}
	A^k \mathbf z(\bar{\omega}) ~~~~~~~~~~~~~~~~~~~~~& \le & \mathbf y^k,~~~~~~~~~~~~~~\bar{\omega} \in \Omega_2^0,~2 < k \le K.
\end{eqnarray}
With necessary permutations of rows and columns, its LHS constraint matrix can be written as
\[
\FCBA^2=(\FCBH^2,\FCBE^2)=\left(\begin{array}{cccc} \FCBA^1_1 & & &\\  & \ddots & & \FCBE^2 \\ & & \FCBA^1_{N_2} & \end{array} \right).
\]
Here $\FCBE^2$ is composed of coefficients of $\y^2$ and $N_2$ is the number of elements in $\Omega_2^1$. Each $\FCBA^1_i$ $(1 \le i \le N_2)$ corresponds to a particular $\bar{\omega}' \in \Omega_2^1$, and
\[
\FCBA_i^1=\left(\begin{array}{cccc} \FCBA^0_1 & & &\\  & \ddots & & \FCBE^1 \\ & & \FCBA^0_{N_1} & \end{array} \right),
\]
where $\FCBE^1$ is composed of coefficients of $\y^1(\bar{\omega}')$, $N_1$ is the number of elements $\bar{\omega}$ in $\Omega_2^0$ such that $\bar{\omega}' \sqsubset \bar{\omega}$. Corresponding to a particular $\bar{\omega}$, $\FCBA^0_i$ ($1 \le i \le N_1$) are the same as $H_i$s in (\ref{eq:defH}), and composed of coefficients associated with $\z(\bar{\omega})$. Thus by definition, $\FCBA^1_i$ $(1 \le i \le n_1)$ is a stage-1 FCBM with characterization numbers $(m,n_1)$. Moreover, since $\FCBA^0_i$ $(1 \le i \le N_1$) contains an identity matrix (constraints (\ref{eq:constz1})), each $\FCBE^1$ $(1 \le i \le N_2)$ contains a negative identity matrix (constraints (\ref{eq:constz2})), $\FCBE^2$ contains a negative identity matrix (constraints (\ref{eq:constz3})), and these submatrices have no overlapping rows, $\FCBA^2$ has full column rank. Hence $\FCBA^2$ is a FCBM with characterization numbers $(m,n_1,n_2)$. 

For general $k$ $(1 \le k \le K)$, the same can be shown by induction. Referring to (\ref{eq:mapping1})-(\ref{eq:mapping4}). Let $N_k$ be the number of elements in $\Omega_k^{k-1}$. Then we can make the induction assumption that coefficients of $\z(\bar{\omega})$ $(\bar{\omega} \in \Omega_k^0)$ and $\y^{k'}(\bar{\omega}')$ $(\bar{\omega}' \sqsubset \bar{\omega}, \bar{\omega}' \in \Omega_k^{k'},\bar{\omega} \in \Omega_k^0, 1 \le k' <k)$ are given by a number of $N_k$ stage-$(k-1)$ FCBMs  with characterization numbers $(m,n_1,\cdots,n_{k-1})$, which we denote by $\FCBA^{k-1}_i$ $(1 \le i \le N_k)$. Coefficients of $\y^k$, which are either -1 or 0, are given by matrix $\FCBE^k$, which has $n_k$ columns.  Thus with necessary permutations of rows and columns, the LHS constraint matrix of $\varphi_M^k(\y^{k+1},\cdots,\y^K,\x)$ can be written as
\[
\FCBA^k=\left(\begin{array}{cccc} \FCBA^{k-1}_1 & & &\\  & \ddots & & \FCBE^{k} \\ & & \FCBA^{k-1}_{N_k} & \end{array} \right).
\]
Since every column in constraints (\ref{eq:mapping1}), (\ref{eq:mapping2}), and (\ref{eq:mapping3}) are covered by an identity or negative identity matrix and rows of these matrices do not overlap,  by definition $\FCBA^k$ is a FCBM  with characterization numbers $(m,n_1,\cdots,n_k)$.

Let $\hat{\FCBA}^k$ be a maximal nonsingular submatrix of $\FCBA^k$ $(1 \le k \le K)$. Then following the same reasoning that leads to (\ref{eq:nsFCBM}), $\hat{\FCBA}^k$ and $\FCBA^k$ must have the same number of columns since $\FCBA^k$ has  full column rank. By (\ref{eq:stagekmat})-(\ref{eq:FCBHk}), for $i=1,\cdots,N_k$, $\hat{\FCBA}^k$ contains either a maximal nonsingular submatrix of $\FCBA_i^{k-1}$ (the rank of the submatrix equals the column rank of $\FCBA_i^{k-1}$) or a submatrix that includes every column of $\FCBA_i^{k-1}$ and has more rows than columns. The total number of extra rows is $n_k$ for $\hat{\FCBA}^k$ to be nonsingular. Hence $\hat{\FCBA}^k$ can be written as
\begin{equation}
\label{eq:nsFCBMgen}
\hat{\FCBA}^k= \left (\begin{array}{lll} \FCBH_1^k & ~~\FCBE^k_1 & ~~\mathbf 0 \\
\mathbf 0 & ~~\FCBE^k_2 & ~~\FCBH_2^k
\end{array}
\right),
\end{equation}
where $\FCBH_1^k$ is a diagonal block matrix, with each block being a nonsingular submatrix of $\FCBA_i^{k-1}$ for some $i$ $(1 \le i \le N_k$),  $\FCBH_2^k$ is a non-square diagonal block matrix, where each block is a submatrix of $\FCBA^{k-1}_i$ for some other $i$ $(1 \le i \le N_k)$, and has more rows than columns, and $\left (\begin{array}{c} \FCBE^k_1 \\ \FCBE^k_2 \end{array} \right)$ is composed of a subset of rows in $\FCBE^k$. The number of blocks $\FCBH_2^k$ ranges from $1$ to $n_k$.

Under the above specifications, the lemma holds for all $k$ $(1 \le k \le K)$ if the following statement is true:

%\begin{lemma}\label{lem:CBMboundgen}
\emph{For $k=1,\cdots, K$, let $\FCBA^k$ be a stage-$k$ FCBM with characterization numbers $(m,n_1,\cdots,n_k)$, where $m, n_1,\cdots,n_k$ can be any positive integers. Let $\hat{\FCBA}^k$ be a maximal nonsingular submatrix of $\FCBA^k$. 
Let $\mathcal V$ be the set of component values in $\FCBA^k$. Let $\mathbf u$ be any vector with the same number of components as the number of rows in $\hat{\FCBA}^k$.
Then there exists a constant $\kappa$, depending only on $\mathcal V$ and $(m,n_1,\cdots,n_k)$, such that
\begin{equation}
\label{eq:CBMboundgen}
||\mathbf u||_1 \le \kappa ||\mathbf u \hat{\FCBA}^k||_1.
\end{equation}
}
%\end{lemma}

%\subsection*{Proof of Lemma \ref{lem:CBMboundgen}:}
We have already shown in (\ref{eq:conk1}) that the above statement is true for $k=1$, and thus can use induction to prove the same for $k>1$:
%Using Lemma \ref{lem:stage1FCBM} as the starting point, we use induction to prove Lemma \ref{lem:CBMboundgen} for cases where $k>1$. 

Assume that (\ref{eq:CBMboundgen}) holds for any maximal nonsingular submatrix of a stage-$(k-1)$ FCBM $(k>1)$ with any finite and positive characterization numbers. 
Let $\FCBA^k$ be a stage-$k$ FCBM with characterization numbers $(m,n_1,\cdots, n_{k-1}, n_k)$. Recall that by permuting rows and columns, a maximal nonsingular matrix of $\FCBA^k$ can be written as
\begin{equation}
\label{eq:blockdiv}
\hat{\FCBA}^k = \left (\begin{array}{lll} \FCBH_1^k & ~~\FCBE^k_1 & ~~\mathbf 0 \\
\mathbf 0 & ~~\FCBE^k_2 & ~~\FCBH_2^k
\end{array}
\right),
\end{equation}
where $\FCBH^k_1$ is a diagonal block matrix, where each block is a nonsingular stage-$(k-1)$ FCBM with characterization numbers $(m,n_1,\cdots,n_{k-1})$, $\FCBH^k_2$  is a diagonal matrix, where each block is also a stage-$(k-1)$ FCBM with characterization numbers $(m,n_1,\cdots,n_{k-1})$ but has more rows than columns, and  $\left (\begin{array}{c} \mathbb E_1^k \\ \mathbb E_2^k \end{array} \right)$ is a matrix with $n_k$ columns and either -1 or 0 as values of its components. 
Let  blocks in $\FCBH^k_1$ be indexed and let $\mathcal M^k_1$ be the index set. Denote these blocks by $\hat{\FCBA}^{k-1,i}$ ($i \in \mathcal M^k_1$), where the first superscript shows the stage and the second one identifies a particular block. 

We partition the vector $\mathbf u$ into $\mathbf u^1$ and $\mathbf u^2$ such that in the product of $\mathbf u$ and $\hat{\FCBA}^k$, components of $\mathbf u^1$ and $\mathbf u^2$ are multiplied with rows in $(\FCBH_1^k~\FCBE^k_1)$ and $(\FCBE^k_2 ~\FCBH_2^k)$ respectively. 
We also partition $\mathbf u^1$ into $\mathbf u^1_i$ $(i \in \mathcal M^k_1)$ such that components of $\mathbf u^1_i$ are multiplied with rows in $\hat{\FCBA}^{k-1,i}$ ($i \in \mathcal M^k_1$).

Since $\hat{\FCBA}^{k-1,i}$ ($i \in \mathcal M^k_1$) is a nonsingular stage-$(k-1)$ FCBM, its maximal nonsingular submatrix is the matrix itself. Therefore by the induction assumption, there exists a constant $\kappa_1$, depending only on $(m,n_1,\cdots,n_{k-1})$ and $\mathcal V$, such that 
\begin{equation}
\label{eq:u1}
||\mathbf u^1||_1 =\sum_{i \in \mathcal M^k_1} ||\mathbf u^1_i||_1 \le \kappa_1 \sum_{i \in \mathcal M^k_1} ||\mathbf u^1_i \hat{\FCBA}^{k-1,i}||_1.
\end{equation}

To bound $\mathbf u^2$, recall that $\FCBH^k_2$ is a diagonal block matrix that contains no more than $n_k$ (non-square) blocks. Suppose that the actual number is $N'$ $(1 \le N' \le n_k)$.  As noted above, each block is a 
%maximal submatrix of a 
stage-$(k-1)$ FCBM with full column rank.
% and thus is itself a stage-$(k-1)$ FCBM. 
Let $\mathcal M^k_2$ be  the index set of these blocks and denote them by $\FCBA^{k-1,i}$ ($i \in \mathcal M_2^k$). 
Following the definition in (\ref{eq:stagekmat})-(\ref{eq:FCBHk}),
\[
\FCBA^{k-1,i}=(\FCBH^{k-1,i}~\FCBE^{k-1,i}),~~i \in \mathcal M^k_2. 
\]
For the ease of presentation, we let indexes in $\mathcal M^k_2$ take integer values $i=1,\cdots, N'$. Thus the lower right submatrix in (\ref{eq:blockdiv}) can be written as
\begin{align*}
\begin{split}
(\FCBE^k_2~\FCBH^k_2)&=\left (
	\begin{array}{c c c c c c c c c c }
		& \FCBA^{k-1,1} ~& & & & & & & &  \\
		& & & & &~~~~~~ ........& &  & &   \\
		\FCBE^k_2~~~ & & & & & & & \FCBA^{k-1,i}   & &  \\
		& & & &  & & & & ~~~~~~~~...........&    \\
		& &  & & & & & & & ~\FCBA^{k-1,N'}  \\
	\end{array}
	\right)
	\\
	&=
	\left (
	\begin{array}{c c c c c c c c c c }
		& \FCBH^{k-1,1}~ &~\FCBE^{k-1,1}&  & &  & & &   \\
		& & & .......& &  & & &  \\
		 \FCBE^k_2~~~ & & & &~~ \FCBH^{k-1,i} ~& \FCBE^{k-1,i}& & & & \\
		& & & & &  & ........& & \\
		& &  & & & & &~ \FCBH^{k-1,N'} & &\FCBE^{k-1,N'} \\
	\end{array}
	\right).
\end{split}
\end{align*}
Applying (\ref{eq:stagekmat})-(\ref{eq:FCBHk}) recursively, every $\FCBH^{k-1,i}$ is a diagonal block matrix, where each block is a stage-$(k-2)$ FCBM with characterization numbers $(m,n_1,\cdots,n_{k-2})$. Let $M(i)$ be the number of blocks in $\FCBH^{k-1,i}$ $(i \in \mathcal M^k_2)$. Denote a block in $\FCBH^{k-1,i}$ by $\FCBA_j^{k-2,i}$ $(1 \le j \le M(i),i \in \mathcal M^k_2)$. Then with necessary permutations of rows and columns, $(\FCBE^k_2~\FCBH^k_2)$ can be written in a more expanded form as
\begin{align*}
\begin{split}
\left (
\begin{array}{c c c c c c c c c c c |c c c c c c}
\FCBA^{k-2,1}_1~&     &    & & & & & & & & & & & & & &  \\
			  & ....&    & & & & & & & & &  \FCBE^{k-1,1} & & & & &   \\
                &     & \FCBA^{k-2,1}_{M(1)} & & & & & & & & & & & & & &  \\ 
                 &     &                     & ....& & & & & & & & & ....& & & & \\ 
                 &     &    &     & \FCBA^{k-2,i}_1 & & & & & & & & & & & & \\ 
                  &     &    &     &                   &.... & & & & & & &  & \FCBE^{k-1,i} & & &~~~~ \FCBE^k_2 \\
			  & & & & & &  \FCBA^{k-2,i}_{M(i)} & & & & & & & & & &  \\ 
 & & & & & & &..... & & & & &  & &... & &  \\ 
& & & & & & & & \FCBA^{k-2,N'}_1 & & & & & & & &  \\
& & & & & & & &   &....  & & &  & & & \FCBE^{k-1,N'} & \\
   & & & & & & &  & &  & \FCBA^{k-2,N'}_{M(N')} & & & & &  & \\
\end{array}
\right).
\end{split}
\end{align*}
The submatrix to the left of $|$ is a diagonal block matrix, where each block is a stage-$(k-2)$ FCBM. The submatrix on the right has either $0$ and $-1$ as its component values. Since $\FCBE^k_2$ has $n_k$ columns, each $\FCBE^{k-1,i}$ $(i \in \mathcal M_2^k)$ has $n_{k-1}$ columns, and the number of components of $\mathcal M^k_2$ is $N' \le n_k$, the latter submatrix has a finite number of $N'n_{k-1}+n_k$ columns. Also by (\ref{eq:blockdiv}), $(\FCBE^k_2~\FCBH^k_2)$ has full rank. Hence $(\FCBE^k_2~\FCBH^k_2)$ fits the definition as a nonsingular stage-$(k-1)$ FCBM with characterizations numbers of $(m,n_1,\cdots, n_{k-2},N'n_{k-1}+n_k)$.

By the induction assumption, there exists a constant $\kappa_2$, which depends only on $\mathcal V$ and  $(m,n_1,\cdots,n_{k-2},N'n_{k-1}+n_k)$ (and thus $(m,n_1,\cdots,n_{k-1},n_k$) since $1 \le N' \le n_k$), such that 
\begin{equation}
\label{eq:u2}
||\mathbf u^2||_1  \le \kappa_2 ||\mathbf u^2 \left (\FCBE^k_2~\FCBH^k_2\right)||_1= \kappa_2 \left(||\mathbf u^2 \FCBE^k_2||_1+||\mathbf u^2 \FCBH_2^k||_1\right).
\end{equation}
The proof can then be completed by observing that 
\begin{eqnarray*}
	||\mathbf u \hat{\FCBA}^{k}||_1 &=& \sum_{i \in \mathcal M^k_1} ||\mathbf u^1_i \hat{\FCBA}^{k-1,i}||_1 + ||\mathbf u^1 \FCBE^k_1+ \mathbf u^2 \FCBE^k_2||_1+||\mathbf u^2 \FCBH^k_2||_1
	\\
	&\ge & 
	\sum_{i \in \mathcal M^k_1} ||\mathbf u^1_i \hat{\FCBA}^{k-1,i}||_1 + ||\mathbf u^2 \FCBE^k_2||_1-||\mathbf u^1\FCBE^k_1||_1+||\mathbf u^2 \FCBH^k_2||_1
	\\
	&\ge&
	\sum_{i \in \mathcal M^k_1} ||\mathbf u^1_i \hat{\FCBA}^{k-1,i}||_1 + ||\mathbf u^2 \FCBE^k_2||_1-n_k ||\mathbf u^1||_1+||\mathbf u^2 \FCBH^k_2||_1
	\\
	&\ge&
	\sum_{i \in \mathcal M^k_1} ||\mathbf u^1_i \hat{\FCBA}^{k-1,i}||_1 + ||\mathbf u^2 \FCBE^k_2||_1-n_k \kappa_1 \sum_{i \in \mathcal M^k_1} ||\mathbf u^1_i \hat{\FCBA}^{k-1,i}||_1+||\mathbf u^2 \FCBH^k_2||_1
	\\
	&\ge&
	\frac{||\mathbf u^1||_1}{\kappa_1}+\frac{||\mathbf u^2||_1}{\kappa_2} - n_k \kappa_1 ||\mathbf u \hat{\FCBA}^k||_1,
\end{eqnarray*}
where the second inequality holds because $\FCBE^k_1$ has $n_k$ columns, the third inequality comes from (\ref{eq:u1}), and the last inequality is implied by (\ref{eq:u1}), (\ref{eq:u2}), and the first equality in the above.
Thus
\[
||\mathbf u||_1=||\mathbf u^1||_1+||\mathbf u^2||_2  \le \kappa_1 \vee \kappa_2 (1+n_k \kappa_1) ||\mathbf u \hat{\FCBA}^k||_1.~~~~~~~~~~~~\blacksquare
\]
\section*{Appendix B: Illustration of the Inventory Policy}
\label{sec:illustration}

We use a very simple example, the $N$ system, to illustrate our inventory policy developed in Section \ref{sec:policy}. Figure \ref{fig:spesys} in Section \ref{sec:formulation} shows that the system has two products and two components. One unit of component $1$ is used by both products, and one unit of component $2$ is used by product $2$ only. Components have different lead times and $L_1<L_2$. Thus $K=2$. We also assume that $c_1>c_2$.

The formulation of the SP (\ref{eq:MSP2}) for this system specializes to 
\begin{align*}
\begin{split}
\varphi^2 &=\min_{y_2} \left \{h_2 y_2 + \E[\varphi^1(y_2,\BD^2)]\right \} \\
\varphi^1(y_2,\x)&= \min_{y_1} \left \{h_1 y_1 +\E[\varphi^0(y_1,y_2, \x+\BD^1)]\right\}
\\
\varphi^0(y_1,y_2,\x)&=-\max_{z_1,z_2} \left \{c_1 z_1+c_2 z_2 \left | \right . z_1 \le x_1, z_2 \le x_2, z_1+z_2 \le y_1, z_2 \le y_2\right\}.
\end{split}
\end{align*}
To maximize $\varphi_0(y_1,y_2,\x)$,
\[
z_1^*=x_1 \quad \mbox{and} \quad z_2^*=x_2 \wedge (y_1-x_1) \wedge y_2.
\]
Correspondingly, with $y_2$ and $\x$ given, $y_1^*$ is chosen to minimize
\begin{equation}
\label{eq:s1spe}
h_1 y_1 -c_1 \E[x_1+D_1^1]-c_2 \mathbf E\left[(x_2+D_2^1) \wedge (y_1-x_1-D_1^1) \wedge y_2\right],
\end{equation}
where the expectation is taken with respected to $\BD^1=(D_1^1,D_2^1)$. Substituting $\x$ with the realization of $\BD^2$ in the resulting optimal objective value $\varphi^1(y_2,\x)$ and taking the expectation over $\BD^2$, $y_2^*$ is chosen to minimize
\begin{equation}
\label{eq:s2spe}
h_2 y_2 +\E[\varphi^1(y_2,\mathbf D^2)].
\end{equation}
It follows that the replenishment policy prescribed in Table \ref{tbl:replenishment} specializes to the following actions:
\begin{enumerate}
\item
Let $y_2^*$ be the (constant) inventory position target for component $2$ and follow a base stock policy to keep the actual inventory position at that level.
\item
At time $t$, where $t$ is either $-L_1$ or any time afterwards when there is a demand arrival:
\begin{enumerate}
\item 
\label{itm:rep}
choose $y_1^*$ to minimize (\ref{eq:s1spe}) with $\x=\BD(t-L_2+L_1,t)$,set component $1$'s inventory position target to
\[
\tIP_1(t)=y_1^*- [D_1(t-L_2+L_1,t)+D_2(t-L_2+L_1)],
\]
and order $(\lceil \tIP_1(t) \rceil-IP^-_1(t))^+$ units of that component;
\item
schedule a future update of the inventory position target and a possible new replenishment at time $t'=t-L_1+L_2$. Repeat Step \ref{itm:rep} when the system reaches  $t=t'$  
\end{enumerate}
\end{enumerate}

\bigskip

For the allocation policy, the procedure prescribed in Table \ref{tbl:allproc} specializes to the following actions:

\medskip
\noindent
At each time $t$, where $t$ is either $0$ or any time afterwards when there is a demand arrival or some previously-ordered component is received:
\begin{enumerate}
\item
set backlog targets at 
\[
\bB_1(t)=0 \quad \mbox{and} \quad \bB_2(t)=(D_1(t-L_1,t)+D_2(t-L_1,t)-IP_1(t-L_1))^+ \vee (D_2(t-L_2,t)-y_2^*)^+,
\]
which correspond to the optimal solution to
\[
\min_{\mathbf B \ge 0} \left \{c_1 B_1+c_2 B_2 \left | \right. B_1+B_2 \ge Q_1(t), B_2 \ge Q_2(t)\right\},
\]
using
\[
Q_1(t)=D_1(t-L_1,t)+D_2(t-L_1,t)-IP_1(t-L_1)\quad \mbox{and} \quad Q_2(t)= D_2(t-L_2,t)-y_2^*.
\]
\item
use component $1$ to serve demand $1$ until either the component runs out or $B_1(t)=0$. In the latter case, use the remaining amount of component 1 and the available amount of component 2 to serve as much demand $2$ as possible.  
\end{enumerate}
Notice that in this process, the backlog target $\bB_1(t)$ does not depend on system state and $\bB_2(t)$ is the higher shortage level of the two components. Thus the above procedure can be simplified to a priority policy, i.e., use available components to serve as much demand as possible, and when there is not sufficient component $1$ to serve both demands, use the component for product $1$ first. 

\section*{Appendix C: Formulation and Perturbation of LP (\ref{eq:finitLPobj})-(\ref{eq:mapping4})}
The SP in (\ref{eq:MSPapp}) has $K+1$ stages. Correspondingly, we develop a scenario tree with $K+1$ levels to encode information available at each stage. The top level ($k=K$) has a single node, which is the root of the tree.  A node at lower levels $k=K-1, \cdots,1$ is the root node of a subtree that starts from that level. On that subtree, the path from the root node to a descendant node at level $k'$ ($0 \le k'<k$) is encoded by a string
\[
\bar{\omega}=\omega^k \cdots\omega^{k'+1}
\]
and $\BD_M^l(\omega^l)$ is a realization of demand $\BD_M^l$ ($k'< l \le k$). Since $\BD_M^l$ are independent, this specification applies to all subtrees that start from level $k$.  Let $\Omega_{k}^{k'}$ be the collection of these strings ($0 \le k' < k \le K$). Note that  $\Omega_{k}^{k'}$ depends on $M$, but we suppress this to ease the notational burden. For $k'<k$, the probability associated with a path encoded by string $\bar{\omega} \in \Omega_{k}^{k'}$ is
\begin{equation}
\label{eq:proscaling}
\mathbb P(\bar{\omega})= \mathbb P^k_M(\omega^k)\times \cdots \times \mathbb P^{k'+1}_M(\omega^{k'+1}),
\end{equation}
where  $\mathbb P^l_M(\omega^l)$ is the probability that $\BD_M^l=\BD_M^l(\omega^l)$ ($k'<l \le k$). For convenience, we also allow $k'=k$, where $\Omega_k^k$ contains a single element, $\bar{\omega}$, corresponding to an empty string, and
$
\mathbb P(\bar{\omega})=1.
$
The total demand realized on the path $\bar{\omega}=\omega^k\cdots\omega^1$ is then denoted by 
\[
\underline{\BD}_M^k (\bar{\omega})=\BD_M^k(\omega^k)+\cdots+\BD_M^1(\omega^1). 
\]

For any two strings $\omega_1$ and $\omega_2$, write $\omega_1 \sqsubset \omega_2$ if $\omega_1$ is a prefix substring of string $\omega_2$. 
On any subtree that starts from level $k$ ($1 \le k \le K$), let  $\bar{\omega}'$ ($\bar{\omega}' \in \Omega^{k'}_k$) encode a path between its root node and a descendant node at level $k'$ ($0<k' \le k$).
Let $\bar{\omega}''$ ($\bar{\omega}'' \in \Omega_k^{k''}$) encode a path between the root node and a descendant node at a lower level $k''$ ($0 \le k'' <k$).
Then the former path is a segment of the latter one if and only if $\bar{\omega}' \sqsubset \bar{\omega}''$. 

A tree starting from level $k$ ($0 \le k \le K$), as specified above, and associated with data $(\y^{k+1},\cdots,\y^K,\x)$ ($\x \ge \mathbf 0$), allows us to formulate $\varphi_M^k(\y^{k+1},\cdots,\y^K,\x)$ as the following LP (recall that $\mathbb P(\bar{\omega})=1$ for $\bar{\omega} \in \Omega_k^k$ as the set contains a single element):
\begin{equation}
\label{eq:finitLPobj+}
\varphi_M^k(\y^{k+1},\cdots,\y^K,\x)=\min_{\y^k,\cdots,\y^1, \mathbf z} \left \{\sum_{k'=1}^{k}  \sum_{\bar{\omega} \in \Omega_k^{k'}}  \mathbb P (\bar{\omega}) \mathbf h^{k'} \cdot \y^{k'}(\bar{\omega})
 -\sum_{\bar{\omega} \in \Omega_k^0} \mathbb P(\bar{\omega}) \mathbf c \cdot \z(\bar{\omega}) \right \}
\end{equation}
subject to: 
\begin{eqnarray}
\label{eq:mapping1+}
\z(\bar{\omega}) ~~~~~~~~~~~~~~~~~~~~~~~~~~~~~~&\le&~~~ \x+\underline{\BD}^k_M (\bar{\omega}),~~~\bar{\omega} \in \Omega_k^0,
\\
\label{eq:mapping2+}
A^{k'} \mathbf z(\bar{\omega}) ~~~- \y^{k'}(\bar{\omega}')~~~~~~~~~~~~~~&\le&~~~ \mathbf 0,~~~~~~~~~~~~~~~\bar{\omega} \in \Omega_k^0,~\bar{\omega}' \in \Omega_k^{k'},~\bar{\omega}' \sqsubset \bar{\omega}, ~1 \le k' <k,
\\
\label{eq:mapping3+}
A^k \mathbf z(\bar{\omega}) ~~~~~~~~~~~~~~~~~~~~~- \y^k~~~ &\le&~~~ \mathbf 0,~~~~~~~~~~~~~~~\bar{\omega} \in \Omega_k^0,
\\
\label{eq:mapping4+}
A^{k'} \mathbf z(\bar{\omega}) ~~~~~~~~~~~~~~~~~~~~~~~~~~~~~~ &\le&~~~ \y^{k'},~~~~~~~~~~~~~\bar{\omega} \in \Omega_k^0, ~k < k' \le K.
\end{eqnarray}
The above is the same formulation as (\ref{eq:finitLPobj})-(\ref{eq:mapping4}) given in Section \ref{subsec:finitedef}.  When $k=K$, (\ref{eq:finitLPobj+})-(\ref{eq:mapping4+}) is the same problem defined in (\ref{eq:MSPapp}). When $k<K$, it defines subproblems of (\ref{eq:MSPapp}) at stage $k$ ($0 \le k <K$), with $\x$ determined by possible realizations of $\BD_M^K+\cdots+\BD_M^{k+1}$ and $(\y^{k+1},\cdots,\y^K)$ given by decisions already taken at stages $K,\cdots,k+1$ ($0 \le k <K$).

For this problem, we can formulate another approximating LP that always has a unique optimal solution by perturbing coefficients in (\ref{eq:finitLPobj+}).

Let $\tilde{\mathbb P}(\bar{\omega})$ be perturbed values of $\mathbb P(\bar{\omega})$ ($\bar{\omega} \in \Omega_k^l, 1 \le l \le k \le K)$, $\tilde{\mathbf h}^l$ be perturbed values of $\mathbf h^l$ $(1 \le l \le k)$, and $\tilde{\mathbf c}$ be the perturbed value of $\mathbf c$.  The goal of the perturbation is to ensure that for all $k=1,\cdots,K$, there does not exist any $\mathbf r^l(\bar{\omega})$ $(\omega \in \Omega_k^l,1\le l \le k)$ and $\mathbf r^0(\bar{\omega})$ $(\bar{\omega} \in \Omega_k^0)$, unless all of them are zeros, such that 
\begin{equation}
\label{eq:rational}
\sum_{l=1}^k \sum_{\bar{\omega}  \in \Omega_k^l} \tilde{\mathbb P}(\bar{\omega}) \tilde{\mathbf h}^l \cdot \mathbf r^l(\bar{\omega}) +\sum_{\bar{\omega} \in \Omega_k^0} \tilde{\mathbb P}(\bar{\omega}) \tilde{\mathbf c} \cdot \mathbf r^0(\bar{\omega})=0,
\end{equation}
which is easily achievable with a negligible perturbation error. 

For instance, we can perturb values of $\mathbb P_M^k(\omega^k)$ ($1 \le k \le K)$ and every component of $\mathbf h^k$ $(1 \le k \le K)$ and $\mathbf b$ to make each of them a sum of a rational number and a unique irrational number, where the latter number can be the square root of a unique prime number divided by an integer, and then determine values of $\mathbb P(\bar{\omega})$ ($\bar{\omega} \in \Omega_k^l, 1 \le l \le k \le K)$ and $\mathbf c$ accordingly. In each case, we can make the rational part of the perturbed value arbitrarily close to the original value, and the irrational part arbitrarily small by dividing it with a large enough integer. Moreover, by Theorem \ref{thm:solbound}, the finiteness of $\mathbf D_M^k$ $(1 \le k \le M)$ implies the finiteness of the optimal solution. Therefore we can always keep the difference of the optimal objective values between the original problem in (\ref{eq:finitLPobj+})-(\ref{eq:mapping4+}) and its perturbed version sufficiently small to keep the loss of optimality by perturbation negligible.

At each stage $k$ $(1 \le k \le K)$, the optimal solution of the perturbed LP must be unique. If not, then the problem will have two extreme-point solutions, which must be rational-valued since $A^k$ and $\underline{\BD}_M^k(\bar{\omega})$ ($1 \le k \le K$) in (\ref{eq:mapping1+})-(\ref{eq:mapping4+}) are integers, and by induction from $k=K$ downwards, $\mathbf y^{k'}$ $(k<k' \le K)$ on the RHS of (\ref{eq:mapping4+}) are rational vectors. The two solutions yield the same objective value, which means that (\ref{eq:rational}) can be satisfied with rational vectors $\mathbf r^l(\bar{\omega})$ and  $\mathbf r^0(\bar{\omega})$ if we let them be respectively the differences in $\y^l(\bar{\omega})$ $(\bar{\omega} \in \Omega_k^l, 1 \le l \le k)$ and $\mathbf z(\bar{\omega})$  $(\bar{\omega} \in \Omega_k^0)$ between the two solutions. This contradicts with the design of the perturbation scheme. 

%For systems with larger lead times, a large $M$ is needed to keep the (discrete) approximation sufficiently accurate. This will lead to more elements in $\Omega_k^l$ ($0 \le l \le K$, $1 \le k \le K$), and thus require more unique values of $\tilde{\mathbb P}(\bar{\omega})$. This requirement can be accommodated as the set of prime numbers is an infinite set. Moreover, we can use an arbitrage large value for the aforementioned integer to divide these unique values to keep the perturbation error negligible. 

%In Sections \ref{subsec:finitedef} and \ref{subsec:FCBM}, we slightly abuse the notation by continuing to refer to the perturbed version of (\ref{eq:finitLPobj})-(\ref{eq:mapping4}), which has a unique optimal solution and a negligible difference in the optimal objective value  from the original problem, as  $\varphi^k_M(\mathbf y^{k+1},\cdots,\mathbf y^K,\mathbf x)$ $(1 \le k \le K)$. 

\bigskip

\noindent
{\bf Acknowledgment:}  {
\begin{itemize}
\item
This study is based upon work supported by the National Science Foundation under Grant No. CMMI-1363314.
\item
The authors are grateful to the associate editor and anonymous referees for their constructive comments. 
\end{itemize}

\bibliography{ATO_MWW}
\bibliographystyle{siam}

\end{document}